\documentclass[11pt]{amsart}
\usepackage[T1]{fontenc}
\usepackage{lmodern}
\usepackage{amsmath,amssymb,amsthm,mathtools}
\usepackage[margin=1.03in]{geometry}
\usepackage{microtype,booktabs,longtable,array,enumitem}
\usepackage{xcolor}
\definecolor{paperlink}{RGB}{0,0,114}
\usepackage{hyperref}
\hypersetup{colorlinks=true,linkcolor=paperlink,citecolor=paperlink,urlcolor=paperlink,pdfstartview=Fit,pdftitle={Periods, isomonodromy, and modularity of ALG/ALG* gravitational instantons}}
\setlist{itemsep=.2em,topsep=.35em}
\numberwithin{equation}{section}
\newtheorem{theorem}{Theorem}[section]
\newtheorem{maintheorem}[theorem]{Theorem}
\newtheorem{proposition}[theorem]{Proposition}
\newtheorem{lemma}[theorem]{Lemma}
\newtheorem{corollary}[theorem]{Corollary}
\theoremstyle{definition}
\theoremstyle{remark}\newtheorem{remark}[theorem]{Remark}
\DeclareMathOperator{\tr}{tr}
\DeclareMathOperator{\diag}{diag}
\DeclareMathOperator{\Res}{Res}
\DeclareMathOperator{\Id}{Id}
\DeclareMathOperator{\Rea}{Re}
\DeclareMathOperator{\Ima}{Im}

\DeclareMathOperator{\End}{End}
\DeclareMathOperator{\Pic}{Pic}

\DeclareMathOperator{\SU}{SU}
\DeclareMathOperator{\SL}{SL}
\newcommand{\C}{\mathbb C}\newcommand{\R}{\mathbb R}
\newcommand{\Z}{\mathbb Z}\newcommand{\PP}{\mathbb P}
\newcommand{\cO}{\mathcal O}\newcommand{\cM}{\mathcal M}
\newcommand{\cD}{\mathcal D}\newcommand{\dd}{\mathrm d}
\newcommand{\ddc}{\dd\dd_J^c}
\newcommand{\Op}{O'}\newcommand{\sfm}{\mathrm{sf}}
\title[Periods, isomonodromy, and modularity]{Periods, isomonodromy, and modularity of ALG/ALG$^*$ gravitational instantons}
\author{Shiyu Cao}
\address{Chern Institute of Mathematics and LPMC, Nankai University,
Tianjin 300071, China}
\email{shiyucao@126.com}
\date{29 November 2026}
\subjclass[2020]{Primary 53C26; Secondary 14D21, 34M56, 58J05}
\keywords{Wild nonabelian Hodge theory, gravitational instantons, Hitchin metric, isomonodromic deformation, Torelli theorem}
\begin{document}
\begin{abstract}
We prove that every smooth ALG or ALG$^*$ gravitational instanton of type $A_0,A_1,A_2,D_0,D_1,D_2,D_3$ is, up to homothety and hyperk\"ahler rotation, isometric to a rank-two wild Hitchin moduli space on $\PP^1$ with its $L^2$ metric. We compute the compact hyperk\"ahler periods and asymptotic invariants of these moduli spaces and give explicit inverses on the admissible parameter domains. For the three ALG types, global isomonodromic transport realizes the two real parameters left undetermined by the compact periods. A Hamiltonian variation formula for the $L^2$ metric and uniform estimates for its difference from the semiflat metric identify the data required by the Torelli and uniqueness theorems for ALG and ALG$^*$ gravitational instantons.
\end{abstract}
\maketitle
\enlargethispage{4pt}
\vspace{-3pt}
\tableofcontents
\section{Introduction}
\label{sec:introduction}

\subsection{Nonabelian Hodge theory and metric asymptotics}

Hitchin introduced the self-duality equations on a Riemann surface and endowed their moduli spaces of solutions with a hyperk\"ahler structure~\cite{Hitchin87a}. The correspondence between Higgs bundles and flat connections was developed through the harmonic-map and gauge-theoretic methods of Donaldson, Corlette, and Simpson~\cite{Donaldson87,Corlette88,Simpson88,Simpson92}. In complex structure $I$, the characteristic polynomial of the Higgs field defines an algebraic integrable system~\cite{Hitchin87b}, whose smooth generic fiber parametrizes line bundles on the spectral curve~\cite{BNR89}. The holomorphic deformation complex carries the complex symplectic form. The harmonic metric induces the $L^2$ pairing on harmonic representatives of tangent vectors.

For punctured curves, the filtration at a marked point records the growth of the harmonic metric. In the tame case, Simpson studied harmonic bundles on noncompact curves and Konno constructed moduli spaces of parabolic Higgs bundles~\cite{Simpson90,Konno93}. At an irregular singularity one also prescribes the formal irregular type. Sabbah proved the existence of harmonic metrics for irregular connections~\cite{Sabbah99}. For unramified irregular data with diagonal polar parts, Biquard--Boalch established the wild nonabelian Hodge correspondence and constructed $L^2$ hyperk\"ahler metrics by reduction in weighted Sobolev spaces. These metrics are complete when the full moduli space has no reducible points~\cite[Theorems~0.1 and~5.4]{BB04}. Mochizuki proved the existence and uniqueness results for harmonic metrics that we use here, including the case of good filtered bundles with ramified irregular data~\cite{Moch11,Moch21}. Eper--Szab\'o later used invariant weighted Sobolev spaces on a ramified cover to construct the hyperk\"ahler metrics for rank-three Painlev\'e realizations~\cite[Theorem~2.1 and Remark~2.2]{EperSzabo25}. Here the construction and completeness of the gauge quotient require arguments in addition to the existence of an adapted harmonic metric.

The semiflat metric is determined by the special K\"ahler geometry of the regular locus of the Hitchin base and the dual tori~\cite{Freed99}. Gaiotto--Moore--Neitzke proposed a description of the difference between the Hitchin metric and this semiflat metric in terms of spectral data and wall crossing~\cite{GMN13}. Analytically, Mazzeo--Swoboda--Weiss--Witt constructed rank-two limiting configurations and their desingularizations, and developed the asymptotic expansion of the Hitchin metric~\cite{MSWW16,MSWW19}. Exponential decay of the difference between the Hitchin and semiflat metrics was proved on the Hitchin section by Dumas--Neitzke~\cite{DN19} and along generic rays in higher rank by Fredrickson~\cite{Fredrickson20}. Mochizuki determined these asymptotics in arbitrary rank on the locus of smooth spectral curves~\cite{Moch26}. In these analyses, harmonic metric estimates control the $L^2$ norms of harmonic representatives of tangent vectors to the moduli space.

In the presence of singularities, the local models and admissible parameter ranges depend on the pole orders and ramification. Fredrickson--Mazzeo--Swoboda--Weiss proved asymptotic results for the Hitchin metric on the moduli space of rank-two parabolic Higgs bundles~\cite{FMSW22}. Chen--Li studied rank-two irregular Higgs bundles on the projective line and determined the models at infinity of the ALG and ALG$^*$ moduli spaces~\cite{CL24}. Mochizuki--Szab\'o and Sagman--Smillie extended earlier results on harmonic metric asymptotics and estimates for high-energy harmonic maps to settings with weaker regularity assumptions~\cite{MochSz25,SS26}. For a family of strongly parabolic rank-three Higgs bundles on the three-punctured sphere with weights parametrized by $0<\delta<1/3$, exponential decay of the difference between the Hitchin and semiflat metrics was established in~\cite[Theorems~1.1--1.3]{Cao26}. These moduli spaces carry complete ALG-$E_6$ metrics with $\beta=1/3$ and $\tau=e^{2\pi i/3}$. 

Related degenerations arise when parabolic weights tend to zero. Heller--Heller--Traizet used loop groups to construct twistor lines in a shrinking bounded-energy region for small equal weights on the four-punctured sphere. After rescaling, they identified the limiting metric as the Eguchi--Hanson metric modulo a $\Z_2\times\Z_2$ action~\cite{HHT25}. Heller--Heller--Meneses proved real-analytic convergence of the rescaled Hitchin metrics to the hyperpolygon metric on compact subsets of the corresponding bounded-energy domains for small generic weights, and obtained a convergent higher-order expansion~\cite[Theorems~5.3--5.4]{HHM25}. Fredrickson--Yae proved pointwise convergence of rescaled strongly parabolic Hitchin metrics to the hyperpolygon metric; with four punctures, their construction gave a degeneration from ALG-$D_4$ to ALE-$D_4$~\cite{FY26}.

Wong interpreted certain quadratic Hitchin Hamiltonians through differences of isomonodromic flows~\cite{Wong12}. Boalch constructed symplectic isomonodromic connections for wild character varieties~\cite{Boalch14}, and Korman studied the corresponding curvature terms in parabolic hyperk\"ahler quotients~\cite{Korman17}. Over the universal family of compact Riemann surfaces, Hitchin related the variation of harmonic-map energy to the variation of the moduli-space metric~\cite{Hitchin26}. The isomonodromic geometry of universal Hitchin moduli spaces was studied in~\cite{AGGT26,CTW25}.

\subsection{Gravitational instantons and Boalch's conjecture}

Kronheimer constructed ALE spaces as moduli spaces and classified their metrics by periods~\cite{Kronheimer89a,Kronheimer89b}. Hein used rational elliptic surfaces to construct gravitational instantons of quadratic volume growth~\cite{Hein12}. Chen--Chen developed asymptotic models, improved decay estimates, and compactification under curvature-decay hypotheses~\cite{CC21a,CC19,CC21b}. Sun--Zhang established the asymptotic structure of complete hyperk\"ahler four-manifolds with $L^2$ curvature~\cite{SZ24}; Li--Sun gave a recent account of these developments~\cite{LiSun26}.

We use Chen--Viaclovsky's construction and classification results for ALG and ALG$^*$ gravitational instantons of quadratic volume growth. In the ALG case with cone-angle parameter $\beta>1/2$, their description includes two additional real parameters~\cite{CV24}. Chen--Viaclovsky--Zhang proved Torelli theorems for ALG and ALG$^*$ spaces with order-$2$ decay~\cite{CVZ24}. For ALG$^*$ ends they used Hodge theory to improve the decay estimate while preserving the asymptotic constants~\cite{CVZH23}. Lee--Lin determined the corresponding integral period domains and proved surjectivity of the period maps~\cite[Theorem~3.13]{LL23}.

For the tame type $D_4$, Fredrickson--Mazzeo--Swoboda--Weiss computed the compact periods and integral topology in Part I of their Torelli program~\cite[Theorem~4.1]{FMSW26}. Their Theorem~9.2 showed that the period map for the standard parabolic Hitchin moduli spaces has open dense image in the period domain. They proposed an enlarged definition involving partial flags in Conjecture~1.5 and deferred the analysis of the ends to a companion paper. In the tame parabolic setting, Dimakis--Fredrickson--Yae proved affine variation of hyperk\"ahler periods under parallel transport on the generic parameter locus, under their coprimality and nonemptiness hypotheses~\cite[Theorem~5.5]{DFY26}.

Boalch conjectured that every ALG or ALG$^*$ gravitational instanton is isometric to a Hitchin moduli space with its hyperk\"ahler metric~\cite[Conjecture~1.2]{FMSW26}. His list of H3 surfaces associated four-dimensional meromorphic Hitchin systems to affine Dynkin types~\cite[Section~4.1]{Boalch18}. We prove the conjecture for the seven irregular rank-two types $A_0,A_1,A_2,D_0,D_1,D_2,D_3$. Here modularity means realization as a Hitchin moduli space with its $L^2$ metric, up to homothety and hyperk\"ahler rotation. The companion paper~\cite{CaoExceptional26} treats the exceptional ALG types $E_6,E_7,E_8$ using parabolic Higgs bundles on the three-punctured sphere.

Chen--Li included these seven types in their list of four-dimensional wild Hitchin moduli spaces, with nine pole configurations because $A_1,D_2$ each have two realizations~\cite[Tables~1--2]{CL24}. We choose one spectral curve equation for each type. The compact periods determine the weights and residues; a remaining leading coefficient determines the ALG$^*$ logarithmic constant. For $A_0,A_1,A_2$, filtered isomonodromic transport, defined for all values of $a$, realizes the two real parameters that the compact periods do not determine.

For one of these four-dimensional moduli spaces $M$, let $H_2=H_2(M;\Z)$ and let $F$ be the class of a smooth elliptic Hitchin fiber with its complex orientation. We normalize the ordered hyperk\"ahler forms by $\int_F(\omega_I,\omega_J,\omega_K)=(4\pi^2,0,0)$. In Table~\ref{tab:types}, $u$ is the Hitchin base coordinate and the eigenvalue one-forms satisfy $\lambda^2=q(z)\dd z^2$.

{\small
\setlength{\tabcolsep}{4pt}
\renewcommand{\arraystretch}{1.2}
\begin{longtable}{@{}p{.07\textwidth}p{.88\textwidth}@{}}
\caption{The seven spectral families, compact lattices, and end data.}\label{tab:types}\\
\toprule
Type & Spectral coefficient $q(z)$; poles and geometric data\\
\midrule
\endfirsthead
\multicolumn{2}{l}{\small Table~\thetable\ continued.}\\
\toprule
Type & Spectral coefficient $q(z)$; poles and geometric data\\
\midrule
\endhead
\bottomrule
\endfoot
$A_0$ & \refstepcounter{equation}\label{eq:familyA0}$q(z)=z^3+az+u$.\hfill(\theequation)\\*
& Pole: ramified wild at $\infty$ (order $4$). ALG, $\beta=5/6$; $H_2/\Z F=0$; $c_{\mathrm{CV}}\in\C$.\\
\midrule
$A_1$ & \refstepcounter{equation}\label{eq:spectral}$q(z)=z^4+az^2+bz+u$.\hfill(\theequation)\\*
& Pole: unramified wild at $\infty$ (order $4$). ALG, $\beta=3/4$; $H_2/\Z F=\langle-2\rangle$; $c_{\mathrm{CV}}\in\C$.\\
\midrule
$A_2$ & \refstepcounter{equation}\label{eq:familyA2}$q(z)=z^{-2}\bigl(z^4+az^3+(a^2/4+2\mu)z^2+uz+c^2\bigr)$.\hfill(\theequation)\\*
& Poles: tame at $0$, unramified wild at $\infty$ (order $3$). ALG, $\beta=2/3$; $H_2/\Z F=-A_2$; $c_{\mathrm{CV}}\in\C$.\\
\midrule
$D_0$ & \refstepcounter{equation}\label{eq:familyD0}$q(z)=p_0z^{-3}+uz^{-2}+p_\infty z^{-1}$, $p_0p_\infty\ne0$.\hfill(\theequation)\\*
& Poles: ramified wild at $0,\infty$ (order $2$). ALG$^*$, $\nu=4$; $H_2/\Z F=0$; $\kappa_0\in\R$.\\
\midrule
$D_1$ & \refstepcounter{equation}\label{eq:familyD1}$q(z)=p^2z^{-4}+2pmz^{-3}+uz^{-2}+bz^{-1}$, $pb\ne0$.\hfill(\theequation)\\*
& Poles: unramified wild at $0$, ramified wild at $\infty$ (both order $2$). ALG$^*$, $\nu=3$; $H_2/\Z F=\langle-4\rangle$; $\kappa_0\in\R$.\\
\midrule
$D_2$ & \refstepcounter{equation}\label{eq:familyD2}$\begin{aligned}[t]
q(z)={}&p_0^2z^{-4}+2p_0m_0z^{-3}+uz^{-2}\\
 &+2p_\infty m_\infty z^{-1}+p_\infty^2,\qquad p_0p_\infty\ne0.
\end{aligned}$\hfill(\theequation)\\*
& Poles: unramified wild at $0,\infty$ (order $2$). ALG$^*$, $\nu=2$; $H_2/\Z F=\langle-2\rangle\oplus\langle-2\rangle$; $\kappa_0\in\R$.\\
\midrule
$D_3$ & \refstepcounter{equation}\label{eq:familyD3}$\begin{aligned}[t]
q(z)={}&p^2+2pmz^{-1}+\frac{u}{z(z-1)}\\
 &+c_0^2z^{-2}+c_1^2(z-1)^{-2},\qquad p\ne0.
\end{aligned}$\hfill(\theequation)\\*
& Poles: tame at $0,1$, unramified wild at $\infty$ (order $2$). ALG$^*$, $\nu=1$; $H_2/\Z F=-A_3$; $\kappa_0\in\R$.\\
\end{longtable}
}
{\footnotesize The pole orders refer to the meromorphic Higgs field in the chosen holomorphic lattices; a tame pole has order one. The coefficients other than $u$ specify the principal parts. Among these family parameters, $a$ denotes the ALG isomonodromic parameter, while $p_0,p_\infty$ denote the $D_0$ leading coefficients. The notation $-A_k$ denotes the negative root lattice and $\langle-d\rangle$ a rank-one lattice generated by a class of square $-d$. The ALG cone angle is $2\pi\beta$; $c_{\mathrm{CV}}$ is the Chen--Viaclovsky deformation parameter. For ALG$^*$, $V(r)=\kappa_0+(\nu/\pi)\log r$ with scale $L=1$. Their precise normalizations appear in Table~\ref{tab:end-parameters} and~\eqref{eq:CVpotentialmatch}. The fiber $F$ has periods $(4\pi^2,0,0)$ throughout.\par}

An admissible compact period triple has these fiber periods and a nonzero period vector on every integral class of square $-2$. The parameter domains in Table~\ref{tab:compact-periods} include the permitted equal and integer weights and retain the parabolic flags when a tame residue vanishes; the additional end data are listed in Table~\ref{tab:types}.

\begin{maintheorem}
\label{thm:realization}
Up to homothety and hyperk\"ahler rotation, every smooth gravitational instanton of one of the seven types in Table~\ref{tab:types} is isometric to a moduli space in the corresponding filtered Hitchin family equipped with its $L^2$ metric. In the normalized integral marking, every admissible compact period triple and every additional end datum listed there is realized. Explicit inverses of the compact period maps are given in Table~\ref{tab:compact-periods}; the parameters at infinity are given in Table~\ref{tab:end-parameters} and~\eqref{eq:CV-inverse}.

After suitable admissible changes of integral marking and of the normalized coordinates at infinity, the isometries preserve the ordered hyperk\"ahler triples.
\end{maintheorem}

\subsection{Hamiltonian variation and metric comparison}

The variation formula holds on an arbitrary compact marked Riemann surface. We use it on the seven moduli surfaces of Table~\ref{tab:types} and on the eight-dimensional example in Appendix~\ref{app:variation}.

Let $X$ be such a marked surface, with the following prescribed singular data. At an unramified wild pole, the leading coefficient is regular semisimple. In rank two, we also allow tame poles with residues $(c,-c)$ on the graded quotients and weights $(\alpha,1-\alpha)$ satisfying
\begin{equation}
 0\leq\alpha\leq\tfrac12,\qquad c\ne0\quad\text{or}\quad0<\alpha<\tfrac12,
 \label{eq:tame-domain}
\end{equation}
and ramified poles for which the pulled-back Higgs field has regular semisimple leading coefficient in its formal normal form on a specified finite cover. When $c=0$, the nilpotent strongly parabolic residue is allowed to vanish. We fix the filtered formal data, a compatible flat parabolic determinant line bundle, and weights for which the filtered bundle has parabolic degree zero. On the smooth moduli space of stable filtered Higgs bundles we take the trace-free deformation complex. Shifting a real representative of a parabolic weight by an integer changes the holomorphic lattice by an elementary modification and induces the corresponding change in the determinant extension. For the de Rham complex structure $J$, we use $\dd_J^c=-J^*\dd$, so $\dd\dd_J^c=2i\partial_J\bar\partial_J$ on functions.

\begin{maintheorem}
\label{thm:variation}
For unramified wild data in arbitrary rank or for rank-two bundles with the pole types specified above, solutions to Hitchin's equations admit smoothly varying local gauge slices, and the $L^2$ horizontal projections vary smoothly as the complex structure $j_s$ on $X$ varies in a smooth finite-dimensional family whose members agree with $j_0$ near every marked point. The maps between the moduli spaces that keep the filtered flat connection fixed preserve the de Rham complex structure $J$ and its holomorphic symplectic form.

Write the filtered flat connection as $\cD=A+\Psi$, where $A$ is unitary for the harmonic metric $h$ and $\Psi$ is a one-form with values in the $h$-self-adjoint endomorphisms of $E$. Along any smooth one-parameter path $j_s$ in this family, set $\Phi_s=\Psi^{1,0}_{j_s}$ and define the Beltrami differential $\mu_s$ by the following convention: in a $j_s$-holomorphic coordinate $z$, the holomorphic cotangent line for $j_{s+\varepsilon}$ is generated locally by $\dd z+\varepsilon\mu_s\dd\bar z+O(\varepsilon^2)$. Then the $L^2$ K\"ahler form satisfies
\begin{equation}
 \frac{\dd}{\dd s}\omega_{J,s}
 =\dd\dd_J^c\left(2\Ima\int_X\tr(\Phi_s^2)\mu_s\right).
 \label{eq:main-variation}
\end{equation}
The integrand has compact support away from the marked points. The formula holds throughout the smooth locus of stable filtered Higgs bundles, even when the spectral curve is singular.
\end{maintheorem}

For the $a$-variation in the three ALG families, evaluation of the residues gives
\begin{equation*}
 \dot\omega_J=-K_j\ddc\Rea(\dot a\,u),\qquad
 (K_0,K_1,K_2)=(8\pi,4\pi,4\pi).
\end{equation*}
Section~\ref{subsec:quarticvariation} computes these constants in the spectral normalization $\Phi=i\Phi_{\mathrm{sp}}$.

For each family in Table~\ref{tab:types}, $u:\cM\to\C$ denotes the Hitchin map on its moduli space $\cM$, and $g_{\sfm}$ denotes the semiflat metric normalized so that both torus angles have period $2\pi$. The admissible parameter domains are described in Section~\ref{sec:integral}.

\begin{maintheorem}
\label{thm:comparison}
In each of the seven families, at every point on a Hitchin fiber with sufficiently large $|u|$, the $L^2$ and semiflat metrics satisfy the quadratic-form bounds
\begin{equation}
 (1+C\varepsilon(u))^{-1}g_{\sfm}\leq g_{L^2}\leq(1+C\varepsilon(u))g_{\sfm},
 \label{eq:main-comparison}
\end{equation}
with the constants and the threshold for $|u|$ chosen uniformly on compact subsets of the space of admissible weights and prescribed singular data. One can take
\begin{equation*}
 \varepsilon(u)=
 \begin{cases}
 e^{-c|u|^{5/6}},&A_0,\\
 e^{-c|u|^{3/4}},&A_1,\\
 |u|^{-1/10},&A_2,D_3,\\
 e^{-c|u|^\sigma}\text{ for some }\sigma>0,&D_0,D_1,D_2.
 \end{cases}
\end{equation*}
The corresponding estimates hold for covariant derivatives of the metric difference on the moduli space at every order, with constants depending on the order and possibly a smaller positive polynomial decay exponent or exponential rate, as detailed in Section~\ref{subsec:spatial-comparison}. The comparison includes all equal and integer weights in the admissible filtered domains, and all strongly parabolic tame data allowed by~\eqref{eq:tame-domain}.
\end{maintheorem}

Chen--Li prove the comparison along generic curves and obtain estimates uniform on the four-dimensional ends for their fixed pole data, with equal weights $1/2$ at weakly parabolic poles~\cite[Theorem~1.1, Assumption~6.8 and Section~6.4]{CL24}. In Theorem~\ref{thm:comparison} the estimates for each spatial derivative order are uniform on compact sets of admissible pole data, including equal and integer weights at unramified wild poles and the full tame range~\eqref{eq:tame-domain}.

In an ALG family $A_j$, the weights and all coefficients of the prescribed principal parts other than $a$ remain constant. We denote the moduli space by $\cM_a$, its $L^2$ metric by $g_a$, and its Hitchin map by $u_a:\cM_a\to\C$. All distances, norms and covariant derivatives are taken with respect to the metric under consideration. We write $\Op(\rho^q)$ for a tensor whose $k$th covariant derivative is $O(\rho^{q-k})$ for every $k\geq0$, where $\rho$ is one plus the distance from a chosen basepoint. In Theorem~\ref{thm:transport} we use $g_0$ on $\cM_0$. We set $(\beta_0,\beta_1,\beta_2)=(5/6,3/4,2/3)$ and use the Hamiltonian constants $K_j$ above. The constants $\ell_j>0$ are the leading coefficients in the normalized spectral period coordinates $U\sim\ell_j u^{\beta_j}$; their values are given below.

\begin{maintheorem}
\label{thm:transport}
For every admissible choice of compact periods in $A_j$, $j=0,1,2$, the filtered isomonodromic transport $T_a:\cM_0\to\cM_a$ constructed in Section~\ref{sec:finite} exists for all $a\in\C$ and preserves $J,\omega_I,\omega_K$. The metric at $a=0$ has an ALG end with order-$2$ decay. There is a smooth real function $\phi_a$ on $\cM_0$ such that
\begin{equation*}
 T_a^*\omega_{J,a}=\omega_{J,0}+i\partial_J\bar\partial_J\phi_a\,.
\end{equation*}
Its leading term is $-2K_j\Rea(a u_0)$, with remainder $\Op(\rho^{2/5})$ in type $A_0$, $\Op(\rho^{2/3})$ in type $A_1$, and $\Op(\rho^{3/2-\delta})$ for some $\delta>0$ in type $A_2$, uniformly for $a$ in compact subsets of $\C$.

With respect to the continuously varying asymptotic coordinates chosen in Section~\ref{sec:periodend}, the rotated triples $(\omega_I,\omega_K,-\omega_J)$ form the two-parameter family of~\cite[Proposition~5.3]{CV24}, with
\begin{equation*}
 c_{\mathrm{CV}}=i\gamma_j\bar a,\qquad\gamma_j=K_j\ell_j^{-1/\beta_j}.
\end{equation*}
\end{maintheorem}

The positive leading coefficients in the asymptotic expansions of the base coordinates defined by the spectral periods are
\begin{equation}
 \ell_0=\frac{4\sqrt{3\sqrt3/2}}3B(1/3,3/2),\qquad
 \ell_1=\frac83\mathbf K(1/\sqrt2),\qquad
 \ell_2=\frac{4\sqrt{\sqrt3/2}}3B(1/6,3/2),
 \label{eq:ell-constants}
\end{equation}
where $B$ is the beta function and $\mathbf K(k)$ is the complete elliptic integral of the first kind with modulus $k$~\cite[Sections~5.12 and~19.2]{DLMF}. The choices of unit complex factors and integral cycles on the spectral curve in Section~\ref{sec:periodend} determine the continuously varying choice of coordinates at infinity used in Theorem~\ref{thm:transport}.

\subsection{Proof strategy and organization}

Sections~\ref{sec:families}--\ref{sec:variation} express the variation of the $L^2$ K\"ahler form under isomonodromic transport as the complex Hessian of a Hamiltonian. In a local gauge slice the vertical contribution cancels the variation along the complex gauge orbit. Residues determine the Hamiltonians for the three ALG families. On closed parameter surfaces, Chern--Weil transgression and changes of the Liouville form give the real and complex compact periods in terms of the discriminant of the universal bundle and the degrees and residues of the graded lines at the poles (Section~\ref{sec:compact}).

The estimates in Section~\ref{sec:comparison} apply uniformly on the distant Hitchin fibers. In the tame local model the difference between the projected harmonic and decoupled connections extends to a holomorphic one-form on the compactified spectral normalization $\PP^1$ and hence vanishes. On an annulus whose monodromy exchanges the sheets, the angular Poincar\'e inequality gives the coercivity needed to compare the glued and harmonic metrics. Sections~\ref{sec:integral}--\ref{sec:periodend} compute the integral marking and the models at infinity. For the ALG types, isomonodromic transport supplies the two real parameters invisible to the compact periods. Sections~\ref{sec:finite}--\ref{sec:realization} recover them from the leading potential and apply Torelli in a common end chart.

\subsection*{Acknowledgment}
The author is grateful to Professors Sebastian Heller and Szil\'ard Szab\'o for their interest in an earlier version of this work and for their helpful comments and suggestions.

\par\smallskip\noindent\textit{Declaration of AI use.}\enspace
The research program and mathematical ideas are the author's. The author used ChatGPT (OpenAI) for literature searches, language editing and calculations, and to assist in developing parts of the metric comparison proof and checking arguments throughout the paper. The author takes responsibility for the manuscript.

\section{Filtered Higgs bundles and harmonic families}
\label{sec:families}

We use weighted local slices for regular semisimple wild poles, for tame poles with a parabolic flag, and, after passage to a local cover, for ramified wild poles. In the last case stability on the original curve determines the kernel on invariant sections.

\subsection{Metric normalization and notation}
In the seven rank-two families in Table~\ref{tab:types}, the Higgs fields are trace-free, the determinant line bundle with its parabolic structure is fixed, and the eigenvalue one-forms are denoted by $\pm\lambda$. In the general variation theorem the determinant is a prescribed flat parabolic line bundle, and we work with the trace-free deformation complex. We regard each parabolic weight as a real representative of its class modulo $\Z$, determined by the chosen holomorphic lattice in the eigenline bundle. At wild poles these representatives need not lie in $[0,1)$. Lowering a weight by an integer changes the holomorphic extension by an elementary modification. The filtered bundle and its harmonic metric remain unchanged under this operation. In these spectral equations, $u$ is the coordinate on the Hitchin base and the other coefficients are determined by the prescribed principal parts of the eigenvalue one-forms.

The complex structure $I$ is the one in which the Higgs field is holomorphic, and the hyperk\"ahler forms satisfy
\begin{equation}
 \Omega_I=\omega_J+i\omega_K,\qquad
 \int_F(\omega_I,\omega_J,\omega_K)=(4\pi^2,0,0).
 \label{eq:fiberperiod}
\end{equation}
Here $F$ is a smooth elliptic Hitchin fiber with its complex orientation. For a compact two-cycle $\Gamma$ in the four-dimensional moduli space, the normalized periods are
\begin{equation}
 r_\Gamma=\frac1{4\pi^2}\int_\Gamma\omega_I,
 \qquad Q_\Gamma=\frac1{4\pi^2}\int_\Gamma\Omega_I.
 \label{eq:compact-normalization}
\end{equation}
Fix the formal data and a contractible sector $U_0$ of regular values of the Hitchin map. Choose a parallel integral basis $(A,B)$ of the homology of the compact spectral curve, with $A\cdot B=1$, and choose parallel lifts to the curve with the poles of $\lambda$ removed. For either lifted cycle put $Z_\gamma=\oint_\gamma\lambda$. Changing a lift adds $2\pi i$ times an integral linear combination of the residues, which is constant in $u$ for fixed formal data. We set $\tau(u)=(\partial_uZ_B)/(\partial_uZ_A)$, with the weights and the prescribed singular parts of the eigenvalues held fixed. The derivative $\partial_u\lambda$ is a nonzero holomorphic differential, so $\Ima\tau(u)>0$. 

On the smooth spectral curve $\Sigma_u$, denote by $S_u$ the union of the ramification points and the preimages of the marked poles, counting coincident points once. The unitary connection on the decoupled spectral eigenline is flat on $\Sigma_u\setminus S_u$. Its peripheral holonomies are fixed by the weights and the local ramification factors. Choose at a base point the character of one such eigenline and transport this reference character in the local system of the punctured spectral curves. The ratio of the character of any spectral eigenline to the reference character is trivial on the meridians and defines a character of $H_1(\Sigma_u;\Z)$. We fix the signs of the angles by requiring the holonomy ratios on $A,B$ to be $\exp(i\theta_e),\exp(i\theta_m)$, respectively. For the connection convention $\dd+A$, changing these angles by the periods of a closed real form $\beta$ changes the connection by $-i\beta$. Horizontality holds the character fixed under this transport. Two parallel choices of reference change the angles by constants on $U_0$.

To see the corresponding local lift, extend the reference character to the total punctured spectral family over the contractible sector. Its flat line bundle carries a holomorphic structure. Extend it at the marked sections with the prescribed logarithmic lifts, using the ramification twist when passing between the eigenline subbundle and the spectral quotient line. On the exterior, a horizontal variation keeps the flat transitions constant. On a simple ramification disk the same assertion holds after removal of the half-ramification factor. If a tame residue $c$ vanishes, choose the compact cycles outside a disk containing the three points $0,c,-c$ in the local spectral coordinate. The meridian around the disk has the same character through the collision: the eigenline residues $1/2,\alpha,1-\alpha$ add to the limiting residue $3/2$, as verified in Proposition~\ref{prop:globalmodel}. On the common exterior of this disk the compact cycles can be transported through the collision. A unitary character extends over the punctured disk with prescribed meridians exactly when their product equals its boundary character, since these meridians generate its first homology with this boundary relation. This identifies the exterior character coordinates across $c=0$. These descriptions give the horizontal lifts used to compare the tangent pairings in Section~\ref{sec:comparison}. Ratios of extensions with fixed local data have degree zero and the norm map takes values in $\Pic^0(\PP^1)$, which is trivial. Thus the determinant places no further restriction on the angle lattice.

\paragraph{Conventions.}
For endomorphism-valued forms, a superscript $0$ denotes the trace-free part. We denote the spectral Higgs field by $\Phi_{\mathrm{sp}}$ and write
\begin{equation}
 \cD=A+i\Phi_{\mathrm{sp}}-i\Phi_{\mathrm{sp}}^{*h}=A+\Psi,
 \qquad\Phi=\Psi^{1,0}=i\Phi_{\mathrm{sp}},
 \qquad\tr(\Phi_{\mathrm{sp}}^2)=2\lambda^2=-\tr(\Phi^2).
 \label{eq:phase}
\end{equation}
The de Rham complex structure is $J(a,\psi)=(i\psi,ia)$. We use $\dd_I^c=-I^*\dd$ and
\begin{equation}
 \dd_J^c=-J^*\dd,\qquad \ddc=2i\partial_J\bar\partial_J.
 \label{eq:dc}
\end{equation}
In Section~\ref{sec:variation}, $\Phi$ is used for the connection decomposition. Once the residue formula is converted by~\eqref{eq:phase}, all spectral and period calculations use $\Phi_{\mathrm{sp}}$ and $\lambda$.

The lattice scale in the normalized ALG model is $L=2\pi/\sqrt{\Ima\tau_0}$, where $\tau_0$ is the limiting modular parameter in the chosen integral basis; for ALG$^*$ it is $L=1$ in the notation of~\cite[Section~2]{CVZ24}. The leading coefficients in the expansions of the spectral periods as functions of the base coordinate will be denoted by $\ell_j$.

The estimates are uniform when the admissible weights and singular parts of the eigenvalue one-forms range over a compact subset $K$, with constants allowed to depend on $K$ and the derivative orders. Bounds for inverses on local gauge slices also depend on the compact subset of the smooth moduli family under consideration. Comparisons near infinity hold at every point of each Hitchin fiber.

\subsection{Local pole types and filtered data}

At an unramified wild pole of order $m\ge2$, the regular semisimple leading coefficient determines the formal eigenlines, including at coincident weights.

At a tame pole we use the domain~\eqref{eq:tame-domain}. The distinguished line of the parabolic flag has weight $1-\alpha$ and the residue acts on it by $-c$. If $c=0$, a nonzero nilpotent residue uniquely determines the parabolic flag, whereas the zero residue permits every flag in the fiber. The moduli problem records the flag, and global stability determines which choices occur in the stable quotient. For the compact flag curves used in Section~\ref{sec:integral} we check stability for every line. The admissible parameter domains and integral markings for the seven families are specified in Section~\ref{sec:integral}.

At the ramified pole of $A_0$, we read the filtration from the leading model
\begin{equation}
 \Phi_{\mathrm{mod}}=\begin{pmatrix}0&w^{-3}\\w^{-4}&0\end{pmatrix}\dd w\,,
 \qquad h\sim\diag(|w|^{1/2},|w|^{3/2})\,.
 \label{eq:ramified-local-model}
\end{equation}
On the cover $w=\zeta^2$, the meromorphic eigenframe~$(1,\zeta^{-1}),(1,-\zeta^{-1})$ of the model gives weight $1/2$ on both pulled-back eigenline bundles. The deck involution exchanges the two bundles and changes the sign of the trace-free diagonal generator, so the component of an invariant section along this generator has vanishing angular average at a ramified point. At a ramified pole of order two, the corresponding matrix is $\left(\begin{smallmatrix}0&w^{-1}\\w^{-2}&0\end{smallmatrix}\right)\dd w$; the eigenvectors and weights are the same. The trace-free eigenvalue one-forms are exchanged with their negatives by the involution. Their residues vanish, since the residue of a meromorphic one-form is invariant under a holomorphic change of coordinate. This accounts for both residues in the three trace-free ramified families. For a general fixed determinant the common scalar residue is prescribed by the determinant. The determinant extension accounts for the scalar weights and does not alter the trace-free deformation problem.

\subsection{Families of solutions for unramified data}
\label{subsec:unramified-families}

We first establish the smooth-family part of Theorem~\ref{thm:variation} for unramified rank-$n$ data over $X$. The formal data remain fixed as $s$ varies in the family of complex structures on $X$. Near a puncture of order $m$, we use a coordinate $w$ with $r=|w|$. We write $\cD_0$ for the model connection with the prescribed singular terms. The trace-free endomorphisms commuting with its diagonal polar part form the subbundle $\mathfrak t_0$. Its orthogonal complement for the model metric is $\mathfrak t^\perp$. The two summands have degrees $0,m$ in the grading of~\cite[Section~2]{BB04}. Write $u_i$ for the component of degree $i$ and $\nabla$ for the derivative induced by the unitary model connection and the smooth background metric on tensors. To display the boundary values separately, put $d_0=1$, $d_m=m$, and use the local norms
\[
 L^2_\sigma=\{u:r^{-\sigma-1}u\in L^2\},\qquad
 \widehat L^{k,2}_\sigma
 =\{u:r^{-d_i(k-\ell)}\nabla^\ell u_i\in L^2_\sigma,\quad
          i=0,m,\ 0\le\ell\le k\}.
\]
The unweighted $L^2$ norm here uses smooth area and the model bundle metric. Use the ordinary Sobolev norm off the polar disks. For $0<\delta<1$ sufficiently small, the perturbation space and the unitary gauge Lie algebra are
\[
 \begin{aligned}
 \mathcal A&=\cD_0+L^{1,2}_{-2+\delta}(\Omega^1\otimes\End_0 E),\qquad
 L^{1,2}_{-2+\delta}=\widehat L^{1,2}_{-2+\delta},\\
 \mathfrak g&=\Big(\bigoplus_p\chi_p\mathfrak t_{0,p}^{\mathrm u}\Big)
       \oplus\widehat L^{2,2}_{-2+\delta}(\mathfrak{su}(E)).
 \end{aligned}
\]
Here $\mathfrak t_{0,p}^{\mathrm u}$ is the skew-Hermitian part of the trace-free centralizer at $p$, extended as constant diagonal matrices in the model frame, and $\chi_p$ is supported in the polar disk and equals one near $p$. The sum describes the full domain $L^{2,2}_{-2+\delta}$ in the convention of Biquard--Boalch. Their weighted Hardy inequality identifies the two norms for one-form perturbations, giving $b\in L^2_{-1+\delta}$. For gauge parameters the trace estimate gives a boundary value in $\mathfrak t_{0,p}^{\mathrm u}$, and subtracting its cutoff extension gives the indicated decaying remainder. These are~\cite[(2.4)--(2.7), (3.3)--(3.4)]{BB04}. The complex gauge Lie algebra is $\mathfrak g\otimes\C$. Since $j_s=j_0$ near every pole, the model and its domains are independent of $s$.

\begin{proposition}
\label{prop:unramified-slices}
Let $X$ be a compact marked curve with prescribed unramified wild formal data having regular semisimple leading coefficients, prescribed weights, and a compatible parabolic determinant. Let $j_s$ be a smooth finite-dimensional family of complex structures that agree near the marked points. On the smooth locus of stable filtered Higgs bundles with the prescribed determinant, using the trace-free deformation complex, solutions to Hitchin's equations admit local gauge slices $q_s(x)$ in the weighted spaces above, and the $L^2$ horizontal projections depend smoothly on $(x,s)$. The gauge Laplacian is invertible on each such slice. Bounds for its inverse are uniform on compact subsets of a smooth local family.
\end{proposition}

\begin{proof}
We use the decomposition $\cD=A+\Psi$. We work with trace-free fields and treat the determinant contribution in Section~\ref{sec:variation}. The flat locus is denoted by $\mathcal Z\subset\mathcal A$. Multiplication by $i$ on complex connections gives $J(a,\psi)=(i\psi,ia)$, independently of $j_s$, and preserves the linearized flatness equation. The quotient of the smooth stable locus of $\mathcal Z$ by the complex gauge group is the de Rham moduli space, independent of $s$~\cite[Section~6 and Theorem~6.1]{BB04}.

For the Hodge star $*_s$ on $X$ and the trace-free projection $\operatorname{pr}_0$, we define the map
\begin{equation*}
 \widehat m_s(\cD)
 =i *_s\operatorname{pr}_0(\dd_A *_s\Psi)
 \in L^2_{-2+\delta}(\mathfrak{su}(E)).
\end{equation*}
This map is smooth on the weighted configuration space. The product estimates needed for its quadratic terms are those of \cite[Lemma~3.5]{BB04}; variation of $*_s$ only changes smooth bounded coefficients on a fixed compact subset away from the punctures. The complex gauge action is smooth by \cite[Section~6]{BB04}.

We choose a smooth local lift $q_0:U\to\mathcal Z$ assigning a solution of the harmonic metric equation to each point of a finite-dimensional moduli chart at $s=0$ and consider the smooth map
\begin{equation*}
 \mathcal F(x,s,\eta)
 =\widehat m_s\bigl(\exp(i\eta)\cdot q_0(x)\bigr),\qquad
 \eta\in\mathfrak g.
\end{equation*}
Here $G_s$ is the $L^2$ inner product on the configuration space, and $\xi^\#=(-\dd_A\xi,[\xi,\Psi])$ is the infinitesimal unitary gauge action. At a solution to Hitchin's equations the derivative with respect to $\eta$ satisfies
\begin{equation*}
 D_\eta\mathcal F=-L_q,\qquad
 \langle L_q\eta,\xi\rangle_{L^2}
       =G_s(\eta^\#,\xi^\#),
\end{equation*}
so the derivative, like the weighted gauge Laplacian, is Fredholm of index zero by~\cite[Lemma~5.1]{BB04}. An element of the kernel is parallel for the unitary connection and commutes with the self-adjoint field, so its eigenspaces give an orthogonal decomposition into parallel filtered subbundles of parabolic degree zero. By filtered stability the trace-free kernel is zero, and $D_\eta\mathcal F$ is an isomorphism from $\mathfrak g$ to $L^2_{-2+\delta}(\mathfrak{su}(E))$. By the implicit function theorem, the solution in the local gauge slice is
\begin{equation*}
 q_s(x)=\exp(i\eta(x,s))\cdot q_0(x)
 \in\mathcal Z\cap\widehat m_s^{-1}(0),
\end{equation*}
with smooth dependence on $x,s$ and on any additional parameters varying the complex structure of the marked curve. Since the same model is used at the poles, this argument applies after differentiation with respect to the parameters to any fixed finite order. Interior elliptic regularity gives smooth representatives away from the punctures. Every compact subset of the smooth total moduli family is covered by finitely many such local slices, with local bounds allowed to depend on that compact subset.

The same inverse $L_q^{-1}$ gives the smooth $L^2$ projection onto unitary gauge directions and hence the horizontal splitting. The horizontal space is $J$-invariant and its induced complex structure on the quotient tangent space is the de Rham complex structure $J$. By \cite[Theorem~5.4]{BB04}, the quotient metric is the $L^2$ metric.

\end{proof}

The real weights can also vary at an unramified wild pole without changing the local operator domain.

\begin{lemma}
\label{lem:weight-domains}
At an unramified rank-two pole of order $m\ge2$ with regular semisimple leading coefficient, the graded weighted domains of the gauge Laplacian have equivalent graph norms on every compact interval of real representatives of the weights. This includes equal and integer weights. After identifying the affine configuration spaces by translation by the model connections, the local slices and horizontal projections vary smoothly with the weights wherever the trace-free kernel vanishes.
\end{lemma}

\begin{proof}
To vary the real weights at an unramified wild pole of order $m\ge2$, we use the following rank-two model. In a fixed unitary bundle on the punctured curve, the connection model is
\begin{equation}
 A_\alpha=A_{\alpha_0}
 +i(\alpha-\alpha_0)\diag(1,-1)\dd\theta,
 \label{eq:weightmodel}
\end{equation}
with the diagonal principal part of the Higgs field unchanged. The angular perturbation is $O(r^{-1})$, whereas on the off-diagonal summand the commutator with the regular semisimple leading term has norm bounded below by $cr^{-m}$. On $r\le\epsilon$, the inequality $r^{-1}\le\epsilon^{m-1}r^{-m}$ makes that perturbation of lower order in the weighted graph norm. The operator on the diagonal summand is independent of the weight difference. On each compact weight interval we can therefore choose a single small positive $\delta$ for which the operator domains have equivalent graph norms, including at equal and integral weights. For variation of the pole data we translate the affine configuration spaces by the corresponding differences of the model connections, using cutoffs to extend the models. We apply the implicit function theorem to all three moment maps together with the Coulomb condition at a reference solution $q$. Indeed, let $R\xi=\xi^\#$ at $q$, let $I_1,I_2,I_3$ denote the ambient complex structures, and set $I_0=\Id$. After choosing the signs of the three moment-map components, the linearization is $P b=((I_jR)^*b)_{j=0}^3$. The quaternionic identities and equivariance of the moment maps at their zero level give
\[
 (I_jR)^*(I_kR)=\delta_{jk}L_q .
\]
For $j\ne k$ the pairing is a moment-map value on a bracket and hence vanishes. To integrate by parts with a constant centralizer boundary value, insert a cutoff $\chi_\epsilon$ whose derivative is supported in $\epsilon<r<2\epsilon$ and satisfies $|\dd\chi_\epsilon|\le C/\epsilon$. The gauge parameter is bounded and $R\xi\in L^2_{-1+\delta}$, so Cauchy--Schwarz bounds the boundary integral by $O(\epsilon^\delta)\|R\xi\|_{L^2_{-1+\delta}}$. Thus $(f_j)\mapsto\sum_{j=0}^3 I_jR L_q^{-1}f_j$ is a bounded right inverse. In the translated spaces the exact diagonal models satisfy the equations on the punctured disks. The remaining terms have smooth compactly supported coefficients or are products controlled by the graded norms. In particular the commutator of a model variation with an off-diagonal field costs at most the factor $r^{-m}$ in the definition of the norm. The parameterized implicit function theorem now gives smooth solution charts, submersive over the pole parameters, and smooth horizontal projections. This argument also applies to the tame and invariant domains below once their inverses are established. Shifting a weight by an integer changes only the holomorphic extension of this filtered family.

\end{proof}

\subsection{Tame poles and completeness in the mixed case}
\label{subsec:mixed}
At a tame point the allowed datum is
\begin{equation*}
 \mathcal P_{\mathrm{tame}}=\{(c,\alpha):0\leq\alpha\leq1/2,
            \ c\ne0\text{ or }0<\alpha<1/2\}.
\end{equation*}
The weights are $(\alpha,1-\alpha)$, and the parabolic line has weight $1-\alpha$ and induced residue $-c$. When $c=0$, we allow the strongly parabolic nilpotent residue to vanish. At zero residue the local incidence condition allows every flag. A limit in the weighted space of solutions will represent a polystable filtered Higgs bundle. The irreducibility hypothesis below makes such a limit stable.

\begin{proposition}
\label{prop:mixed-slices}
Consider stable rank-two filtered Higgs bundles with unramified wild poles as above and tame data in $\mathcal P_{\mathrm{tame}}$. On a compact set $K$ of admissible pole data, there is a common positive lower bound $\nu_K$ for the absolute values of the off-diagonal tame indicial roots. For a sufficiently small weight $0<\delta<\nu_K$, the gauge Laplacian is invertible at each stable point. The local harmonic slices and horizontal projections vary smoothly on these finite-dimensional parameter charts, including points over singular spectral curves and zero residues with their flags. If every polystable filtered Higgs bundle with the prescribed data is irreducible, the $L^2$ metric is complete.
\end{proposition}

\begin{proof}
\emph{Step 1: The indicial roots and inverse.}
In an adapted unitary frame, the squares of the off-diagonal indicial roots of the gauge Laplacian are
\begin{equation*}
 (n+2\alpha-1)^2+16|c|^2,\qquad n\in\Z.
\end{equation*}
Indeed the weight difference shifts the angular Fourier mode by $2\alpha-1$. With $H=\diag(1,-1)$, the model self-adjoint form is $icH\,\dd w/w-i\bar c H\,\dd\bar w/\bar w$. On either root line its radial and angular commutators contribute $16((\Ima c)^2+(\Rea c)^2)$ in the cylindrical metric. At $\alpha=0$ shifting the weight-one lattice reindexes the angular mode by an integer, so the absolute roots are $\sqrt{n^2+16|c|^2}\ge4|c|$. For $c=0$ and $0<\alpha<1/2$, the admissible nilpotent residue term has size $O(r^{-2\alpha})$, which is of lower order than the $r^{-1}$ term at a regular singularity. On each compact subset $K\subset\mathcal P_{\mathrm{tame}}$, a finite cover by charts with nonzero semisimple residue and charts with $0<\alpha<1/2$ gives a common positive lower bound $\nu_K$ for the absolute values of the off-diagonal indicial roots. We choose $0<\delta<\nu_K$ small enough to lie below all the remainder exponents. The diagonal part of the indicial operator is the scalar Laplacian, with constant centralizer boundary values allowed independently at the punctures. The tame off-diagonal summand has degree one. Both summands use $d_i=1$ in the decaying norm, and only the diagonal summand admits constant gauge values. For $c=0$ the first positive off-diagonal indicial root is $\min(2\alpha,1-2\alpha)$. In particular $\delta$ may be chosen below both exponents on a compact chart meeting this stratum.

The local parametrices at the regular and irregular singularities give a Fredholm operator after patching. To see compactness of the lower-order terms in the tame model, write $r^2L$ on the cylinder $t=-\log r$. The nilpotent perturbation has coefficients of positive decay order on the first- and zeroth-order terms. On $t>T$ its operator norm from the decaying weighted $W^{2,2}$ domain to weighted $L^2$ tends to zero as $T\to\infty$. On the remaining compact cylinder the inclusion into $W^{1,2}$ is compact by Rellich's theorem, and the cutoff constants span a finite-dimensional space. A vector in the cokernel solves the adjoint equation in the dual space $L^2_{-\delta}$ for the unweighted pairing. There are no off-diagonal indicial roots between the dual weights. The diagonal expansion allows $C_p+D_p\log r$ at each puncture; Green's formula paired with an independent centralizer value there forces $D_p=0$. The bounded constant and the positive-decay remainder belong to the gauge domain. Formal self-adjointness thus identifies kernel and cokernel. From the quadratic form identity for the linearized gauge operator
\begin{equation*}
 \langle L\upsilon,\upsilon\rangle
 =\|\dd_A\upsilon\|_2^2+\|[\Psi,\upsilon]\|_2^2
 =\|\dd_A\upsilon\|_2^2+2\|[\Phi,\upsilon]\|_2^2
\end{equation*}
we see that a kernel element is parallel and commutes with $\Phi$. Stability makes every such trace-free endomorphism zero. This includes coincident weights with distinct residue eigenvalues.

At tame poles the positive Sobolev weight gives the multiplication estimates in the weighted spaces. At irregular poles we use the graded bounds of \cite[Section~3]{BB04}. These inverse and multiplication estimates allow us to apply the implicit function theorem as in the unramified case. The local gauge slices of solutions and horizontal projections so obtained depend smoothly on the finite-dimensional parameters on the stable locus in the mixed case, including its singular spectral fibers. The same domain admits a zero nilpotent coefficient, since zero belongs to the same space of perturbations with positive decay exponent.

\emph{Step 2: The configuration space and completeness.} We apply~\cite[Theorem~5.4]{BB04} after matching the tame configuration spaces. For $c=0$ we write $\varepsilon=1-2\alpha>0$. A residue that maps the fiber into the parabolic line and annihilates it has norm $O(r^{-1+\varepsilon})$ in the adapted unitary frame. Its derivative and its product with the model connection have norm $O(r^{-2+\varepsilon})$. For $0<\delta<\varepsilon$, the weighted norm is finite because
\begin{equation*}
 \int_0^{r_0}r^{-1+2(\varepsilon-\delta)}\,\dd r<\infty\,.
\end{equation*}
Thus every admissible nilpotent residue term, including zero, defines a perturbation in $L^{1,2}_{-2+\delta}$ around the diagonal model. In an adapted holomorphic frame an ordinary holomorphic entry in the opposite direction is $O(r^{-1+2\alpha})$ in the unitary frame and its first derivative is $O(r^{-2+2\alpha})$. It too lies in the configuration space for the choice $\delta<2\alpha$ above. Taking the diagonal power metric in this local frame, and extending the metric smoothly off the disks, therefore puts the holomorphic Higgs pair in the analytic Dolbeault configuration space. When $c\ne0$ we may instead take a holomorphic residue eigenframe, diagonalizing the Higgs endomorphism near the puncture. At an unramified wild pole the construction is that of~\cite[Lemma~7.3]{BB04}. Conversely, the holomorphic lattice in the weighted complex gauge construction recovers the extension and flag from the growth filtration. In an adapted holomorphic frame, the opposite off-diagonal residue has norm of order $r^{-1-\varepsilon}$ in the adapted unitary frame and is excluded by the domain, as are changes in the induced residues on the graded quotients and in the higher polar terms. At $c\ne0$, holomorphic diagonalization gives the prescribed residue eigenspaces, also at coincident weights. Integer weights are handled by the elementary modifications described above. At a point with $c=0$ the power metric puts the pair in the weighted space, but its contracted Hitchin curvature can grow as $r^{-2+2\varepsilon}$ when the nilpotent residue is nonzero. We construct an initial metric with bounded contracted curvature, as used in the proof of~\cite[Lemma~9.2]{BB04}. This construction works throughout a parameter chart meeting $c=0$. Use a holomorphic frame in which the second vector spans the flag at the pole. The residue in this frame has diagonal entries $c,-c$ and just one possible off-diagonal entry, from the quotient to the flag. Choose the scalar part of the power metric to be the prescribed harmonic determinant metric and take the diagonal connection and the diagonal Higgs field with these two induced residues as the model. Locally a holomorphic determinant frame removes the smooth harmonic scalar factor. On a relatively compact chart with $0<\alpha<1/2$ let $t_*>0$ be a lower bound for $\min(1-2\alpha,2\alpha)$, and decrease $\delta$ below $t_*$. Under $w=R\zeta$, in the unitary frame of the power metric, the model is dilation invariant. The rescaled lower nilpotent coefficient in the one-form has size $O(R^{1-2\alpha}|\zeta|^{-1+(1-2\alpha)})$ and the opposite regular coefficient has size $O(R^{2\alpha}|\zeta|^{-1+2\alpha})$. One covariant derivative costs at most $C/|\zeta|$, and the regular diagonal remainder is $O(R)$. The perturbation on the fixed unit disk is consequently $O(R^{t_*})$ in $L^{1,2}_{-2+\delta}$. In the lower entry this estimate applies also to the flag-preserving off-diagonal residue for $c\ne0$. Derivatives in the weights introduce logarithms, absorbed by using any exponent between $\delta$ and $t_*$. These estimates do not require an eigenframe as $c$ tends to zero.

On the unit disk impose zero Dirichlet data at $|\zeta|=1$ for a trace-free self-adjoint correction to the metric, allowing the usual diagonal limits at the puncture. The derivative of the real Hitchin equation in this variable at the diagonal solution is the model gauge Laplacian. Its Dirichlet inverse is bounded by~\cite[Lemma~4.4]{BB04}; on the present tame chart the common indicial gap and a small perturbation of the coefficients keep this inverse bounded locally in $(c,\alpha)$. The weighted products above make the moment map smooth. The implicit function theorem, applied to the small rescaled perturbation, gives a positive metric in the same weighted complex gauge orbit solving Hitchin's equation on the disk. Its trace-free logarithm relative to the power metric is smooth on the regular annulus by interior elliptic regularity. Cut off that logarithm on the annulus, keeping the corrected metric on a smaller disk. We thus obtain a global adapted reference metric with zero contracted curvature near each such puncture and bounded contracted curvature on the intervening compact region. Near a tame parameter with $c\ne0$ bounded away from zero, including the endpoints of the allowed weight interval, the exact holomorphic eigenframe already gives a diagonal harmonic reference metric near the pole.

The extension of invariant subbundles and equality of analytic and filtered parabolic degrees in~\cite[Lemma~7.5]{BB04} use the growth of the bundle metric. They apply with the same diagonal power metric here. The analytic correspondence~\cite[Theorem~6.1]{BB04} thus supplies a harmonic representative of every stable filtered pair in these weighted spaces. In applying the analytic theorem one may multiply the trace-free Higgs field by a single unit complex number, keeping the prescribed scalar part fixed, so that the two rotated trace-free residues at each $c\ne0$ have distinct imaginary parts. The corresponding model flat residues then have distinct eigenvalues at these points by~\cite[(1.9)--(1.11)]{BB04}: the difference of their imaginary parts is twice the imaginary part of the difference of the rotated Higgs residues. At $c=0$ the two distinct interior weights give distinct real parts modulo $\Z$. At $\alpha=0$, the weight-one eigenline may first be shifted to weight zero in its elementary-modified lattice. The two induced residues still distinguish the eigenlines. The phase is chosen for the existence argument at each fixed datum. The parameter-dependent slice uses the unrotated equations and the common indicial gap. Multiplication by a unit complex number preserves Hitchin's equations and the unitary gauge norms, so it may be undone after constructing the representative. Thus the stability condition on the analytic quotient imposes the global stability condition also when the residue vanishes. In the absence of reducible solutions, completeness on the entire quotient follows from~\cite[Theorem~5.4]{BB04}. This includes the points over singular spectral fibers. We determine the reducibility conditions for the seven families in Section~\ref{sec:integral}.

\end{proof}

\subsection{Ramified poles and the invariant quotient}
\label{subsec:ramified}

Eper--Szab\'o constructed ramified hyperk\"ahler moduli spaces by restricting the weighted configuration space and gauge group to invariant connections and gauge transformations on a ramified cover~\cite[Theorem~2.1, Remark~2.2, and their proof]{EperSzabo25}. Their completeness argument used smoothness at generic parameters. We give the invariant-domain argument needed for the admissible rank-two data here, using stability on the original curve.

\begin{proposition}
\label{prop:ramified-slices}
For the rank-two ramified formal data of Theorem~\ref{thm:variation}, the gauge Laplacian on invariant weighted sections is invertible at every stable filtered Higgs bundle on the original curve. The solutions, local harmonic slices, and horizontal projections depend smoothly on finite-dimensional parameter charts after translating the affine spaces by their model connections. On each polar disk the contribution to the $L^2$ metric and symplectic forms is the integral on its local cover divided by the local degree.
\end{proposition}

\begin{proof}
For the general rank-two data of Theorem~\ref{thm:variation}, we work on local finite Galois covers on which the eigenvalue one-forms are meromorphic and the leading coefficient of the pulled-back Higgs field is regular semisimple. The weighted Sobolev domain of the gauge Laplacian and the dual of the target, $L^2_{-\delta}$, carry the induced finite-group action. Their invariant parts descend to the punctured curve. Averaged local parametrices and invariant cutoffs give a Fredholm operator on the invariant spaces. An invariant cokernel vector solves the adjoint equation against every test section: its pairing with a test section equals the pairing with that section's average. The calculation of the adjoint domain in Proposition~\ref{prop:mixed-slices} applies. The finite-group action preserves the inner product on the centralizer. Green's pairing between the allowed invariant constants and the invariant logarithmic coefficients is therefore nondegenerate. It follows that the logarithmic coefficient vanishes, and the remaining bounded constant term and the terms decaying as a positive power of the radius lie in the original domain.

By the quadratic form identity for the gauge Laplacian, an invariant element of the kernel or cokernel is parallel and commutes with the Higgs field. It descends to a trace-free endomorphism of the filtered Higgs bundle on the original curve, where it vanishes by stability. Existence and uniqueness of the adapted harmonic metric follow from the good filtered correspondence~\cite[Proposition~2.22 and Theorem~2.23]{Moch21}. To see that this metric has a representative in the stated weighted spaces, fix the filtered pair and choose a global double cover of $X$ ramified at its genuinely ramified poles and at an ordinary point where its eigenvalues are distinct. Add another ordinary branch point if necessary to make the branch divisor even in degree. Such a cover exists because a divisor of even degree admits a square root as a line bundle on a compact Riemann surface. The quadratic function field defined by the eigenvalues is ramified at a genuine ramified pole and unramified at the chosen ordinary branch point. It is therefore distinct from the function field of this auxiliary double cover. The eigenvalues remain conjugate of degree two after pullback, so the pulled-back filtered Higgs bundle has no invariant line and is stable. All its marked singularities are now unramified. The unramified or mixed analytic correspondence just established gives a harmonic representative in the weighted spaces, with the pulled-back determinant metric. Start with an invariant adapted reference metric, using the exchanged eigenframes near a ramified pole. Uniqueness with fixed determinant makes the harmonic metric invariant, and the positive gauge relating the two metrics is invariant as well. Near each genuinely ramified pole this gives the required weighted representative on its minimal double cover. A further local cover gives the pullback of that representative, with a sufficiently small positive weight: in the centralizer the smooth model remainder and in its complement the rapid decay of~\cite[Lemmas~4.6 and~5.3]{BB04} are preserved by pullback. The argument supplies a reference solution for the invariant slice; the kernel calculation applies equally on covers where the pullback is only polystable. The inverse on invariant sections, the multiplication estimates, and the moment-map right inverse in the proof of Lemma~\ref{lem:weight-domains} now give smooth local gauge slices of solutions and horizontal projections. For a polar disk at $p$ with local covering degree $d_p$, use $d_p^{-1}$ times the squared weighted norm of the invariant lift there and the ordinary squared norm on the complement, with an invariant partition of unity. The same local averaging gives the $L^2$ metric and symplectic forms on the original curve. On a global cover of degree $d$ it amounts to division of the entire integral by $d$. These local constructions establish the smooth-family assertion of Theorem~\ref{thm:variation} in the ramified case on a compact marked curve.

\end{proof}

For the three families of four-dimensional Hitchin moduli spaces with ramified poles, a global double cover also permits a direct completeness argument.

\begin{proposition}
\label{prop:ramified-quotient}
The filtered Hitchin families~$A_0,D_0,D_1$ have smooth complete $L^2$ hyperk\"ahler moduli spaces for $a\in\C$ in $A_0$, for $p_0,p_\infty\in\C^*$ in $D_0$, and for $\alpha\in\R$, $m\in\C$, $p,b\in\C^*$ in $D_1$. After translation of the affine configuration spaces by the local models, the solutions to Hitchin's equations, the local gauge slices in the weighted Sobolev spaces, and the horizontal projections vary smoothly with the finite-dimensional parameters, including at stable points over singular spectral fibers.
\end{proposition}

\begin{proof}
\emph{Step 1: Invariant spaces on a global cover.} For~$A_0$, we use the coordinate~$w=z^{-1}$ at infinity and the double cover~$w=\zeta^2$. On the cover,
\begin{equation}
 \lambda=-\left(w^{-7/2}+\frac a2w^{-3/2}+\frac u2w^{-1/2}+\cdots\right)\dd w
 =-(2\zeta^{-6}+a\zeta^{-2}+u+\cdots)\dd \zeta\,.
 \label{eq:A0-cover-expansion}
\end{equation}
We work on a fixed double cover of~$\PP^1$ branched at infinity and at an auxiliary unmarked point outside the support of the variations of the complex structure. On this cover the leading coefficient of the pulled-back Higgs field at the pole is regular semisimple. For~$D_0$ and~$D_1$ the cover is~$z=\zeta^2$, branched at zero and infinity. In~$D_0$ the two poles of order two become unramified. In~$D_1$ the pullback unramifies the pole at infinity and takes the pole of order two at $z=0$ to one of order three.

The filtration and eigenframe are those of~\eqref{eq:ramified-local-model}. The deck involution exchanges the eigenline bundles. Its action on the trace-free diagonal generator changes sign, so an invariant diagonal section has zero constant angular coefficient at a ramified point.

The deck group is~$G$, and the covering degree is~$d=2$. We choose invariant model connections at the poles, auxiliary metrics and cutoffs, and take the~$G$-invariant parts of the graded weighted spaces in~\cite[Sections~3--5]{BB04}. The degree of each off-diagonal summand in the weighted norm is determined by the pole order of the corresponding eigenvalue difference, and the gauge domain includes the centralizer boundary values. The metric on this quotient is~$d^{-1}$ times the $L^2$ metric on the cover.

In each of the three spectral equations the quadratic differential has odd order at a ramified point and cannot be the square of a rational one-form. For any~$u$ this rules out a meromorphic $\Phi_{\mathrm{sp}}$-invariant line subbundle. Every filtered Higgs bundle in these families is stable, including those with singular spectral curves. The inverse on the invariant space uses this stability on the original curve. The obstruction space is dual to the infinitesimal automorphism space in the invariant deformation complex, so by the implicit function theorem with parameters the quotient is smooth and its points have representatives in local families of solutions. The hyperk\"ahler reduction identities restrict to invariant tensors and produce the $L^2$ metric.

\emph{Step 2: The auxiliary branch point.} At the ramification point over an auxiliary branch point away from the marked poles, we choose the equivariant structure for which the stabilizer acts trivially on the fiber of the bundle. The holomorphic bundle and Higgs field then descend. We use the ordinary unweighted $W^{1,2}$ domain for connection perturbations at this unmarked point of the cover. In the pullback of a frame on the base, the smooth positive metric on the cover and its inverse are bounded, so the descended metric is adapted to weight zero. By \cite[Theorem~4.1]{Moch21}, the descended stable good filtered Higgs bundle with the auxiliary point left unmarked admits a harmonic metric smooth at that point. Regard both metrics as adapted to the filtration with one additional marked point of weight zero and trivial residue. This addition preserves stability and parabolic degree. Uniqueness in~\cite[Theorem~4.1]{Moch21}, applied with that point included and with the same prescribed determinant, gives equality on the punctured curve. Consequently the metric extends smoothly in the base coordinate, and the local connection has the ordinary unweighted regularity at this point.

\emph{Step 3: Completeness.} For any unit-speed geodesic defined on $0\le s<\ell<\infty$, we take a horizontal lift $q_s=q_0+b_s$ to the invariant configuration space. The operator $D_{q_s}$ is the de Rham differential induced by the flat connection $q_s$ on trace-free endomorphism-valued forms; if $q_s=A_s+\Psi_s$, then $D_{q_s}=\dd_{A_s}+[\Psi_s,\cdot]$. Its $L^2$ adjoint is denoted by $D_{q_s}^*$. Averaging over the deck group shows that the derivative of the invariant lift along the geodesic is orthogonal to every infinitesimal gauge direction on the cover. Hence
\begin{equation}
 (D_{q_s}+D_{q_s}^*)\dot b_s=0\,,\qquad
 \|\dot b_s\|_{L^2}^2=d\,.
 \label{eq:ramified-geodesic}
\end{equation}
We apply the proof of~\cite[Theorem~5.4, equations~(5.9)--(5.12)]{BB04} to this lift, choosing invariant conformal metrics and cutoffs so that the differential operators preserve invariant fields. For $\delta/2<\delta''<\delta$ as in the cited estimates, the estimate following (5.12) bounds $\dot b_s$ and hence $b_s$ in $L^4_{-1+\delta''}$ uniformly. Writing the first equation of~\eqref{eq:ramified-geodesic} relative to $q_0$ gives $(D_{q_0}+D_{q_0}^*)\dot b_s=O(b_s\dot b_s)$. The product is bounded in $L^2_{-2+2\delta''}\subset L^2_{-2+\delta}$. The fixed weighted elliptic estimate for $D_{q_0}+D_{q_0}^*$, with its finite-dimensional kernel controlled by the $L^2$ norm in~\eqref{eq:ramified-geodesic}, now bounds $\dot b_s$ uniformly in $L^{1,2}_{-2+\delta}$ for $0\le s<\ell$. Integration shows that $b_s$ converges in this domain as $s\to\ell$ to an invariant solution of the equations. The quadratic differential has odd order at a ramified point, so the limit on the original curve is irreducible. We can therefore use the inverse on invariant sections to extend the geodesic in the smooth quotient.

\end{proof}
For the pole types of Table~\ref{tab:types}, Proposition~\ref{prop:unramified-slices} and Lemma~\ref{lem:weight-domains} apply to $A_1,D_2$, Proposition~\ref{prop:mixed-slices} to $A_2,D_3$, and Propositions~\ref{prop:ramified-slices} and~\ref{prop:ramified-quotient} to $A_0,D_0,D_1$. Completeness in the unramified case follows from~\cite[Theorem~5.4]{BB04} on the domains without reducible solutions.

\section{Variation of the \texorpdfstring{$L^2$}{L2} metric}
\label{sec:variation}

The variation of the $L^2$ K\"ahler form comes from differentiating the hyperk\"ahler quotient as the complex structure on $X$ changes. The contribution from the varying solution cancels the derivative of the Hamiltonian in the imaginary gauge direction.

\subsection{Compactly supported variation and reduction}

\begin{proof}[Proof of Theorem~\ref{thm:variation}]
On the configuration space the metric and its $J$-K\"ahler form are
\begin{align*}
 G_s((a,\psi),(b,\chi))
 &=\int_X\tr(-a\wedge *_s b+\psi\wedge *_s\chi),\\
 \omega_s^{\mathrm{amb}}((a_1,\psi_1),(a_2,\psi_2))
 &=-i\int_X\tr(\psi_1\wedge *_s a_2-\psi_2\wedge *_s a_1).
\end{align*}
For $\xi\in\mathfrak g$, the moment map is normalized by
\begin{equation*}
 m_s^\xi(\cD)=-i\int_X\tr(\xi\,\dd_A *_s\Psi),\qquad
 \iota_{\xi^\#}\omega_s^{\mathrm{amb}}=-\dd m_s^\xi,\qquad
 \xi^\#=(-\dd_A\xi,[\xi,\Psi]).
\end{equation*}
Since $\dd_A(*_s\Psi)\in L^2_{-2+\delta}\subset L^1$ and the gauge parameters are bounded, the integral defines the moment map for the weighted unitary gauge action. We can also compute the differential of this moment map using radial cutoffs. The embedding in \cite[Section~3]{BB04} implies that $b\in L^2_{-1+\delta}$ for each deformation. For a cutoff with $|\dd\chi_\epsilon|\leq C/r$, the boundary term satisfies
\begin{equation*}
 \left|\int_{r\asymp\epsilon}
          \dd\chi_\epsilon\,\tr(\xi *_s b)\right|
 \leq C\|\xi\|_\infty
 \|r^{-\delta}b\|_{L^2(r\asymp\epsilon)}\,\epsilon^\delta
 \longrightarrow0.
\end{equation*}
The cutoff estimate justifies the moment map identity for tangent vectors in $L^{1,2}_{-2+\delta}$ and gives
\begin{equation*}
 \dd m_s^\xi(J\eta^\#)
 =\omega_s^{\mathrm{amb}}(J\eta^\#,\xi^\#)
 =-G_s(\eta^\#,\xi^\#),
\end{equation*}
which verifies $D_\eta\mathcal F=-L_q$ in Section~\ref{sec:families}.

On the flat locus the de Rham holomorphic symplectic form is
\[
 \omega_K+i\omega_I=-i\int_X\tr(\delta\cD_1^0\wedge\delta\cD_2^0).
\]
Indeed, in the decomposition $\delta\cD=a+\psi$ we have
$I(a,\psi)=(*_s a,-*_s\psi)$ and $K=IJ$. Expanding the integral gives the
two indicated real forms from $G_s$. The tangent vectors are in $L^2$,
and the preceding cutoff estimate shows that integration by parts against a
bounded infinitesimal gauge transformation has zero boundary contribution.
Flatness therefore makes the pairing independent of the representatives.
Its expression on the flat quotient, and multiplication by $i$ there,
are independent of $j_s$.

The variation of the Hodge star has compact support, so we can define
\begin{equation*}
 H_s(\cD)=\frac12\int_X\tr(\Psi\wedge\dot *_s\Psi).
\end{equation*}
It is a smooth function invariant under unitary gauge transformations on $\mathcal A$. Direct differentiation, with $\dd_J^c=-J^*\dd$, gives
\begin{equation*}
 \dd_J^cH_s(a,\psi)=-i\int_X\tr(\Psi\wedge\dot *_s a),\qquad
 \partial_s\omega_s^{\mathrm{amb}}=\ddc H_s.
\end{equation*}
Because $\dot *_s$ has compact support, integration by parts gives the following identity for every $\xi\in\mathfrak g$:
\begin{align*}
 \dd_J^cH_s(\xi^\#)
 &=i\int_X\tr(\Psi\wedge\dot *_s \dd_A\xi)\\
 &=i\int_X\tr(\dd_A\xi\wedge\dot *_s\Psi)\\
 &=-i\int_X\tr(\xi\,\dd_A\dot *_s\Psi)
 =\dot m_s^\xi.
\end{align*}
The second equality uses the symmetry of $(b_1,b_2)\mapsto\int_X\tr(b_1\wedge\dot *_s b_2)$.

We use the smooth lift $q_s:U\to\mathcal Z\cap m_s^{-1}(0)$ constructed in Section~\ref{sec:families}. At fixed $x\in U$, its derivative $V_s=\partial_s q_s$ is tangent to the complex gauge orbit and has the form $V_s=\xi_s^\#+J\eta_s^\#$. For $v\in TU$, we decompose
\begin{equation*}
 (q_s)_*v=h_v+\upsilon_v^\#,\qquad
 (q_s)_*(Jv)=Jh_v+\upsilon_{Jv}^\#,
\end{equation*}
where $h_v$ is $L^2$-horizontal. The horizontal space is orthogonal to both unitary and imaginary gauge directions. Gauge invariance of $H_s$ gives
\begin{align*}
 q_s^*\dd_J^cH_s(v)-\dd_J^c(H_s\circ q_s)(v)
 &=-\dd H_s(J\upsilon_v^\#)\\
 &=\dd_J^cH_s(\upsilon_v^\#)
 =\dot m_s^{\upsilon_v}.
\end{align*}
We differentiate $m_s(q_s)=0$, keeping $\upsilon_v$ independent of $s$ at the point, and obtain
\begin{equation*}
 0=\dot m_s^{\upsilon_v}
   +\dd m_s^{\upsilon_v}(V_s)
  =\dot m_s^{\upsilon_v}
   +\omega_s^{\mathrm{amb}}(V_s,\upsilon_v^\#).
\end{equation*}
The pairing of $V_s$ with $h_v$ vanishes. Consequently
\begin{equation*}
 q_s^*(\iota_{V_s}\omega_s^{\mathrm{amb}})(v)
 =-\dot m_s^{\upsilon_v},
\end{equation*}
which cancels the vertical term in the pullback of $\dd_J^cH_s$. By Cartan's formula we therefore have
\begin{align*}
 \partial_s\omega_{J,s}
 &=\partial_s(q_s^*\omega_s^{\mathrm{amb}})\\
 &=\dd\bigl(q_s^*\dd_J^cH_s+
       q_s^*(\iota_{V_s}\omega_s^{\mathrm{amb}})\bigr)\\
 &=\ddc(H_s\circ q_s).
\end{align*}

At a parameter value $s$, if $\dd z+\varepsilon\mu\,\dd\bar z$ generates the deformed holomorphic cotangent line bundle to first order on $X$, then
\begin{equation*}
 \dot *\dd z=-2i\mu\,\dd\bar z,\qquad
 \dot *\dd\bar z=2i\bar\mu\,\dd z.
\end{equation*}
In the local expressions $\Phi=\varphi\,\dd z$ and $\Psi=\Phi+\Phi^*$, the two-form integrated in the definition of the Hamiltonian is
\begin{equation*}
 \frac12\tr(\Psi\wedge\dot *\Psi)
 =-i\left(\tr(\varphi^2)\mu+
            \overline{\tr(\varphi^2)\mu}\right)\dd z\wedge \dd\bar z.
\end{equation*}
Since conjugation reverses $\dd z\wedge \dd\bar z$, integration gives
\begin{equation}
 H_s(q_s)=2\Ima\int_X\tr(\Phi^2)\mu.
 \label{eq:domainenergy}
\end{equation}
For a nontrivial prescribed flat parabolic determinant line bundle, the reduction calculation is made with $\Psi^0$ and $\Phi^0$. In the decomposition
\begin{equation*}
 \tr(\Phi^2)=\tr((\Phi^0)^2)+\frac1n(\tr\Phi)^2,
\end{equation*}
the last term is constant along the fixed-determinant moduli space, so its $\ddc$ vanishes. The same formula therefore holds when the determinant metric varies with $j_s$, including at stable points with singular spectral curves. At tame poles the annular cutoff estimate above applies with the positive weight of Section~\ref{subsec:mixed}. In the ramified case use the invariant local slices and divide polar-disk pairings by their local degrees as in Section~\ref{subsec:ramified}. Since $\dot *_s$ is supported away from all marked points, $H_s$ is the integral on $X$ displayed above. If a global cover $\pi:\widetilde X\to X$ of degree $d$ is used, the identical expression is $(2d)^{-1}\int_{\widetilde X}\pi^*\tr(\Psi\wedge\dot *_s\Psi)$. This proves~\eqref{eq:main-variation} in each case.
\end{proof}

\begin{corollary}
\label{cor:residuevariation}
In the setting of Theorem~\ref{thm:variation}, we assume that local vector fields $v_p$, holomorphic near the marked points, induce the infinitesimal deformation of the formal data and that pullback by the inverse change of coordinates fixes the formal data to first order. Then, with $\Phi=\Psi^{1,0}$ as in Theorem~\ref{thm:variation},
\begin{equation}
 \dot\omega_J=-4\pi\ddc\Rea
 \sum_p\Res_p\bigl(\iota_{v_p}\tr(\Phi^2)\bigr).
 \label{eq:generalresiduevariation}
\end{equation}
\end{corollary}

\begin{proof}
Choose a smooth field $v$ on $X$ that agrees with $v_p$ near each marked point, using disjoint cutoffs. For this infinitesimal inverse pullback, the Beltrami differential is $\mu=-\bar\partial v$ and is supported away from the marked points. Since $q=\tr(\Phi^2)$ is meromorphic and holomorphic on the complement of the poles, we obtain by Stokes' theorem
\begin{equation}
 \int_Xq\mu=-2\pi i\sum_p\Res_p(\iota_{v_p}q).
\end{equation}
Twice its imaginary part is $-4\pi\Rea\sum_p\Res_p(\iota_{v_p}q)$. Substitution into~\eqref{eq:main-variation} gives~\eqref{eq:generalresiduevariation}.
\end{proof}

With the two Higgs fields related by~\eqref{eq:phase}, we obtain from Corollary~\ref{cor:residuevariation}
\begin{equation}
 \dot\omega_J=4\pi\ddc\Rea\sum_p\Res_p\bigl(\iota_{v_p}(2\lambda^2)\bigr).
 \label{eq:spectral-residue-H}
\end{equation}
A residue evaluated on a local cover on which the eigenvalue one-forms are meromorphic is divided by the local degree of that cover.

\subsection{The three ALG Hamiltonians}
\label{subsec:quarticvariation}

For the Higgs field in~\eqref{eq:spectral}, the trace of its square is
\begin{equation}
 q_{\mathrm{sp}}=\tr(\Phi_{\mathrm{sp}}^2)
 =2\left(w^{-8}+aw^{-6}+bw^{-5}+uw^{-4}\right)\dd w^2.
 \label{eq:qold}
\end{equation}
The coordinate vector field
\begin{equation}
 v=-\frac{\dot a}{2}w^3\partial_w
 \label{eq:jetvector}
\end{equation}
satisfies
\begin{align}
 \mathcal L_v(-w^{-4}\dd w)&=-\frac{\dot a}{2}w^{-2}\dd w,\label{eq:liepolar}\\
 \Res_{w=0}(\iota_v q_{\mathrm{sp}})&=-\dot a\,u.\label{eq:resq}
\end{align}
The inverse pullback generated by~\eqref{eq:jetvector} preserves the prescribed principal parts by~\eqref{eq:liepolar}. Substituting~\eqref{eq:resq} into the residue formula~\eqref{eq:spectral-residue-H} gives the Hamiltonian and variation
\begin{equation}
 H_s(q_s)=-4\pi\Rea(\dot a\,u),\qquad
 \partial_s\omega_J=-4\pi\ddc\Rea(\dot a\,u).
 \label{eq:quarticactual}
\end{equation}
This identity holds for all $a\in\C$ and therefore applies to the K\"ahler forms pulled back by the global transport maps constructed in Section~\ref{sec:finite}.

For $A_0$, the coordinate vector field and its contraction are
\begin{align*}
 v&=-\dot a\,w^3\partial_w,
 &2\lambda^2&=2(w^{-7}+aw^{-5}+uw^{-4})\dd w^2,\\
 \Res_{w=0}\bigl(\iota_v(2\lambda^2)\bigr)&=-2\dot a\,u.
\end{align*}
On the double cover $w=\zeta^2$, the lifted field is $-(\dot a/2)\zeta^5\partial_\zeta$ and the residue is $-4\dot a\,u$. Division by two gives the same result. Thus
\begin{equation}
 \dot\omega_J=-8\pi\ddc\Rea(\dot a\,u).
 \label{eq:A0-Hamiltonian}
\end{equation}

For $A_2$, the Laurent expansion at $w=z^{-1}=0$ is
\begin{equation*}
 \lambda_+=-\left(w^{-3}+\frac a2w^{-2}+\mu w^{-1}+O(1)\right)\dd w.
\end{equation*}
The field $v=-(\dot a/2)w^2\partial_w$ changes the coefficient $a$ and preserves $\mu$. Since
\begin{equation*}
 2\lambda^2=2\left(w^{-6}+aw^{-5}+(a^2/4+2\mu)w^{-4}+uw^{-3}+c^2w^{-2}\right)\dd w^2,
\end{equation*}
its contraction has residue $-\dot a\,u$. Consequently
\begin{equation}
 \dot\omega_J=-4\pi\ddc\Rea(\dot a\,u).
 \label{eq:A2-Hamiltonian}
\end{equation}
For the three ALG families, the common notation is therefore
\begin{equation}
 \dot\omega_J=-K_j\ddc\Rea(\dot a\,u),\qquad
 (K_0,K_1,K_2)=(8\pi,4\pi,4\pi).
 \label{eq:ALG-Hamiltonians}
\end{equation}

\section{Characteristic classes and compact periods}
\label{sec:compact}

We compute real periods by Chern--Weil transgression and complex periods from changes of the Liouville form on overlaps. The formulas~\eqref{eq:CW} and~\eqref{eq:complexclass} express both periods through the discriminant and the labeled eigenline bundle degrees. Section~\ref{sec:integral} applies these formulas to integral bases in the seven families.

\subsection{A Chern--Weil formula for families parametrized by compact surfaces}

We work over a closed oriented smooth surface $S$ equipped with a smooth map into a smooth rank-two wild Hitchin moduli space with fixed unramified regular semisimple formal data. Locally we choose universal vector bundles that lift the universal projective bundle over $X\times S$, and denote each lift by $\mathcal E$. For a marked point $p$, $L_{pj}\to S$ denotes the $j$th labeled formal eigenline bundle of a local universal lift at $p$. We use the characteristic classes
\begin{equation}
 \Delta(\mathcal E)=c_2(\mathcal E)-\frac14c_1(\mathcal E)^2,
 \qquad
 c_1(L_{pj}^0)=c_1(L_{pj})-\frac12c_1(\mathcal E|_{\{p\}\times S}).
 \label{eq:projectiveclasses}
\end{equation}
The notation $c_1(L_{pj}^0)$ denotes the displayed rational characteristic class; $L_{pj}^0$ is used only through this class. Since the classes are unchanged when a local universal vector bundle is tensored with a line bundle from $S$, they descend to the universal projective bundle.

\begin{lemma}
\label{lem:boundarymetric}
The expansion below holds for any finite-dimensional local holomorphic family $(E_x,\Phi_{\mathrm{sp},x})$ of stable rank-$n$ Higgs bundles with unramified irregular poles and regular semisimple leading coefficients, provided that the prescribed principal parts and parabolic weights are independent of $x$. Near a pole of order $m\ge2$, in a holomorphic frame that diagonalizes $\Phi_{\mathrm{sp},x}$, we write
\begin{equation*}
 D_\alpha(r)=\diag(r^{\alpha_1},\ldots,r^{\alpha_n}).
\end{equation*}
The harmonic metrics $h_x$ have the form
\begin{equation}
 h_x=D_\alpha(r)\bigl(H(x)+R_x(w)\bigr)D_\alpha(r),
 \label{eq:leadingmetric}
\end{equation}
where $H(x)$ is smooth, positive and diagonal. The remainder is $O(r^\delta)$ for some sufficiently small $\delta>0$, uniformly on compact subsets of the parameter charts, with the same bound after any fixed finite number of derivatives in the parameters and along circles about the pole on the marked curve. The off-diagonal entries of the remainder vanish to every order. The universal Chern connection consequently satisfies
\begin{equation}
 \mathcal A|_{r=\epsilon}
 =\diag(i\alpha_j\dd\theta+B_j)+O(\epsilon^\delta),
 \label{eq:universalboundary}
\end{equation}
with the same bounds after differentiating any fixed number of times with respect to the parameters; $B_j$ is the Chern connection of the metric $H_j(x)$ on the eigenline bundle.
\end{lemma}

\begin{proof}
Regular semisimplicity gives such a holomorphic eigenframe because $w^m\Phi_{\mathrm{sp},x}/\dd w$ has distinct eigenvalues at $w=0$. In that frame $h_0=D_\alpha^2$ is the local harmonic metric for the diagonal Higgs field. If $(\sigma_j)$ is the chosen holomorphic eigenframe, the associated $h_0$-unitary frame is $e_j=r^{-\alpha_j}\sigma_j$. In the frame $(e_j)$ the model unitary connection is $i\diag(\alpha_j)\dd\theta$; the regular diagonal terms of $\Phi_{\mathrm{sp},x}$ are smooth remainders.

The complex gauge equation for the Higgs bundle has the same positive gauge Laplacian as the reduction in Section~\ref{subsec:unramified-families}. By the weighted implicit function theorem in Proposition~\ref{prop:unramified-slices}, a correcting complex gauge transformation can be chosen to depend smoothly on $x$. In this unitary frame, the embedding of \cite[Section~3]{BB04} gives bounded diagonal limits and a bounded inverse. Continuity of the boundary-value map gives the corresponding bounds after differentiation with respect to the parameters.

On a fixed sufficiently small disk, we obtain local unitary Coulomb gauges by the implicit function argument in the proof of \cite[Lemma~5.3]{BB04}. These gauges vary smoothly on a compact subset of a parameter chart, and the harmonic connection has the form $\cD_0+a$. By \cite[Lemma~4.6]{BB04},
\begin{equation*}
 a^\perp=O(r^N)\quad\text{for every }N,
 \qquad a^{\mathrm{diag}}=O(r^{-1+\delta}).
\end{equation*}
Interior elliptic estimates on balls of radius comparable to $r^m$ give rapid decay of derivatives of the off-diagonal terms on the punctured curve, by first choosing a sufficiently large decay order. In the diagonal part of the Coulomb and harmonic equations, $(\dd+\dd^*)a^{\mathrm{diag}}$ is a smooth form vanishing to every order at the pole. Every nonzero quadratic term in this projection contains two off-diagonal factors. Since $a^{\mathrm{diag}}=O(r^{-1+\delta})$, the diagonal component belongs to $L^p$ for some $p>2$. This integrability rules out a distributional term supported at the puncture in its equation. First-order elliptic removability then gives a smooth extension across the pole.

The derivative of $(\cD_0+\cD_0^*)a=Q(a,a)$ with respect to a parameter has the same model differential operator at the pole. In the off-diagonal equation the smooth diagonal remainder is lower order, and the inhomogeneous terms contain previously controlled derivatives or a rapidly decaying off-diagonal factor. The off-diagonal parameter derivative therefore satisfies the improved weighted estimate of \cite[Corollaries~4.2--4.3 and Lemma~4.4]{BB04}. The equation on the diagonal summand then has a smooth right-hand side, so its solution again extends across the puncture. Finite induction proves these statements with parameter derivatives of every fixed finite order.

We first apply the gauge transformation correcting the Higgs bundle and then the unitary Coulomb gauge. We denote their composition by $g$ and the transformed Higgs field by $\Phi_{\mathrm{sp}}'$. Its labeled eigenvalues remain $\lambda_j$. The polynomial formula for its spectral projectors and the separation $|\lambda_i-\lambda_j|\ge cr^{-m}$ give
\begin{equation*}
 P_j=E_{jj}+O(r^N)\qquad\text{for every }N>0.
\end{equation*}
The remainder and its derivatives of any fixed finite order in the local coordinate and the parameters vanish to infinite order at $w=0$. For the matrix $W=(P_1e_1,\ldots,P_ne_n)$, the same holds for $W-\Id$, and $\Phi_{\mathrm{sp}}'W=W\diag(\lambda_j)$. The transformation $g$ sends each original eigenline bundle to the corresponding labeled eigenline bundle spanned by a column of the matrix just constructed. Thus $g=Wd$ for a diagonal matrix $d$. All these expressions are in the $h_0$-unitary frame.

Here $A'$ denotes the unitary connection of the transformed harmonic pair and $A_0$ the model unitary connection. The $\bar\partial$ gauge equation becomes
\begin{equation*}
 \bar\partial d=-Cd,
 \qquad C=W^{-1}\bar\partial W+W^{-1}A'^{0,1}W-A_0^{0,1}.
\end{equation*}
The diagonal part of $A'-A_0$ is smooth, while its off-diagonal part and the commutator of $W-\Id$ with the $r^{-1}$ model connection vanish to every order. Thus $C$ is smooth, and the equation for the invertible diagonal matrix $d$ shows that $C$ is diagonal. Since $d$ is bounded, its $\bar\partial$ equation holds distributionally across the puncture, where elliptic regularity gives a smooth extension of $d$. The differentiated equations give the same conclusion for parameter derivatives. The inverse $d^{-1}$ remains bounded, so $d(0,x)$ is invertible.

Unitary left multiplication has not changed the harmonic metric. The matrix of the harmonic metric relative to the model metric is
\begin{equation*}
 g^*g=d^*W^*Wd,
 \qquad H(x)=d(0,x)^*d(0,x).
\end{equation*}
The relative metric matrix has smooth diagonal entries with positive limits at the pole, and its off-diagonal entries vanish to every order. In the holomorphic frame the metric therefore has the form~\eqref{eq:leadingmetric}. The formula for the Chern connection gives $i\alpha_j\dd\theta$ from $r^{2\alpha_j}$ and $B_j$ from $H_j(x)$; the remaining terms on $r=\epsilon$ tend to zero with their tangential and parameter derivatives. This proves~\eqref{eq:universalboundary}.
\end{proof}

\begin{proposition}
\label{prop:CW}
For any smooth map from a closed oriented surface $S$ to the smooth moduli space of rank-two Higgs bundles with prescribed unramified regular semisimple formal data and weights, let $\mathcal E$ denote the pullback of the universal projective bundle to $X\times S$. The period in the normalization~\eqref{eq:fiberperiod} is given by
\begin{equation}
 \frac1{4\pi^2}\int_S\omega_I
 =\int_{X\times S}\Delta(\mathcal E)
  -\sum_{p,j}\alpha_{pj}\int_Sc_1(L_{pj}^0).
 \label{eq:CW}
\end{equation}
\end{proposition}

\begin{proof}
We pull back to $S$ the local complex universal families on moduli charts, their harmonic metrics, and the Chern connections $\mathcal A$ on the local universal vector bundles $\mathcal E$. For a connection $\mathcal A$ with anti-Hermitian curvature, its discriminant Chern--Weil form is $\Delta(\mathcal A)=(8\pi^2)^{-1}\tr(F_{\mathcal A}^0\wedge F_{\mathcal A}^0)$. We compute the differential forms on the complex charts before pulling them back, so the map from $S$ need not be holomorphic. By Lemma~\ref{lem:boundarymetric}, the boundary expansion~\eqref{eq:universalboundary} also holds after differentiation with respect to the parameters.

In local coordinates $w_p$ at the punctures, the truncated surface is $X_\epsilon=X\setminus\bigcup_p\{|w_p|<\epsilon\}$. On this surface, the curvature of the universal connection gives
\begin{equation}
 \omega_I^\epsilon
 =4\pi^2\int_{X_\epsilon}\Delta(\mathcal A)+\dd_S\vartheta_\epsilon,
 \label{eq:cutoffCW}
\end{equation}
where $\omega_I^\epsilon$ denotes the integral of the harmonic tangent pairing over $X_\epsilon$. To identify the exact term, we compute on a holomorphic disk in the moduli chart and then use polarization. In a local complex parameter $v$, we write $\eta=h^{-1}\partial_v h$. The universal connection represents the tangent vector by
\begin{equation}
 (F_{v\bar z},\nabla_v\Phi_{\mathrm{sp}})
 =(-\bar\partial_z\eta,\partial_v\Phi_{\mathrm{sp}}+[\eta,\Phi_{\mathrm{sp}}]).
 \label{eq:universaltangent}
\end{equation}
Write $\Phi_{\mathrm{sp}}^0=\varphi\,\dd z$. Holomorphicity in $v$ gives $\nabla_v\varphi^*=0$. Taking the trace-free part of the Bianchi identity gives $\nabla_v F^0_{z\bar z}=\nabla_z F^0_{v\bar z}$. Differentiation of the coefficient equation $F^0_{z\bar z}+[\varphi,\varphi^*]=0$ therefore yields
\begin{equation}
 \nabla_z F^0_{v\bar z}+[\nabla_v\varphi,\varphi^*]=0.
\end{equation}
By the K\"ahler identities this is the harmonic gauge condition for~\eqref{eq:universaltangent}. The mixed components of the curvature of the universal connection and $\nabla_v\Phi_{\mathrm{sp}}$ transform tensorially under parameter-dependent holomorphic changes of frame, so the calculation patches between the local universal families. Write $v=x+iy$, $B=F^0_{v\bar z}$, $f=\nabla_v\varphi$, and $\dd A_z=(i/2)\dd z\wedge\dd\bar z$. The real tangent in the $x$ direction has unitary connection component $B\,\dd\bar z-B^*\dd z$ and self-adjoint component $if\,\dd z-if^*\dd\bar z$. Its $I$-image is the $y$ tangent. The metric of Section~\ref{sec:variation} therefore gives
\[
 \omega_I^\epsilon(\partial_x,\partial_y)
       =4\int_{X_\epsilon}\tr(BB^*+ff^*)\,\dd A_z.
\]
Since $F^0_{z\bar v}=B^*$ and $\dd v\wedge\dd\bar v=-2i\,\dd x\wedge\dd y$, the Chern--Weil term is
\begin{equation}
 4\pi^2\!\left(\int_{X_\epsilon}\Delta(\mathcal A)\right)(\partial_x,\partial_y)
   =4\int_{X_\epsilon}\tr(BB^*-F^0_{z\bar z}F^0_{v\bar v})\,\dd A_z.
 \label{eq:curvatureterms}
\end{equation}
The Chern connection and Hitchin's equation give
\begin{equation}
 \partial_v\partial_{\bar v}\tr(\varphi\varphi^*)
   =\tr(ff^*)-\tr(F^0_{v\bar v}[\varphi,\varphi^*])
   =\tr(ff^*+F^0_{v\bar v}F^0_{z\bar z}).
 \label{eq:higgsnormsecond}
\end{equation}
Thus the difference of the two pairings is $\dd\dd_I^cN_\epsilon$, where
\[
 N_\epsilon=\frac i2\int_{X_\epsilon}
      \tr\bigl(\Phi_{\mathrm{sp}}^0\wedge(\Phi_{\mathrm{sp}}^0)^*\bigr).
\]
Here $\dd\dd_I^cN_\epsilon(\partial_x,\partial_y)=4\partial_v\partial_{\bar v}N_\epsilon$. The cutoff function $N_\epsilon$ is globally defined by gauge invariance. Pulling the one-form $\dd_I^cN_\epsilon$ back to $S$ gives $\vartheta_\epsilon$. Since $S$ is closed,
\begin{equation}
 \int_S\dd_S\vartheta_\epsilon=0
 \label{eq:exactcutoff}
\end{equation}
for every $\epsilon>0$. We take the limit after integrating over $S$. Since the $L^2$ norms of the tangent representatives on the deleted disks around the punctures tend to zero uniformly on the compact surface $S$, we have $\int_S\omega_I^\epsilon\to\int_S\omega_I$ as $\epsilon\downarrow0$.

Choose a smooth reference connection on the universal projective bundle whose restriction near each pole is the projectivization of the direct sum of the $B_j$. In a local vector-bundle lift we denote a lift of this connection by $\mathcal A_0$. Differences of its trace-free part and the trace-free universal connection, and their invariant trace polynomials, patch globally. The difference of the trace-free parts of the two connections is $a=\mathcal A^0-\mathcal A_0^0$. Here $F_0$ and $\dd_0$ denote the curvature and covariant derivative of $\mathcal A_0^0$. The transgression formula is
\begin{align}
 \Delta(\mathcal A)-\Delta(\mathcal A_0)&=\dd\operatorname{CS}(\mathcal A,\mathcal A_0),\label{eq:CStrans}\\
 \operatorname{CS}(\mathcal A,\mathcal A_0)
 &=\frac1{8\pi^2}\tr\left(2a\wedge F_0+a\wedge\dd_0a+\frac23a^3\right).
 \label{eq:CS}
\end{align}
On $r=\epsilon$, with $\bar\alpha=\frac12\sum_j\alpha_j$, we have
\begin{equation}
 a=i\diag(\alpha_j-\bar\alpha)\dd\theta+O(\epsilon^\delta).
\end{equation}
The leading quadratic and cubic terms of~\eqref{eq:CS} vanish, and
\begin{equation}
 \operatorname{CS}(\mathcal A,\mathcal A_0)
 =\frac{i}{4\pi^2}\sum_j(\alpha_j-\bar\alpha)\dd\theta\wedge F_{B_j}
  +O(\epsilon^\delta).
 \label{eq:CSboundary}
\end{equation}
The clockwise orientation of each inner boundary circle and $c_1(L_j)=\frac{i}{2\pi}[F_{B_j}]$ therefore give
\begin{align}
 \lim_{\epsilon\downarrow0}\int_{X_\epsilon\times S}\Delta(\mathcal A)
 &=\int_{X\times S}\Delta(\mathcal E)
   -\sum_j(\alpha_j-\bar\alpha)\deg L_j\notag\\
 &=\int_{X\times S}\Delta(\mathcal E)
   -\sum_j\alpha_j\deg L_j^0.
 \label{eq:CSlimit}
\end{align}
Summing the contributions from the poles and applying~\eqref{eq:cutoffCW} and~\eqref{eq:exactcutoff} proves~\eqref{eq:CW}.
\end{proof}

\subsection{Normalization on a smooth fiber and boundary terms at tame poles}
\label{subsec:calibration}
For a smooth genus-one spectral double cover $\pi:\Sigma\to\PP^1$, we write $F=\Pic^d(\Sigma)$. A choice of line bundle $L_0$ of degree $d$ gives, by tensoring with $L_0$, an isomorphism from the Jacobian to $F$. We use the degree-zero Poincar\'e line bundle $P$ on $\Sigma\times F$, trivialized on $\{p_*\}\times F$ and on $\Sigma\times\{[L_0]\}$ under this identification. With the orientation determined by a symplectic basis,
\begin{equation*}
 \int_{\Sigma\times F}c_1(P)^2=-2.
\end{equation*}
We apply Grothendieck--Riemann--Roch~\cite[Section~7]{BS58} to compute the discriminant integral for the pushforward
\begin{equation*}
 \mathcal E=(\pi\times\Id_F)_*(\operatorname{pr}_\Sigma^*L_0\otimes P).
\end{equation*}
We obtain
\begin{equation*}
 \int_{\PP^1\times F}\left(c_2(\mathcal E)-\frac14c_1(\mathcal E)^2\right)=1.
\end{equation*}
Only $c_1(P)^2/2$ contributes in degree two along $F$ to this integral. Its integral is $-1$, giving the integral of $\operatorname{ch}_2$. Up to a line bundle of degree zero the determinant is pulled back from $\PP^1$, so the square of its first Chern class vanishes. Evaluation of $P$ at a fixed point of $\Sigma$ gives a line of degree zero on $F$. The induced eigenlines have that degree, including the graded lines at a simple ramification point, where taking a first jet also tensors with a fixed one-dimensional vector space. In the unramified wild case~\eqref{eq:CW} now gives $\int_F\omega_I=4\pi^2$. The tame and ramified versions use the boundary and descent calculations below.

\paragraph{Boundary estimate at a tame pole.}
Fix the pole parameters $c,\alpha$ throughout the family parametrized by $S$. Write $L_{pj}$ for the graded lines of its parabolic flag. For $c\ne0$ these are canonically identified with the two residue eigenlines. When $c=0$, the two graded lines are the parabolic line and the quotient of the fiber by that line, whether the full residue is zero or a nonzero nilpotent. We verify the boundary estimate used in the transgression. Locally on the parameter chart, choose a holomorphic frame diagonalizing the Higgs field for $c\ne0$, and a holomorphic frame adapted at the puncture to the flag for $c=0$. Use $D_\alpha^2$ as reference metric. In its unitary frame the model angular connection is $i\diag(\alpha_j)\dd\theta$. In the first case the remainder of $w\Phi_{\mathrm{sp}}/\dd w$ beyond its constant diagonal residue is $O(r)$. In the second case put $e=1-2\alpha\in(0,1)$. The nilpotent residue from the quotient to the flag is $O(r^e)$ in this rescaled coefficient, and the opposite regular entry is $O(r^{1-e})$. These estimates hold with derivatives in the parameters and in $r\partial_r,\partial_\theta$.

In the chosen holomorphic frame put $M=D_\alpha^{-1}hD_\alpha^{-1}$. By Step~2 in the proof of Proposition~\ref{prop:mixed-slices}, the harmonic metric $h_{x_0}$ at a basepoint $x_0\in S$ lies locally in the weighted complex gauge orbit of $D_\alpha^2$. The preceding coefficient bounds place the holomorphic family in the fixed weighted configuration space. At $h_{x_0}$ the derivative of the metric equation in the self-adjoint gauge direction is the invertible gauge Laplacian of that proposition. The implicit function theorem therefore gives correcting gauges that depend smoothly on the parameters in its weighted domain. The weighted embedding on the regular singular summands~\cite[Section~3]{BB04} gives $M=H+O(r^{\delta'})$, where $H$ is positive diagonal and $0<\delta'<\delta$ may be decreased. The same bound holds for the remainders of all fixed parameter derivatives. Here the off-diagonal summand has no constant term, and the diagonal constants are precisely the boundary values admitted by the domain.

For the normal derivative needed in the angular Chern connection, rescale $w=r\zeta$ on $1/2<|\zeta|<2$. The weighted domain bounds $M-H$ in $W^{2,2}$ on each such annulus by $O(r^{\delta'})$. In this unitary frame the leading tame coefficients and their derivatives are uniformly bounded on the fixed annulus. The constant diagonal matrix $H$ solves the leading equation. Substitution in the full equation produces an $O(r^\sigma)$ error, with $\sigma=1$ for $c\ne0$ and $\sigma=\min(e,1-e)$ for $c=0$. Decrease $\delta'$ below $\sigma$. Writing the metric equation with its scalar Laplacian as principal part, the remaining terms contain smooth uniformly bounded coefficients and at most quadratic products of first derivatives. The embedding $W^{2,2}\subset W^{1,p}$ for finite $p$ on the rescaled annulus and the interior $W^{2,p}$ estimate give $C^1$ control of $M-H$ by $O(r^{\delta'})$ on a smaller annulus. Differentiation and the interior estimates give the same bounds for a further derivative on the annulus and for any fixed parameter derivatives. Consequently the tangential part of the universal Chern connection on $r=\epsilon$ is
\[
 \mathcal A=\diag(i\alpha_j\dd\theta+B_j)+O(\epsilon^{\delta'}),
\]
with the same bound for the tangential derivatives in $\theta$ and along $S$ needed in~\eqref{eq:CS}. This formula is in the $D_\alpha^2$-unitary frame $r^{-\alpha_j}\sigma_j$ associated with the chosen holomorphic frame. Here $B_j$ is the Chern connection of the induced metric on the graded line. For $c\ne0$ the transition between labeled eigenframes is diagonal. For $c=0$ an adapted holomorphic transition preserves the parabolic line at $w=0$. In the corresponding $D_\alpha^2$-unitary frames its entry from the quotient to the line is $O(r^e)$, and the opposite entry is $O(r^{1-e})$. The limiting transitions are therefore those of the graded lines. The formula patches on $S$, and~\eqref{eq:CSboundary} gives $-\sum_j\alpha_j c_1(L_{pj}^0)$ at the tame pole. At an integer weight, $i:\{p\}\times S\longrightarrow X\times S$ denotes the inclusion and $\mathcal L$ the eigenline bundle whose weight changes from one to zero. The modification satisfies $[\mathcal E']=[\mathcal E]+i_*\mathcal L$. The integral of its discriminant class changes by $-\deg\mathcal L^0$, which cancels the removed integer boundary term. The real characteristic formula thus depends only on the filtered Higgs bundle.

\subsection{The complex period}
\label{subsec:complex-period}

\begin{proposition}
\label{prop:complex-period}
On the smooth rank-two moduli space with the unramified wild and tame data of Section~\ref{sec:families}, let $r_{pj}$ be the eigenvalue residues on the labeled formal summands at a wild pole and on the graded lines of the parabolic flag at a tame pole. Let $c_1(L_{pj}^0)$ denote the corresponding normalized line classes of~\eqref{eq:projectiveclasses}. Then
\begin{equation}
 [\Omega_I]=8\pi^2i\sum_{p,j}r_{pj}c_1(L_{pj}^0).
 \label{eq:complexclass}
\end{equation}
The formula also holds at a tame point with zero residue, with its parabolic flag retained.
\end{proposition}

\begin{proof}

The $L^2$ holomorphic symplectic normalization is
\begin{equation}
 \Omega_I(v_1,v_2)=2i\int_X\tr\left(
 \delta A^{0,1}_1\wedge\delta\Phi_{\mathrm{sp},2}-
 \delta A^{0,1}_2\wedge\delta\Phi_{\mathrm{sp},1}\right).
 \label{eq:complexL2}
\end{equation}
We use holomorphic representatives in local trivializations on disks around the punctures. The integrand below then vanishes on these disks, and the local Liouville form
\begin{equation}
 \Theta=2i\int_X\tr(\Phi_{\mathrm{sp}}\wedge\delta A^{0,1})
\end{equation}
satisfies $\dd_\cM\Theta=\Omega_I$. Let $g$ be the transition between two such frames and use primes in the second frame. Then $(A')^{0,1}=g^{-1}A^{0,1}g+g^{-1}\bar\partial g$ and $\Phi'_{\mathrm{sp}}=g^{-1}\Phi_{\mathrm{sp}}g$. With $\eta=g^{-1}\dd_\cM g$, the equation $\bar\partial_{A'}\Phi'_{\mathrm{sp}}=0$ gives
\begin{align}
 \Theta'-\Theta
 &=2i\int_X\tr(\Phi'_{\mathrm{sp}}\wedge\bar\partial_{A'}\eta)\notag\\
 &=-2i\int_{\partial X_\epsilon}\tr(\Phi'_{\mathrm{sp}}\eta)
 =-4\pi\sum_p\Res_p\tr(\Phi'_{\mathrm{sp}} g^{-1}\dd_\cM g).
 \label{eq:liouvilletransition}
\end{align}
The last sign again comes from orienting the inner boundary circles clockwise.

For a rank-one formal summand, the transition has the expansion
\begin{equation}
 \log g_j(w)=\log g_j(0)+\sum_{k\ge1}c_{jk}w^k.
\end{equation}
The coefficients of the higher-order poles pair with the additive cocycles $\dd_\cM c_{jk}$ for the coefficients of positive powers of $w$ in this expansion. Since the additive structure group is contractible, a smooth splitting represents each cocycle as a coboundary, so their contributions to the de Rham class vanish. At a tame point with nilpotent residue the strictly flag-lowering part has zero trace pairing with a flag-preserving infinitesimal transition. Thus only the induced residues on the parabolic graded lines enter, including when the full residue vanishes. The remaining class comes from multiplying the residue $r_j$ by $\dd\log g_j(0)$. Since $[F_{B_j}]=-2\pi i c_1(L_j)$, equation~\eqref{eq:liouvilletransition} gives~\eqref{eq:complexclass}.

The integral of the holomorphic symplectic form over the fiber vanishes. \end{proof}

On a regular genus-one spectral fiber, represent a vertical change of character by $-i\beta$, where $\beta$ is a closed real one-form with periods $(\delta\theta_e,\delta\theta_m)$. The sign follows from parallel transport for $\dd+A$. By subtracting differentials of cutoff local primitives we may take $\beta$ to vanish near the ramification points and the poles. For a horizontal change the diagonal Higgs variation is the variation $\delta\lambda$ of the spectral differential at fixed $z$. The trace on the two sheets in~\eqref{eq:complexL2} gives the mixed pairing $2\int_\Sigma\delta\lambda\wedge\beta$. The Riemann bilinear relation evaluates it as $2(\delta Z_A\,\delta\theta_m-\delta Z_B\,\delta\theta_e)$. The vertical pair vanishes on the Picard fiber. The horizontal pair also vanishes, since fixing the transported character gives the holomorphic section over a complex curve constructed in Section~\ref{sec:families}. Therefore
\begin{equation}
 \Omega_I=2(\dd Z_A\wedge\dd\theta_m-\dd Z_B\wedge\dd\theta_e).
 \label{eq:factor2}
\end{equation}

\subsection{Descent at ramified poles}

\begin{corollary}
\label{cor:ramified-periods}
The real and complex period formulas apply to the ramified families of Section~\ref{subsec:ramified} by computing on a finite Galois cover and dividing by its degree. A ramified pole contributes its pulled-back eigenline terms divided by the local degree. In $A_0,D_0,D_1$, these poles contribute zero to the complex period formula. An auxiliary unmarked branch point contributes no boundary term.
\end{corollary}

\begin{proof}

For a family with ramified poles, we compute the characteristic classes on the finite cover of Section~\ref{sec:families}, where the eigenvalue one-forms are meromorphic. The universal connection, its discriminant form, and the cutoff boundary expression are invariant under the deck action. Division by the covering degree gives the $L^2$ forms on the moduli space for the original curve. Thus the contribution of a ramified pole is the sum of the contributions from the eigenline bundles on the cover divided by the local covering degree. This prescription uses the pulled-back filtration, including the determinant extension. On a smooth Picard fiber, passing from the ordinary pullback of the spectral pushforward to its adapted eigenline lattices changes the bundles by modifications at marked sections. Indeed, the two lattices coincide off those sections and are commensurable there. Filter their finite-length difference by powers of a local spectral parameter. The graded quotients on the Picard fiber are evaluations of the spectral line at the fixed points, tensored with fixed jet factors, and hence have degree zero. The modification calculation of Section~\ref{subsec:calibration} then leaves the discriminant integral unchanged; ordinary pullback multiplies it by the covering degree. The boundary line degrees on the cover vanish on that fiber for the same reason. Thus its downstairs area is $4\pi^2$ also in the ramified cases. For $A_0,D_0,D_1$, the pulled-back Higgs residues at every ramified pole vanish, so these poles contribute zero to~\eqref{eq:complexclass}.

By Proposition~\ref{prop:ramified-quotient}, the smooth descended metric has weight zero at an auxiliary branch point away from the marked poles. There is therefore no boundary contribution at that point, and the Chern--Weil and Liouville formulas apply for the same families parametrized by surfaces as in the unramified case.
\end{proof}

\section{Comparison with the semiflat metric near infinity}
\label{sec:comparison}

On a distant regular Hitchin fiber we glue the decoupled metric to local harmonic metrics at ordinary branch points and tame poles. Chen--Li carried out this comparison for the pole data in~\cite[Sections~5--6]{CL24}. Here the monodromy about a branch point gives a lower bound for the gauge quadratic form uniformly over the Picard fiber. Closed representatives of the tangent classes then give an upper bound for the $L^2$ metric, including at the tame parameters in~\eqref{eq:tame-domain}. The real volume form determined by $\Omega_I$ fixes the product of the two relative eigenvalues, giving the reciprocal lower bound.

\subsection{General comparison and derivative estimates}
\label{subsec:spatial-comparison}

\begin{lemma}
\label{lem:comparison-bootstrap}
Let $g$ and $g_{\sfm}$ be the hyperk\"ahler and semiflat metrics on the regular locus near infinity of one of the seven families, with their common complex structure $I$ and holomorphic symplectic form $\Omega_I$. Suppose that, on compact sets of admissible pole parameters,
\begin{equation*}
 g\le(1+C\varepsilon(u))g_{\sfm},\qquad \varepsilon(u)\longrightarrow0.
\end{equation*}
Then $g\ge(1+C\varepsilon(u))^{-1}g_{\sfm}$. On the rescaled holomorphic charts described below, let $E_r$ be the difference of the Hermitian matrices of the rescaled metrics and let $\varepsilon_r$ bound $\varepsilon(u)$ on the larger chart. If the rescaled semiflat coefficients and their derivatives are uniformly bounded and the metric is uniformly positive there, then, on a smaller chart,
\begin{equation*}
 \|E_r\|_{C^k}\le C_{k,N}\varepsilon_r^{1-k/N}\qquad(N>k).
\end{equation*}
Consequently, a polynomial $C^0$ decay rate yields a positive polynomial decay rate for covariant derivatives of every order, while an exponential or stretched-exponential rate retains its form with a possibly smaller positive constant in the exponent.
\end{lemma}

\begin{proof}
The eigenvalues of $g$ relative to $g_{\sfm}$ satisfy $\lambda_1\lambda_2=1$, because the metrics are K\"ahler in the same complex structure on a complex surface and have the same volume form. The upper bound therefore implies the stated lower bound.

On a cover of the regular part of the end on which the parallel character angles lift to real functions, we use common holomorphic coordinates
\begin{equation}
 z_1=Z_A,\qquad z_2=\theta_m-\tau(z_1)\theta_e,
 \qquad\Omega_I=2\dd z_1\wedge\dd z_2.
 \label{eq:holomorphic-end-coordinates}
\end{equation}
For an ALG end, the base and lifted-fiber radii are comparable to $\rho$, and we rescale the metric by $\rho^{-2}$. The elliptic-period expansions in Section~\ref{sec:periodend} bound the rescaled semiflat coefficients and their derivatives. For an ALG$^*$ end, the radial variable and the imaginary part of the period ratio are $r=|2Z_A|$ and $V=\Ima\tau\asymp\log r$. We take a base radius $\varepsilon r$, a lifted-fiber radius $\varepsilon rV$, and divide the metric by $r^2V$. At the chart center, first choose $B\mapsto B+nA$ so that $\Rea\tau$ is bounded, then lift both angles to $[-\pi,\pi]$ by a deck translation. The identity
\begin{equation*}
 \dd z_2=\dd\theta_m-\tau\dd\theta_e-\theta_e\tau'(z_1)\dd z_1
\end{equation*}
shows that the mixed base--fiber coefficient of the rescaled metric remains bounded because $\tau'=O(r^{-1})$ and the range of the lifted angle coordinate has length $O(r)$. The derivatives of this coefficient in the rescaled coordinates are bounded as well. The $C^0$ comparison gives uniform ellipticity for the rescaled $L^2$ metric on these charts.

A closed positive $(1,1)$-form with uniformly bounded coefficients on an $I$-holomorphic coordinate ball has a local potential with bounded $C^0$ norm on a smaller ball. To bound the potential in terms of the coefficients of the form, we solve the Euclidean Poisson equation for the trace of the form with zero boundary data. The solution has bounded $W^{2,p}$ norm for every finite $p$. Its complex Hessian differs from the original form by a closed trace-free form. In complex dimension two, a real trace-free $(1,1)$-form is anti-self-dual for the Euclidean metric. Closedness therefore implies coclosedness in distributions. Weyl's lemma gives smoothness, and the interior harmonic estimates bound its derivatives by its $L^p$ norm. A local $\dd\dd_I^c$ primitive of this smooth harmonic remainder gives the bounded potential. Since $\Omega_I=2\dd z_1\wedge\dd z_2$, the common real volume form makes both metric determinants constant in~\eqref{eq:holomorphic-end-coordinates}. On these rescaled balls, the bounded potential and uniform metric bounds give a uniform Laplacian bound and a positive lower bound for this determinant. By the interior estimate of~\cite[Theorem~1.1]{Wang12}, the potential is therefore bounded in $C^{2,\beta}$ on a smaller ball for some $\beta>0$. Applying Schauder estimates to the differentiated complex Monge--Amp\`ere equation $\log\det(g_{i\bar j})=\mathrm{const}$ gives all higher interior bounds.

On nested balls the preceding estimates give $\|E_r\|_{C^0}\leq C\varepsilon_r$ and $\|E_r\|_{C^N}\leq C_N$. Interior interpolation gives, for $N>k$,
\begin{equation}
 \|E_r\|_{C^k}\leq C_{k,N}\|E_r\|_{C^0}^{1-k/N}\|E_r\|_{C^N}^{k/N}
 \leq C_{k,N}\varepsilon_r^{1-k/N}.
 \label{eq:spatialinterpolation}
\end{equation}

For a polynomial $C^0$ error and $k\ge1$, taking $N=2k$ in~\eqref{eq:spatialinterpolation} gives decay at half the original exponent for derivatives of order $k$. For $k=0$, we use the original $C^0$ bound. One positive decay exponent thus applies to covariant derivatives of every order on the moduli space, with constants depending on their order. Exponential and stretched-exponential errors remain of the same form with a smaller positive constant in the exponent. Returning to the physical metric multiplies the $k$th covariant derivative bound by $\Lambda^{-k}$, where $\Lambda\asymp\rho$ on an ALG end and $\Lambda=r\sqrt V\asymp\rho$ on an ALG$^*$ end. In the latter case logarithmic factors in the conversion of the $|u|$-decay to distance decay are absorbed by taking a slightly smaller exponent.

\end{proof}

\subsection{Families with only wild poles}

\subsubsection{Gluing and coercivity lemmas}

We first estimate the energies of trial representatives near ordinary branch points. Two further lemmas compare the glued and harmonic metrics on a punctured surface and give the coercivity needed on the sphere. These last two estimates will also be used in the families with tame poles.

\begin{lemma}[Trial representatives at the ordinary ramification points]
\label{lem:ad-disk-comparison}
Consider a regular fiber, with fixed polar data, for which the normalized spectral curve is smooth. On a simply connected set of regular values use the horizontal lift defined by the flat character of the decoupled spectral line. Around each ordinary simple zero choose nested disks with coordinate $z$ on the base curve. The eigenvalues on the two sheets are $\pm T_j\sqrt{z}\,\dd z$, where $T_j>0$ after a unitary change of phase. The disks, including a slightly larger disk, have fixed radii in these coordinates. Suppose the approximate metric $h_*$ used on the fiber is the decoupled metric $h_d$ outside the larger disks, is the exact simple-zero metric $h_j$ on the smaller disks, and is obtained on fixed annuli by multiplying the self-adjoint logarithm of $h_d^{-1}h_j$ by a fixed smooth cutoff and exponentiating relative to $h_d$. Use the same decoupled eigenline frames in the three regions, and suppose the exterior spectral differentials are admissible in the fixed-pole tangent domains. If $T_j\ge T_*>1$, then every horizontal or vertical tangent class on the fiber admits a closed representative whose squared norm for $h_*$ is at most $1+C e^{-cT_*}$ times the semiflat squared norm of the class. Here $C,c$ are uniform for a bounded number of disks of the prescribed radii.
\end{lemma}

\begin{proof}
Write $D''=\bar\partial_E+[\Phi,\cdot]$ and $D'_k=\partial_{E,k}+[\Phi^{*k},\cdot]$ on endomorphisms. Multiplication by a spectral function or form is denoted by $M$. On the unramified set this uses the eigenline splitting, and we regard $M_\nu$ for a holomorphic spectral differential as a Higgs component and $M_\xi$ for an antiholomorphic spectral differential as a Dolbeault component. The fixed unit complex factor between $\Phi$ and $\Phi_{\mathrm{sp}}$ is included in $\nu$. Normalize these differentials by the Hodge metric on the spectral curve, using the same trace and degree convention as for the semiflat metric. Trace-free forms are anti-invariant under the spectral involution. On the unit disk with $s^2=z$ they have primitives $\alpha,\beta$, vanishing at $s=0$, with $\dd\alpha=\nu$ and $\dd\beta=\xi$. The former is holomorphic and the latter antiholomorphic. Both are odd. On a smaller disk the coefficients of the forms, of the primitives, and of their derivatives of any fixed order are bounded by a constant times their global spectral norm, by the mean value inequality and Cauchy's estimates in $s$. These bounds do not use compactness of the full spectral family.

For the branch model the factor $h_d$ in a spectral generator consists of the prescribed branch power times the exponential of a regular harmonic function. Removing the real part of its holomorphic primitive by a change of generator puts it in the flat model frame. Fix the determinant and use the corresponding cyclic frame on the direct image. In this frame the simple-zero solution $h_j$ has constant determinant and agrees with $h_d$, with an $O(e^{-cT_j})$ error through any fixed number of derivatives on a fixed annulus. This is the ordinary fiducial solution. For example, in a coordinate $x$ with $Y^2=Q(x)$ one can take $z=(\frac32\int_{x_0}^x\sqrt Q\,\dd x)^{2/3}$. Then $Y=J\sqrt z$ with $J=\dd z/\dd x$; in the cyclic frame $(J^{1/2},YJ^{-1/2})$ the Higgs field has the indicated companion form. The conjugations just described are isometries between the corresponding metric frames, so these annular estimates are independent of the point of the Picard fiber. Polynomial factors in $T_j$ are absorbed by decreasing $c$.

On a core disk define
\[
 H_j=D'_{h_j}M_\alpha,\qquad
 V_j=D''S_j,\qquad S_j=(M_{\bar\beta})^{*h_j}.
\]
Here $M_\alpha$ and $M_{\bar\beta}$ extend as holomorphic matrices on the base disk by spectral multiplication. In particular $H_j,V_j$ are smooth across $z=0$. Since $h_j$ solves Hitchin's equation, $D''D'_{h_j}+D'_{h_j}D''=0$. Holomorphic spectral multiplication commutes with $\Phi$. Hence $D''H_j=D'_{h_j}H_j=D''V_j=D'_{h_j}V_j=0$. The K\"ahler identities on endomorphism-valued one-forms then give $(D'_{h_j})^*H_j=(D'')^*V_j=0$. These are the adjoint equations used in the boundary integrals below.

We recall the ordinary-branch calculation of Mochizuki~\cite[Lemma~3.42]{MochRankTwoLocal}. If $C_j$ is the nondegenerate symmetric pairing for the decoupled local model, $H_j^{1,0}$ is $C_j$-symmetric and $H_j^{0,1}$ is $C_j$-skew. The diagonal part of $H_j$ in the eigenline splitting is $M_\nu$. On the punctured disk,
\[
 H_j-M_\nu=D''q_j ,
 \qquad \|q_j\|_{C^1(\mathrm{ann})}\le C e^{-cT_j}\|\nu\|_{\mathrm{sp}}.
\]
Here $q_j$ is $C_j$-skew. The cited assertion is local to the ordinary simple-zero model, applies at the scale $T_j$ of this disk, and is used only on its fixed overlap annulus. Equivalently, the Higgs component of the difference is off diagonal in the eigenline splitting and determines $q_j$ by inverting the commutator with $\Phi$. On the annulus the eigenvalue gap is bounded below by a positive multiple of $T_j$, and the derivatives of $h_j/h_d-1$ are exponentially small. Closedness and the skew symmetry give the claimed equality for the Dolbeault component as well.

With a cutoff $\chi_j$ equal to one on the core and zero outside the larger disk, set $H=H_j$ on the core, $H=M_\nu+D''(\chi_jq_j)$ on its overlap, and $H=M_\nu$ on the exterior. This class has the prescribed base variation and is the character-horizontal lift. To see the latter statement, on the exterior identify the sheets over the fixed base and parallel-transport their flat eigenlines. This gives precisely the representative $M_\nu$ there. On a punctured spectral disk let $e_{\mathrm{fl}}$ be a multivalued parallel generator of unit norm. Its meridian character is $-1$, and $s^{-1/2}e_{\mathrm{fl}}$ is a holomorphic generator of the spectral quotient line across the ordinary branch. Trivializing $e_{\mathrm{fl}}^2$ by parallel transport identifies the square of this extension with $\cO(\mathcal R)$, where $\mathcal R$ is the moving ramification divisor. In the cyclic basis $s^{-1/2}e_{\mathrm{fl}},s^{1/2}e_{\mathrm{fl}}$, taking the two-sheet trace gives a constant antidiagonal pairing, up to an overall scalar. On the annulus this is the transported pairing because its flat trivialization was used on both sides. This is the local description of a horizontal line family in Mochizuki~\cite[Propositions~3.14 and~3.25]{Moch26}. 

Infinitesimally the Dolbeault component is skew and the Higgs component symmetric for that pairing. Indeed the skew Dolbeault component of the difference between two such local representatives can be gauged to zero by solving $\bar\partial$ on the disk in skew endomorphisms. The remaining difference in the trace-free symmetric Higgs components then has holomorphic matrix $\bigl(\begin{smallmatrix}0&a\\ b&0\end{smallmatrix}\bigr)\dd z$ in the companion frame $\Phi=T_j\bigl(\begin{smallmatrix}0&z\\1&0\end{smallmatrix}\bigr)\dd z$. Equality of the two spectral variations gives $a+zb=0$. The difference is the commutator with the holomorphic skew matrix $\diag(b,-b)/(2T_j)$. There is no skew endomorphism commuting with this companion matrix on the annulus. Thus the skew transition correction is unique, and adds no diagonal cocycle to the transported spectral line. At the wild ends under consideration the exterior identification preserves the fixed filtered lattice and peripheral holonomies. This gives the same horizontal class as the local construction of~\cite[Lemma~3.46]{MochRankTwoLocal}.

For the vertical class take the global closed representative
\[
 V=M_{\xi-\bar\partial\sum_j\chi_j\beta_j}
      +D''\!\left(\sum_j\chi_j S_j\right).
\]
The first term vanishes on the cores and represents the direct image of the same infinitesimal spectral line deformation as $\xi$; the second term is a globally defined exact change of representative. On the cores this is $V_j$. On the fixed annuli $S_j-M_{\beta_j}$ and its derivatives are exponentially small, since in the decoupled metric $(M_{\bar\beta_j})^{*h_d}=M_{\beta_j}$ and $h_j/h_d$ has the stated accuracy. Thus the corrections to $M_\nu$ and $M_\xi$ on the annuli have norm $O(e^{-cT_j})$ times the corresponding spectral norm. On the exterior the metric is the orthogonal direct image and the two representatives are direct images of their spectral forms, so both squared energies there equal the spectral integrals.

On a core disk $D$ use the Green formula for $D'_{h_j},D''$; its normalization below is the trace convention in~\cite[Lemma~3.55 and Corollary~3.56]{MochRankTwoLocal}. Local harmonicity gives
\[
 \|H_j\|_{D,h_j}^2=2i\int_{\partial D}\!\tr
           \bigl(M_\alpha(H_j^{1,0})^{*h_j}\bigr),\qquad
 \|V_j\|_{D,h_j}^2=-2i\int_{\partial D}\!\tr
           \bigl(S_j(V_j^{0,1})^{*h_j}\bigr).
\]
There is no inner boundary since the matrices on the base disk are smooth. On the outer circle one can replace $h_j,H_j^{1,0},S_j,V_j^{0,1}$ by $h_d,M_\nu,M_\beta,M_\xi$, respectively, with an exponentially small error relative to the global spectral squared norm. Taking the trace on the base is integration over the two sheets. By Stokes on the spectral disk these boundary integrals equal $2i\int_{\pi^{-1}D}\nu\wedge\bar\nu$ and $-2i\int_{\pi^{-1}D}\xi\wedge\bar\xi$. These are exactly the core terms of the two semiflat squared norms. This is the local calculation of~\cite[Lemma~3.51 and \S3.5.7]{MochRankTwoLocal}, with the fiducial metric itself on the core.

For the genus-one fibers under consideration the normalizations can also be read off from the periods. Let the holomorphic differential $\omega$ have $A$-period one and $B$-period $\tau$. The bilinear relation gives $2i\int_\Sigma\omega\wedge\bar\omega=4\Ima\tau$. Thus the horizontal squared norm with $A$-period $\dot Z_A$ is $4(\Ima\tau)|\dot Z_A|^2$. A flat character with real angle increments $(a,b)$ has harmonic real one-form $\gamma=k\omega+\bar k\bar\omega$ with periods $(a,b)$, where $2\Rea k=a$ and $2\Rea(k\tau)=b$. With the character convention of Section~\ref{sec:families}, the variation of its Dolbeault operator is $-i\gamma^{0,1}$. Its squared norm is $4(\Ima\tau)|k|^2=|b-\tau a|^2/(\Ima\tau)$. These are the horizontal and vertical semiflat coefficients for angles of period $2\pi$ and the symplectic normalization~\eqref{eq:factor2}. Summing yields the two asserted bounds.

\end{proof}

Suppose also that the determinant-normalized harmonic metric $h$ satisfies $\|\log(h_*^{-1}h)\|_\infty\le Ce^{-cT_*}$. For each closed representative its squared energy for $h$ is at most $\exp(2\|\log(h_*^{-1}h)\|_\infty)$ times its squared energy for $h_*$. Minimizing over the Dolbeault class gives the two upper diagonal bounds for $g_{L^2}$. On the moduli surfaces in this paper, to pass to the full quadratic form choose a unit horizontal and a unit vertical vector as a $g_{\sfm}$-orthonormal complex basis. The two metrics induce the same real volume form, so the matrix of $g_{L^2}$ in this basis has determinant one. If both diagonal entries are at most $1+\varepsilon$, its off-diagonal entry has modulus at most $\sqrt{2\varepsilon+\varepsilon^2}$ and each diagonal entry is at least $(1+\varepsilon)^{-1}$. Both relative eigenvalues are consequently $1+O(\sqrt\varepsilon)$. For $\varepsilon=Ce^{-cT_*}$ this is the desired exponential estimate with a smaller constant in the exponent. No parameter derivative of the correcting metric occurs in this argument.

We record the estimate that transfers the closed representatives of Lemma~\ref{lem:ad-disk-comparison} to the harmonic metric. Compare also Proposition~4.13 of~\cite{CaoExceptional26}. We give the puncture cutoff argument for the fixed determinant metric here. The argument uses the metric variation and logarithmic trace inequalities of Simpson~\cite[Lemma~3.1 and proof of Proposition~5.1]{Simpson88}.

\begin{lemma}[Comparison with a small curvature defect]\label{lem:ad-small-defect}
Let $g_X$ be a smooth conformal metric on a compact surface $X$, and let $D$ be a finite set. Suppose that $h_*$ and $h$ are Hermitian metrics on the same rank-two Higgs bundle over $X\setminus D$, with the same determinant metric, that $h$ is harmonic, and that $s=\log(h_*^{-1}h)$ is bounded for this pair. Write
\[
 Q_*(\upsilon)=\|D_{h_*}\upsilon\|_2^2+2\|[\Phi,\upsilon]\|_2^2,\qquad
 r_*=i\Lambda_{g_X}(F_{h_*}+[\Phi,\Phi^{*h_*}]).
\]
Here $D_{h_*}$ denotes the unitary Chern covariant derivative on endomorphisms. Suppose $Q_*(\upsilon)\ge c_*\|\upsilon\|_2^2$ on compactly supported smooth trace-free $h_*$-Hermitian endomorphisms. For $\eta=\|r_*\|_2$ sufficiently small,
\[
 \|s\|_\infty\le C\sqrt{\eta}.
\]
The smallness threshold and $C$ depend only on $g_X$ and $c_*$.
\end{lemma}
\begin{proof}
Put $M=\|s\|_\infty$ and $h_\theta=h_*e^{\theta s}$. Test the metric variation formula along $0\le\theta\le1$ against $\chi^2s$, where the cutoffs $\chi$ vanish near $D$, tend to one off $D$, and satisfy $\|\dd\chi\|_2\to0$. Such cutoffs are obtained by linear interpolation in the logarithm of the radius in each puncture disk. The error from differentiating $\chi$ is bounded by $CM\|\dd\chi\|_2 I_\chi^{1/2}$, where $I_\chi=\int_0^1\|\chi D''s\|_{h_\theta}^2\,\dd\theta$. Thus $I_\chi\le 2M\|r_*\|_1+CM^2\|\dd\chi\|_2^2$. Since $M<\infty$ for this pair, the metrics along the segment are mutually bounded. It follows that $s$ has finite form norm, and $\chi s\to s$ in that norm. Removing the cutoffs yields
\[
 -\int_{X\setminus D}\tr(sr_*)\,\dd V_{g_X}
   =\int_0^1\|D''s\|_{h_\theta}^2\,\dd\theta
   \ge \frac{Q_*(s)}{2(1+2M)}.
\]
For the inequality, diagonalize $s$ pointwise and use $\int_0^1e^{\theta x}\,\dd\theta\ge(1+|x|)^{-1}$ on each entry, together with $2\|D''s\|_{h_*}^2=Q_*(s)$. Coercivity and Cauchy--Schwarz give $Q_*(s)^{1/2}\le C\eta(1+M)$.

For $f=\log(\tr(e^s)/2)$, Simpson's trace inequality gives $\Delta f\le C|r_*|$ for the nonnegative scalar Laplacian. We have $0\le f\le |s|$ and $|\dd f|\le |D_{h_*}s|$. The same cutoffs extend the inequality distributionally to $X$. The compact Green kernel, shifted by a constant to be nonnegative and uniformly bounded in $L^2$ in its second variable, gives $\sup f\le C(\|f\|_1+\eta)\le C\eta(1+M)$. Here the norm of $s$ in a pointwise eigenvalue bound can be taken to be its operator norm. The two eigenvalues of $s$ are opposite, so $|s|\le C(f+\sqrt f)$. Consequently $M/(1+M)\le C(\eta+\sqrt\eta)$, proving the claim after reducing the threshold.
\end{proof}

\begin{lemma}[Coercivity on the fixed sphere]
\label{lem:observation-form-gap}
Let $g_X$ be a fixed smooth metric on the sphere, and let $A$ be a fixed annulus disjoint from the poles and the gluing regions. Suppose $h_*$ is decoupled on $A$, the two eigenlines are exchanged around its core, and $|\lambda_+-\lambda_-|_{g_X}\ge d>0$ there. Then the quadratic form $Q_*$ of Lemma~\ref{lem:ad-small-defect} satisfies $Q_*(\upsilon)\ge c\|\upsilon\|_2^2$ on compactly supported trace-free Hermitian test endomorphisms. The constant $c$ depends only on $g_X,A,d$.
\end{lemma}
\begin{proof}
On $A$ the Chern connection preserves the orthogonal splitting into diagonal and off-diagonal endomorphisms. The diagonal trace-free part is a flat real line with holonomy $-1$. On each circle of a fixed parametrization of $A$, its coefficient in a parallel frame is antiperiodic, so the angular Poincar\'e inequality controls its squared norm by its covariant derivative. For the off-diagonal part, the commutator with the Higgs field controls its norm by the eigenvalue gap. Integrating over the annulus gives $\|\upsilon\|_{L^2(A)}^2\le C Q_*(\upsilon)$. The scalar Poincar\'e inequality on the sphere, with the mean controlled by observation on $A$, gives
\[
 \|f\|_{L^2(X)}^2\le C_A\bigl(\|\dd f\|_2^2+\|f\|_{L^2(A)}^2\bigr).
\]
For a smooth test endomorphism supported away from the poles, $|\upsilon|$ extends by zero to a Lipschitz function on $X$, hence lies in $W^{1,2}(X)$. The scalar estimate and Kato's inequality $|\dd|\upsilon||\le |D_{h_*}\upsilon|$ give the asserted bound, which extends by density to the closure in the form norm.
\end{proof}

In the five purely wild families, $h_*$ is glued to $h_d$ only at ordinary simple zeros. In a normalized chart with $T_j\asymp R$, the fiducial tail and its first two derivatives give a curvature defect $O(R^N e^{-cR})$ on each collar by the fiducial estimates in Lemma~\ref{lem:ad-disk-comparison}. Local isometric frame changes preserve this bound over the Picard fiber. If changing to the fixed round background costs at most a power of $R$, then $\|r_*\|_2\le Ce^{-c'R}$. Near each pole $h_*=h_d$ is harmonic and induces the same filtered extension as $h$, so the two metrics are mutually bounded for each individual pair by~\cite[Theorem~2.21]{Moch21}. Once the annulus of Lemma~\ref{lem:observation-form-gap} is specified, Lemma~\ref{lem:ad-small-defect} makes the logarithmic ratio uniformly exponentially small. Regularity of the exterior spectral differentials in the polar tangent domains then allows Lemma~\ref{lem:ad-disk-comparison}, the determinant estimate following its proof, and Lemma~\ref{lem:comparison-bootstrap} to give the full metric and its spatial derivative estimates.

\subsubsection{The quartic family}

\begin{proposition}
\label{prop:comparison}
For the family~\eqref{eq:spectral} on the domain~\eqref{eq:domain}, the metrics $g_{L^2}$ and $g_{\sfm}$ satisfy
\begin{equation}
 |g_{L^2}-g_{\sfm}|_{g_{\sfm}}
 \le C\exp(-c|u|^{3/4})
 \label{eq:actualsf}
\end{equation}
for sufficiently large $|u|$. Among the smooth moduli spaces in this family, the constants and the threshold for $|u|$ are uniform for all $\alpha\in[0,1]$, including the endpoints $0,1$, and $(a,b)$ in a fixed compact set. The comparison holds for covariant derivatives of every order on the end, with a possibly smaller constant $c>0$ in the exponent.
\end{proposition}

\begin{proof}

Set $R=|u|^{3/4}$ and $z=u^{1/4}x$ on a sector. The spectral differential is $u^{3/4}Y\,\dd x$, where
\begin{equation}
 Y^2=x^4+A x^2+B x+1,\qquad A=au^{-1/2},\quad B=bu^{-3/4}.
 \label{eq:scaledquartic}
\end{equation}
These quartics approach $Y^2=x^4+1$ uniformly for $(a,b)$ in a compact set. Choose four mutually disjoint fixed-size disks about the simple zeros, with larger disks for the estimates of Lemma~\ref{lem:ad-disk-comparison}. Integration of the normalized square root gives simple-zero coordinates with uniform nested radii. Each local parameter is $R$, up to a phase. Choose a fixed annulus on the $x$-sphere whose inner circle encloses exactly one of the moving zeros together with its gluing collar, and which meets none of the other gluing disks.

Use the real representatives $(\alpha,1-\alpha)$, including at the endpoints. The spectral quotient line in this extension has degree one. Passing to the normalized endpoint lattice is an elementary modification of the same filtered line. Its ordinary ramification weights are all $-1/2$ and its two weights at infinity sum to one. Thus its parabolic degree is zero, and its flat metric exists for every allowed weight. In particular, up to determinant normalization, varying the weight changes its logarithm by
\begin{equation}
 2(\alpha-\alpha_0)(G_{p_+}-G_{p_-}),
 \label{eq:greenweights}
\end{equation}
where $G_p$ has singular part $\log|w|$ at $p$. We use this flat metric to define $h_d$ separately for each spectral line.

Use the round $x$-sphere. On the chosen annulus the sheets are exchanged and the eigenvalue gap is at least $cR$. At fixed $z$ the horizontal differential is $\dot u\,\dd z/(2\sqrt{q(z)})$. At infinity this differential and the conjugate differential representing a vertical class are regular. In an adapted frame the trace-free diagonal connection is $\dd$, so their cylindrical coefficients, including one cylindrical derivative, are $O(|w|)$, in the weighted domain for $0<\delta<1$ at every admissible weight. The common comparison argument gives~\eqref{eq:actualsf}.
\end{proof}

\subsubsection{Admissible weights at wild poles}

\begin{proposition}
\label{prop:purewild}
For $D_2$, throughout the smooth parameter domain for the real representatives of the parabolic weights, the difference of the Hitchin and semiflat metrics and its covariant derivatives of each order decay at a stretched exponential rate. The estimates are uniform on compact subsets of the admissible weight and pole parameters. They also hold for $D_1$ as the weight at its unramified pole varies with the ramified principal parts and the weights at infinity fixed.
\end{proposition}

\begin{proof}

We first consider $D_2$. Write $P=p_0p_\infty$ and choose a square root of $u$ on a sector. Use
\[
 x=z\sqrt u/p_0,\qquad y=\sqrt u/(p_\infty z),\qquad xy=u/P,\qquad R=|u|^{1/2}.
\]
In the $x$ chart the normalized quadratic differential $u^{-1}\lambda^2$ is
\[
 \bigl(x^{-4}+x^{-2}+2m_0u^{-1/2}x^{-3}
          +2Pm_\infty u^{-3/2}x^{-1}+P^2u^{-2}\bigr)\dd x^2 .
\]
The $y$ chart has the same formula with $m_0,m_\infty$ interchanged. On compact subsets away from zero the limiting differential in each chart is $(x^{-4}+x^{-2})\dd x^2$, with simple zeros at $i,-i$. The two charts account for all four zeros for $|u|$ sufficiently large. Their simple-zero coordinates and nested disks are uniform on the compact parameter sets of the proposition. Their model parameters are comparable to $R$.

Use the round $x$-sphere and an annulus enclosing only the zero near $i$, outside its gluing collar. The sheets are exchanged there and their gap is at least $cR$. On a $y$-collar the round metric is comparable to $|u/P|^{-2}|\dd y|^2$, costing only a power of $|u|$ in the $L^2$ curvature defect. The quotient line has degree zero and parabolic degree zero: the four branch weights sum to $-2$ and the two pairs of polar weights sum to $2$. Its flat metric therefore has the prescribed growth at both poles.

For $s^2=z^4q(z)$ one has $\lambda=s\,\dd z/z^2$ and $\partial_u\lambda=\dd z/(2s)$ at fixed $z$. This differential is regular on both sheets over zero and infinity, and a harmonic vertical differential is a multiple of its conjugate. Thus their direct images are diagonal forms in the admissible polar domains as in the quartic. Lemma~\ref{lem:ad-disk-comparison}, the determinant argument and Lemma~\ref{lem:comparison-bootstrap} give the uniform metric and derivative estimates. The identification of the ordered ALG$^*$ model by the two end periods, and the subsequent order-$2$ improvement, are given in Section~\ref{sec:periodend}. For $D_1$ the same comparison, with three ordinary zeros and one ramified pole, is proved in Section~\ref{subsec:D1-comparison}.
\end{proof}

\subsubsection{Ramified families}

The analytic quotient at a ramified pole was constructed in Proposition~\ref{prop:ramified-quotient}. We now estimate the constants in the metric comparison after rescaling the spectral equations for $A_0,D_0,D_1$.

\paragraph{The three ordinary zeros in type $A_0$.}
The rescaling~$z=u^{1/3}x$ gives
\begin{equation}
 \lambda=u^{5/6}\sqrt{x^3+1+\epsilon x}\,\dd x\,,
 \qquad \epsilon=a u^{-2/3}\,.
 \label{eq:A0-rescaling}
\end{equation}
For~$a$ in a fixed compact set and~$|u|$ sufficiently large, the discriminant~$-4\epsilon^3-27$ is bounded away from zero. The three finite zeros remain uniformly separated in the rescaled coordinate. For the normalized differential $u^{-5/6}\lambda$, the ramification at infinity, the leading term of the differential there, and the weights are independent of $\epsilon$; the coefficient of the lower-order pole term varies smoothly with $\epsilon$. The scaled spectral curves and their Jacobians range over compact smooth families.

For the norm estimates we keep a reference coordinate $z=u_0^{1/3}x$ fixed while varying $u$ near $u_0$. Set $L=u_0^{5/6}$, $R=|L|$, $b=u/u_0$, and $Y^2=x^3+au_0^{-2/3}x+b$. Then $\lambda=LY\,\dd x$ and, at $u=u_0$,
\[
 \dot\lambda=\frac{\dot u}{2}u_0^{-1/6}\frac{\dd x}{Y}.
\]
The Hodge norm of $\dd x/Y$ on these elliptic curves stays bounded above and away from zero. Thus $|\dot u|\le C|u_0|^{1/6}$ for a horizontal vector of unit semiflat norm. A unit vertical class is a bounded multiple of $\dd\bar x/\bar Y$. On the spectral disk at infinity, $x=t^{-2}$ gives
\[
 \frac{\dd x}{Y}=-2(1+au_0^{-2/3}t^4+bt^6)^{-1/2}\,\dd t
\]
after choosing the sign of $t$. The exterior spectral representatives therefore extend on the local double cover as regular diagonal one-forms with the prescribed deck action. Their coefficients and first derivatives in this cover coordinate are bounded near the puncture. In particular they belong to the invariant weighted tangent space of Section~\ref{subsec:ramified}, with a sufficiently small positive weight. In fact the induced connection on the trace-free diagonal endomorphism line in the eigenframe is $\dd$; a regular coefficient $F(t)\,\dd t$ with bounded first derivative, or its conjugate, has cylindrical coefficient $O(|t|)$ and satisfies the indicated weighted bounds for $0<\delta<1$.

The local squared energy over $|t|<r$ is $O(r^2)$ for each fixed tangent. In the leading model~\eqref{eq:ramified-local-model} with $w=t^2$, the pushforward of a leading $\dd t$ is $\frac12\bigl(\begin{smallmatrix}0&1\\w^{-1}&0\end{smallmatrix}\bigr)\dd w$. On the local cover the eigenvalue variation is regular and preserves each residue on the graded quotients. The centered vertical primitive is $O(r)$, so cutting it off on $r<|t|<2r$ also costs $O(r^2)$ in squared norm.

For an explicit parallel reference take the spectral quotient line $\cO(O)$, where $O$ is the point at infinity, with squared norm $|1|^2=|Y|^{-1}$ in its meromorphic frame (up to a scalar determinant normalization). At each finite branch point this has the quotient-line factor $|s|^{-1}$; at infinity the holomorphic generator $t^{-1}$ has squared norm comparable to $|t|$. Its unitary holonomy on the punctured curve is the sign determined by the winding of $Y$ and hence is locally constant under transport. Tensoring by a unitary flat degree-zero line gives the other points of the fiber.

Glue on fixed nested simple-zero disks about the three finite roots. All three disks and their collars lie in $|x|<2$. Hence the cover over $2<|x|<3$ exchanges its sheets and the gap there is at least $cR$. The round metric of the $x$-sphere is uniformly comparable to the chart metrics on these collars. The reference line above has the prescribed growth at infinity, and~\eqref{eq:A0-rescaling} and the calculation in $t$ verify the branch-disk and polar hypotheses. The comparison gives
\begin{equation}
 |g_a-g_{\sfm,a}|_{g_{\sfm,a}}
 \le C e^{-c|u|^{5/6}}\,.
 \label{eq:A0-actual-comparison}
\end{equation}
This estimate holds uniformly at every point of each sufficiently distant Hitchin fiber and for~$a$ in compact sets, including sets containing zero.

\paragraph{Two ramified poles in type $D_0$.}
Fix a compact set $K\Subset(\C^*)^2$ of leading coefficients, write
$P=p_0p_\infty$, and put $R=|u|^{1/2}$. The two coordinates
\begin{equation*}
 w_0=uz/p_0,\qquad w_\infty=u/(p_\infty z),\qquad
 w_0w_\infty=u^2/P
\end{equation*}
give the same expression at their respective poles:
\begin{equation*}
 \lambda^2=u\bigl(w^{-3}+w^{-2}+(P/u^2)w^{-1}\bigr)\,\dd w^2 .
\end{equation*}
Each chart contains a simple zero at $w=-1+O_K(|u|^{-2})$. In the original coordinate the two zeros are $-(p_0/u)(1+O_K(|u|^{-2}))$ and $-(u/p_\infty)(1+O_K(|u|^{-2}))$. Each of
the two zero disks, together with a slightly larger disk, can be
chosen with radii independent of the parameters and disjoint from the
pole in its chart. These disks account for all finite simple zeros
away from the poles. The normalized spectral curve
$y^2=z(p_0+uz+p_\infty z^2)$ is smooth and connected when
$u^2\ne4P$, which holds uniformly for sufficiently large $|u|$.
The moving coordinates are used at a fixed reference value $u_0$: in forming an infinitesimal class we hold these charts at their $u_0$ values. Both scaled quadratic differentials, after division by $u$, range on
the disks in a compact family with one simple zero. Their fiducial
parameters are comparable to $R$.

Here the trivial spectral quotient line has the explicit squared norm $|1|^2=|z/y|$ up to a scalar. At the ordinary branch points this has the factor $|s|^{-1}$; at both ramified poles it is comparable to $|t|$ in the local uniformizing coordinate $t$. Its unitary holonomy is given by the winding parity of $z/y$, so it is a parallel character reference on the regular family. Twisting by unitary flat lines sweeps out the degree-zero Picard fiber.

Use the round metric on the $w_0$-sphere. Its coefficients on the remote gluing collar are comparable to $|u|^{-4}|\dd w_\infty|^2$, uniformly on $K$, so $\|r_*\|_2\le C R^N e^{-cR}$. A fixed annulus surrounding $0$ inside $0<|w_0|<1/2$ is decoupled, exchanges sheets about the ramified pole and has gap at least $cR$. The explicit metric on the quotient line gives the prescribed extensions on both polar covers. Lemmas~\ref{lem:observation-form-gap} and~\ref{lem:ad-small-defect} therefore give $\|\log(h_*^{-1}h)\|_\infty\le Ce^{-c'R}$. The slightly larger normalized zero disks control the derivatives and centered primitives of the spectral test forms by their global Hodge norms by the interior estimate in Lemma~\ref{lem:ad-disk-comparison}; degeneration of the full elliptic curve does not affect this local estimate.

To check the behavior at the poles, fix the original base coordinate
$z$ when differentiating $u$. Then
\begin{equation*}
 \lambda=y\,\frac{\dd z}{z^2},\qquad
 \partial_u\lambda=\frac{\dd z}{2y}.
\end{equation*}
This differential is regular in the local spectral coordinates
$z=t^2$ at zero and $z=t^{-2}$ at infinity. The harmonic vertical
form is also regular there. The corresponding diagonal trace-free
one-forms on each local double cover are invariant under the
combined deck and eigenline action. Their coefficients in the
cylindrical coframe are $O(|t|)$ for each fixed tangent vector and hence lie in
the small-weight deformation domain of Section~\ref{subsec:ramified}.
The primitive of a regular spectral form is $O(|t|)$ after centering at the puncture. Its boundary pairing with a bounded regular differential along $|t|=r$ is $O(r^2)$. The two trial energies and the determinant estimate following Lemma~\ref{lem:ad-disk-comparison} now give, uniformly on $K$,
\begin{equation}
 |g-g_{\sfm}|_{g_{\sfm}}\leq C\exp(-c|u|^\sigma)
 \quad\text{for some }\sigma>0.
 \label{eq:D0-actual-comparison}
\end{equation}
As a check, pullback by $z=s^2$ gives $4(p_0s^{-4}+us^{-2}+p_\infty)\dd s^2$, a $D_2$ model with $u'=4u$, zero residues and weights $1/2$, on the excluded walls of~\eqref{eq:D2-domain}.

\subsubsection{Weights in type $D_1$}
\label{subsec:D1-comparison}

Write $\lambda_0=pb$ and choose a branch of $\sqrt u$. We use coordinates
\begin{equation}
 x=z\sqrt u/p,\qquad y=u/(bz),\qquad xy=u^{3/2}/\lambda_0,\qquad R=|u|^{1/2}.
 \label{eq:D1-two-charts}
\end{equation}
The normalized differentials are
\begin{align}
 u^{-1}\lambda^2&=(x^{-4}+x^{-2}+2mu^{-1/2}x^{-3}
                         +\lambda_0u^{-3/2}x^{-1})\,\dd x^2,\label{eq:D1-inner-quad}\\
 u^{-1}\lambda^2&=(y^{-3}+y^{-2}+2\lambda_0mu^{-2}y^{-1}
                         +\lambda_0^2u^{-3})\,\dd y^2 .\label{eq:D1-outer-quad}
\end{align}
The limiting differentials have two simple zeros $i,-i$ in the $x$ chart and one simple zero $-1$ in the $y$ chart. By uniform convergence on fixed neighborhoods of these points, there are three simple-zero disks with fixed nested radii and parameters comparable to $R$. These exhaust the ordinary zeros. The bounds are uniform on compact subsets of the admissible $(\alpha,m,\lambda_0)$ with $\lambda_0\ne0$.

Glue the three exact fiducial metrics to the decoupled metric and use the round metric of the $x$-sphere. A fixed annulus about the zero approaching $i$ encloses no other ramification and lies outside the gluing collar. The two sheets are exchanged around the core, and the eigenvalue gap is at least $cR$. Lemma~\ref{lem:observation-form-gap} gives
\begin{equation}
 Q_*(\upsilon)\ge c_K\|\upsilon\|_2^2.\label{eq:D1-direct-gap}
\end{equation} On the $y$-collar the round metric and the global curvature defect satisfy
\begin{align}
 g_X&\asymp |u^{3/2}/\lambda_0|^{-2}|\dd y|^2,\label{eq:D1-round-collar}\\
 \|r_*\|_{L^2(g_X)}&\le C|u|^N e^{-cR}.\label{eq:D1-direct-defect}
\end{align} The spectral quotient line has degree zero. Its three ordinary ramification weights contribute $-3/2$, the two weights over zero sum to one and the weight at the ramified pole is $1/2$. The decoupled flat metric realizes this filtered extension. It is harmonic near both poles and is mutually bounded there with the actual metric by~\cite[Theorem~2.21]{Moch21}, using the local double cover at infinity. Lemma~\ref{lem:ad-small-defect} therefore gives, uniformly over the Picard fibers and the indicated pole parameters,
\begin{equation}
 \|\log(h_*^{-1}h)\|_\infty\le C e^{-c'R}.\label{eq:D1-direct-ratio}
\end{equation}

Let $s^2=p^2+2pmz+uz^2+bz^3$. At fixed $z$ we have $\partial_u\lambda=\dd z/(2s)$, which is regular over $z=0$. At infinity, putting $z=t^{-2}$ gives
\begin{equation}
 \frac{\dd z}{s}=-2(b+ut^2+2pmt^4+p^2t^6)^{-1/2}\dd t.\label{eq:D1-polar-tangent}
\end{equation} A harmonic vertical form is a multiple of its conjugate. Over infinity both direct images are diagonal and invariant for the simultaneous deck and eigenline action. Their cylindrical coefficients are $O(|t|)$ for each fixed tangent vector and belong to the small-weight domain. The corresponding assertion over zero follows in the unramified eigenframe. Lemma~\ref{lem:ad-disk-comparison} and the determinant argument after it now give
\begin{equation}
 |g-g_{\sfm}|_{g_{\sfm}}\le Ce^{-c|u|^{1/2}}.
 \label{eq:D1-actual-comparison}
\end{equation}
The estimates for all spatial derivatives follow from Lemma~\ref{lem:comparison-bootstrap}. The refinement to order two in appropriate coordinates uses~\cite[Theorem~1.10]{CVZH23} after identification of the semiflat model in Section~\ref{sec:periodend}.

\subsection{Families with tame poles}

For $A_2,D_3$, we compare the local harmonic and decoupled metrics at a tame pole in the same holomorphic eigenframe. Curvature decay controls the noncommuting part of the connection, while extension of a holomorphic one-form to the genus-zero normalization of the local spectral curve determines its commuting part. Closed representatives then bound the norms of all horizontal and vertical tangent classes.

\subsubsection{A normal form at a tame pole}

The rescaled coordinate for $A_2$ is $w=-uz$. For $D_3$ the rescaled coordinates are $w=uz$ at zero and $w=-u(z-1)$ at one. On $|w|<4R$, where $R=|u|^{1/10}$,
\begin{equation*}
 \lambda^2=\frac{F(w)}{w^2}\dd w^2,
 \quad F=c^2-w+E(w),\quad E(0)=0,
\end{equation*}
and, with the corresponding derivative bounds in the rescaled coordinate on the disk,
\begin{equation*}
 E=O(|u|^{-2}w^2)\quad(A_2),\qquad
 E=O(|u|^{-1}(|w|+|w|^2))\quad(D_3).
\end{equation*}
The linear term in the second estimate is necessary: at zero, $2pm/z$ contributes $2pmw/u$ to $F$.

The homotopy $F_s=c^2-w+sE$ factors as $(r_s-w)A_s$, with $A_s$ holomorphic and uniformly close to one. A vector field $v_s=wg_s\partial_w$ satisfies the Moser equation for $F_s\dd w^2/w^2$ if
\begin{equation*}
 2F_sg_s'+F_s'g_s=-E/w.
\end{equation*}
For $g_s=A_s^{-1/2}\widetilde g_s$ and $h_s=-E/(w\sqrt{A_s})$, the equation becomes
\begin{equation*}
 2(r_s-w)\widetilde g_s'-\widetilde g_s=h_s,
 \qquad
 \widetilde g_s(w)=-\frac12\int_0^1
 h_s(r_s+t(w-r_s))t^{-1/2}\,\dd t.
\end{equation*}
Differentiating the formula verifies the equation. For derivatives with respect to $w$, the linear map $h_s^{(k)}\mapsto \widetilde g_s^{(k)}$ has operator norm at most $1/(2k+1)$ in the sup norm, uniformly as $c$ tends to zero. The flow fixes $w=0$. For admissible weights and prescribed coefficients in compact sets, its difference from the identity is small on $|w|<2R$, along with its parameter derivatives of every fixed finite order. The inverse time-one map gives the coordinate for the exact model
\begin{equation}
 \widehat\Phi_{\mathrm{sp},c}=
 \begin{pmatrix}0&c^2-\xi\\1&0\end{pmatrix}\frac{\dd\xi}{\xi}.
 \label{eq:weakmodel}
\end{equation}
When $c=0$, the simple-pole coefficient of $\lambda^2$ is $u$ in $A_2$, and $2pm-u$ at zero and $u$ at one in $D_3$. These coefficients are nonzero on sufficiently distant fibers, uniformly on compact parameter sets. A zero residue would make the Higgs field holomorphic, so the repeated eigenvalue at such a pole has nonzero nilpotent residue. At an ordinary branch point a zero matrix would give a multiple discriminant zero. Thus the Higgs field is cyclic throughout the disk. Its spectral line bundle on the smooth spectral disk admits a holomorphic generator, which together with its spectral-coordinate multiple gives a cyclic frame for $w\Phi_{\mathrm{sp}}/\dd w$. With $\xi=\psi(w)$ set $g=\diag(k,k^{-1})$, where $k^2=w\psi'/\psi$ and $k(0)=1$. Then $g(\psi^*\widehat\Phi_{\mathrm{sp},c})g^{-1}$ has lower-left entry $\dd w/w$ and upper-right entry $F(w)\dd w/w$. These gauge transformations satisfy the same estimates, including the bounds for any fixed finite number of parameter derivatives. We henceforth denote the coordinate $\xi$ in the normal form by $w$.

The odd spectral pole order at infinity excludes an invariant line. The good filtered correspondence~\cite[Theorem~4.1]{Moch21} gives the unique normalized harmonic metric, and the inverse in Section~\ref{sec:families} with the local double cover at infinity gives smooth dependence with bounded derivatives on compact parameter charts. At $c=0$ the nilpotent perturbation has positive decay exponent, while at $\alpha=0,1/2$ the term $16|c|^2>0$ supplies the indicial gap.

\subsubsection{The local model with rational spectral curve}
\label{subsec:globalmodel}

\begin{proposition}
\label{prop:globalmodel}
For every compact $K\subset\mathcal P_{\mathrm{tame}}$, the harmonic metric of \eqref{eq:weakmodel} differs from the decoupled harmonic metric induced from the spectral line bundle by $O_K(|w|^{-N})$ at infinity for every $N$. Both metrics are taken with the same holomorphic identification of the eigenline bundles and the same normalized determinant metric. The estimate also holds after any fixed finite number of derivatives on the local curve and in the parameters.
\end{proposition}

\begin{proof}
The normalization of the compactified model spectral curve is $\Sigma_c=\{\zeta^2=c^2-w\}\simeq\PP^1_\zeta$, also when $c=0$. In a holomorphic eigenframe, $A$ denotes the Chern connection matrix and $\Phi_{\mathrm{sp}}=\diag(\lambda,-\lambda)$. For any $h$, the commutator with the diagonal matrix $\Phi_{\mathrm{sp}}$ has zero diagonal. Hitchin's equation gives
\begin{equation*}
 \bar\partial A_{++}=\bar\partial A_{--}=0.
\end{equation*}
Under a change of holomorphic frame in an eigenline bundle, $A_{++}$ changes by the logarithmic derivative of the transition function. Here $A_{++}$ is the projection of the Chern connection along the holomorphic eigenspace decomposition. The Chern connection of the restricted eigenline bundle metric is a different connection unless the eigenspaces are orthogonal. The difference between the projected connection $A_{++}$ and the corresponding decoupled connection is a holomorphic one-form on the punctured compact spectral curve; the same assertion holds for $A_{--}$.

We verify removability at every point. At a branch point over an unmarked point of the base, we set $z=\zeta^2$ and
\begin{equation*}
 E=\begin{pmatrix}\zeta&-\zeta\\1&1\end{pmatrix},\qquad
 A^{\mathrm{eig}}_{++}
 =\frac12\tr A^{\mathrm{base}}
  +\frac12\left(\frac{A^{\mathrm{base}}_{12}}\zeta+\zeta A^{\mathrm{base}}_{21}\right)
  +\frac{\dd \zeta}{2\zeta}.
\end{equation*}
The original connection is smooth and $\dd z=2\zeta\dd \zeta$, so only the last term has a pole, of residue $1/2$. At an unramified lift of a tame pole with semisimple residue, the adapted expansion is $A_{jj}=\alpha_j\dd z/z+O(r^{-1+\delta})\dd z$. The holomorphic remainder in this expansion extends across the pole, including at the weights $\alpha=0,1/2$.

At the tame point $c=0$, where the spectral cover ramifies, the remainder estimate in the adapted unitary frame for $0<\alpha<1/2$ gives in a holomorphic base frame
\begin{equation*}
 A^{\mathrm{base}}_{12}=O(r^{-2\alpha+\delta})\dd z,
 \qquad A^{\mathrm{base}}_{21}=O(r^{2\alpha-2+\delta})\dd z.
\end{equation*}
The middle term in the preceding eigenframe formula has size
\begin{equation*}
 O\bigl(|\zeta|^{-4\alpha+2\delta}
             +|\zeta|^{4\alpha-2+2\delta}\bigr)\dd \zeta
       =o(|\zeta|^{-2})\dd \zeta.
\end{equation*}
It is holomorphic and anti-invariant under the deck involution. Its only possible pole is simple, but $\dd \zeta/\zeta$ is deck-invariant, so its residue vanishes and this anti-invariant remainder extends holomorphically across $\zeta=0$. Half the determinant connection contributes $\dd \zeta/\zeta$, and the eigenvector derivative contributes $\dd \zeta/(2\zeta)$; thus the residue is $3/2$.

The explicit decoupled connection gives the same residues. In the affine eigenvector frame, with determinant normalized by $\det h=|w|^2$, the decoupled diagonal entry is
\begin{equation*}
 k(\zeta)=2|\zeta|\,|\zeta-c|^{2\alpha}|\zeta+c|^{2(1-\alpha)},
\end{equation*}
with the other entry $k(-\zeta)$. Its connection is
\begin{equation*}
 A_{d,++}=\left(\frac1{2\zeta}+\frac\alpha{\zeta-c}
                         +\frac{1-\alpha}{\zeta+c}\right)\dd \zeta.
\end{equation*}
At $c=0$ it is $3\dd \zeta/(2\zeta)$. The eigenline bundle is a line subbundle of the pulled-back bundle. Its branch contribution to the residue is $1/2$. The twist by the ramification divisor accounts for the difference from the weight $-1/2$ of the spectral line bundle in the quotient presentation of the spectral correspondence~\cite{BNR89}.  At infinity put $q=q_\infty=\zeta^{-1}$. On its ramified cover the off-diagonal perturbations of the flat connection in unitary Coulomb gauge vanish to every order by \cite[Lemmas~4.6 and~5.3]{BB04}. This applies to both the unitary and self-adjoint parts of the connection. The diagonal Higgs coefficients commute with their adjoints, so Hitchin's equation gives $|F_h|_h=O(|q|^N)$ for every $N$.

We transfer the bound to the projected connection in the holomorphic eigenframe. On small disks in $r/2<|q|<2r$ multiply the two eigenvectors by a branch of $q^{3/2}$. The Gram matrix of these holomorphic sections and its inverse are bounded, because the quotient of the adapted metrics is bounded and $|q|^3 k(\pm\zeta)$ is bounded above and below. Applying the integrated Bochner identity to each section with a cutoff supported on a disk of radius comparable to $r$ bounds its squared covariant derivative integral on a smaller disk by $C(r^{-2}+\|F_h\|_\infty)O(r^2)=O(1)$. Hence the rescaled frame's Chern matrix has bounded $L^2$ norm there. Its $\bar\partial$ is the curvature and is bounded in this frame. The interior Cauchy estimate on a disk of radius comparable to $r$ bounds each coefficient by $C(r^{-1}\|A\|_2+r\|F_h\|_\infty)=O(r^{-1})$. The same bound holds for the explicit decoupled connection. Their projected difference is independent of the common frame factor, so it has at most a simple pole at $q=0$. Its deck transform is its negative, since the two determinant connections agree. A simple polar term $\dd q/q$ is invariant and hence its coefficient vanishes.

Consequently
\begin{equation*}
 A_{++}-A_{d,++}\in H^0(\PP^1,K_{\PP^1})=0.
\end{equation*}
Thus the projected connection agrees everywhere with the decoupled connection, identifying the commuting part of the Chern connection.

Now use the intrinsic rapid curvature estimate above in the holomorphic eigenframe. Write $H=\left(\begin{smallmatrix}a&b\\\bar b&d\end{smallmatrix}\right)$ and $\Delta=ad-|b|^2$. Direct matrix multiplication gives
\begin{equation*}
 [\diag(\lambda,-\lambda),
 H^{-1}\diag(\bar\lambda,-\bar\lambda)H]_{12}
       =\frac{4|\lambda|^2db}{\Delta}.
\end{equation*}
Since the metrics are adapted to the filtration, changes of frame and their inverses have at most polynomial growth at the puncture. The leading eigenvalue separation and the rapid curvature bound therefore give $|b|/\sqrt{ad}=O(|q|^N)$ for every $N$. The corresponding first derivative estimate follows in the dyadically rescaled eigenframe above. In that frame $a,d$ are uniformly comparable to one, $A=O(r^{-1})$, and $F=\bar\partial A$ is bounded. From $\partial H=HA$ and $\bar\partial H=A^\dagger H$, with $\dagger$ denoting conjugate transpose of matrices, we obtain $|\Delta_{\mathrm{eucl}}H|\le C r^{-2}$. On an interior ball of radius $\ell\ll r$, the Poisson gradient estimate bounds $|\nabla b|$ by $C(\ell^{-1}\sup|b|+\ell r^{-2})$. Given $N\ge0$, take $\ell=r^{N+3}$ and use the rapid bound for $b$ at order $2N+3$. This gives $|\nabla b|=O(r^N)$ and in particular $|\partial\bar b|/\sqrt{ad}=O(r^N)$ in that frame. Removing the common scalar frame factor adds at most $O(r^{-1})|b|/\sqrt{ad}$ to the normalized estimate.

The ratios $a_\pm=h_{\pm\pm}/h_{d,\pm\pm}$ of the diagonal metric coefficients are positive and bounded, with bounded inverses. The identity $\partial\log h_{++}=A_{++}+(h_{+-}/h_{++})A_{-+}$ in the common rescaled frame, together with the $O(|q|^{-1})$ connection bound just proved, gives for every $N$,
\begin{equation*}
 A_{++}=\partial\log(h_{d,++}a_+)+O(|q_\infty|^N)\dd q_\infty,
 \qquad \partial\log a_+=O(|q_\infty|^N)\dd q_\infty.
\end{equation*}
Integration gives $a_\pm=C_\pm+O(|q_\infty|^N)$ for every $N$. The determinant condition implies $C_+C_-=1$, while invariance under the deck involution near the same ramified puncture implies $C_+=C_-$. Positivity therefore forces $C_+=C_-=1$, so the relative metric difference vanishes to every order at $q_\infty=0$.

For every decay order and every finite set of derivative orders, we obtain uniform constants on $K$ by using the smooth local gauge slice in the weighted Sobolev spaces for the tame and irregular poles on finitely many slightly larger parameter charts. The principal part of the model Higgs field at infinity is independent of $(c,\alpha)$, so we can iterate the off-diagonal weighted estimates and differentiate the equations to bound the parameter derivatives. Equivalently, rapid $C^0$ decay and polynomial bounds for higher parameter derivatives, from the weighted implicit function theorem, give rapid decay of the lower-order derivatives by interpolation. For example, if $|f|\le C_N|q_\infty|^N$ for every $N$ and $|\partial_p^m f|\le C|q_\infty|^{-M_m}$ on the larger chart, the interior interpolation inequality gives
\begin{equation*}
 |\partial_p^k f|
 \le C\|f\|_\infty^{1-k/m}
          \|\partial_p^m f\|_\infty^{k/m}+C\|f\|_\infty,
 \quad k<m.
\end{equation*}
For fixed $k,m,M_m$, we can take $N$ large enough to obtain any desired decay power. The interpolation estimate holds after any fixed number of differentiations in the local coordinate on the cover. At $\alpha=0,1/2$ the same Landau--Kolmogorov interpolation inequality on a one-sided parameter interval bounds the derivatives up to its endpoint. The relation $|q_\infty|\asymp|w|^{-1/2}$ gives the stated decay at infinity.
\end{proof}

\begin{remark}
In this local example $w=\zeta^2$, and the Higgs field and metric are
\begin{equation*}
 \Phi_{\mathrm{sp}}=\begin{pmatrix}0&w^{-1}\\w^{-2}&0\end{pmatrix}\dd w,
 \qquad h_0=\diag(|w|^{1/2},|w|^{3/2}).
\end{equation*}
In the eigenframe $(\zeta,1),(-\zeta,1)$ the metric is $2|\zeta|^3\Id$. Multiplying its two entries by $e^{\pm2\Rea(b\zeta)}$ preserves the determinant and gives a flat harmonic metric descending under the deck involution. For $b\ne0$ its relative difference from $h_0$ has leading order $|w|^{1/2}$. Thus curvature decay alone cannot imply the metric conclusion of~\cite[Lemma~3.6]{CL24}. For the rational model, the projected connections agree because $H^0(\PP^1,K_{\PP^1})=0$. On a positive-genus normalization their difference can lie in a nonzero space of holomorphic one-forms.
\end{remark}

\subsubsection{Closed representatives of horizontal tangent classes at a tame pole}

The differential of the total Higgs deformation complex is $D''=\bar\partial_E+[\Phi_{\mathrm{sp}},\cdot]$, with the graded commutator on endomorphism-valued forms. We write $a_h=\delta u/u$. For the model~\eqref{eq:weakmodel}, the Lie derivative along $w\partial_w$ is
\begin{equation*}
 \mathcal L_{w\partial_w}\widehat\Phi_{\mathrm{sp},c}=-E_{12}\dd w.
\end{equation*}
This deformation preserves the residue. Near $c=0$ use the holomorphic frame $(e_1,(-c,1))$, whose second vector spans the parabolic line at the pole. In this frame $E_{12}$ is unchanged and the residue is $\left(\begin{smallmatrix}c&0\\1&-c\end{smallmatrix}\right)$. Its lower off-diagonal term has unitary norm $O(r^{-2\alpha})$ for all these $c$, so the positive-weight estimates require no division by $c$. On a compact chart with $\alpha\ge\alpha_*>0$ and $\alpha<1/2$ bounded away from $1/2$, the adapted estimates in Section~\ref{subsec:mixed} give $|E_{12}|^2\le C r^{-2+4\alpha_*}$ for $r=|w|<1$. The estimate is uniformly integrable through $c=0$. The remaining parameter charts can be chosen with $|c|$ bounded away from zero. Choose a fixed small disk with $|w|<|c|^2/2$, on which the two eigenvalues are distinct. In a holomorphic eigenframe for the varying family use the two labeled eigenvectors to trivialize the moving lines, with their chosen filtered lattices. The inner tangent then consists of the regular commuting eigenvalue variation $\varphi_{\sfm}$ below, bounded in the adapted metric also at endpoint weights. On a collar avoiding the pole its difference from the cyclic-frame deformation $-E_{12}\dd w$ is the commutator with the transition $S$ of~\eqref{eq:weak-commutator}. Cut off this transition on the collar, where the eigenvalue gap and frame changes are uniformly controlled. This constructs the required representative in the fixed tangent domain. On the outer part of the model, Proposition~\ref{prop:globalmodel} and the explicit decoupled metric give $|E_{12}|^2=O(|w|^{-1})$. The contribution to the squared energy over $|w|<2R$ satisfies
\begin{equation}
 E_{h,\mathrm{tame}}\le C|a_h|^2(1+R).
 \label{eq:weakenergy}
\end{equation}

On the annulus $R<|w|<2R$, we use
\begin{align}
 \varphi_{\sfm}&=-\frac{w}{2(c^2-w)}\widehat\Phi_{\mathrm{sp},c},
 \qquad S=\frac{w}{4(c^2-w)}\diag(1,-1),\notag\\
 [\widehat\Phi_{\mathrm{sp},c},S]&=-E_{12}\dd w-\varphi_{\sfm}.
 \label{eq:weak-commutator}
\end{align}
On this annulus $S$ and $R\dd S$ are bounded. The squared cutoff energies are at most $C|a_h|^2$ for the connection component and $C|a_h|^2R$ for the Higgs component. To differentiate the normal form at fixed original base coordinate, write $w=\psi_u(w_0)$, with $w_0$ the initial rescaled coordinate. The induced vector field for $u\partial_u$ has coefficient
\[
 (u\partial_u\psi_u)(w_0)+w_0\psi_u'(w_0)=wg(w).
\]
Here we differentiate in the varying holomorphic cyclic frame of the normal form. In the original frame, differentiating $\diag(k,k^{-1})$ adds an infinitesimal change of frame. The parameter estimates for $k^2=w_0\psi_u'/\psi_u$ allow its gauge term to be cut off on the matching annulus. The first summand differentiates the normalizing map at fixed $w_0$, and the second pushes forward the dilation of $w_0$. With $\epsilon=|u|^{-j}$, $j=2$ for $A_2$ and $j=1$ for $D_3$, the integral formula for the normalizing flow gives $\psi_u(w_0)/w_0=1+O(\epsilon R)$ on the larger disk. This quotient is holomorphic in $w_0$ and in $u$ on a sector, with the same estimate when $u$ varies in a disk of radius proportional to $|u|$ and $w_0$ is fixed. Cauchy's estimates on nested disks give $g=1+O(\epsilon R)$ and $g'=O(\epsilon)$. Writing $B=\left(\begin{smallmatrix}0&c^2-w\\1&0\end{smallmatrix}\right)$, the Lie derivative adds the commuting term $g'B\dd w$. Its squared energy, including the factor $a_h$, is $O(|a_h|^2\epsilon^2R^3)\le C|a_h|^2R$.

We identify the tangent class and estimate the full transition on the overlap. At fixed $c,\alpha$, extend the transported unitary spectral character over the total punctured spectral family above a contractible neighborhood of the chosen regular base value as in Section~\ref{sec:families}. In the exact coordinate $w$, the three marked spectral coordinates $0,c,-c$ are fixed. The punctured local family is a product after restricting to a common disk in $\zeta$. The meridians generate its first homology, so its unitary flat line is isometrically isomorphic, over the total local family, to the explicit model with the same meridians. When $c=0$ there is one meridian with the sum of the three logarithmic lifts. The flat isomorphism extends holomorphically over the marked sections in the chosen extensions, since their logarithmic factors agree. The associated isomorphism identifies the corresponding spectral quotients and hence the cyclic direct images. A constant phase can be chosen to respect the determinant. In these flat frames the isomorphism is constant throughout the horizontal variation. Its derivative in another frame contributes precisely that change of frame to the deformation cocycle.

More explicitly, put
\[
 m(\zeta)=\zeta^{1/2}(\zeta-c)^\alpha(\zeta+c)^{1-\alpha},
 \qquad E=\begin{pmatrix}\zeta&-\zeta\\1&1\end{pmatrix}.
\]
Up to constant factors, parallel eigenvectors are the columns of $P=E\diag(m(\zeta)^{-1},m(-\zeta)^{-1})$. Branch choices change the columns by constant phases. With $\dot w=a_hwg$ as above, set $X\zeta=-a_hwg/(2\zeta)$. The trace-free change of the parallel frame relative to the cyclic frame is
\[
 -\big((\dot P)P^{-1}\big)^0
 =a_hgS+\frac{X\log m(\zeta)-X\log m(-\zeta)}{2\zeta}B
 =a_hgS+\frac{(2\alpha-1)a_hgc}{2\zeta^2}B.
\]
The first summand produces $a_hg$ times~\eqref{eq:weak-commutator}. The second, denoted by $S_{\parallel}$, is holomorphic on the collar and commutes with the Higgs field. Its eigenvalues $\pm(2\alpha-1)a_hgc/(2\zeta)$ are $O_K(|a_h|R^{-1/2})$ in the orthogonal eigenframe of the decoupled metric. For a cutoff on the scale $R$, $D''(\chi S_{\parallel})=((\bar\partial\chi)S_{\parallel},0)$. The bound $|\dd\chi|=O(R^{-1})$ and the Euclidean collar area $O(R^2)$ give squared energy $O_K(|a_h|^2/R)$ for this term alone, also in the uniformly comparable model metric. It vanishes at $c=0$. The holomorphic core variation is $a_h(-gE_{12}+g'B)\dd w$. The fixed-disk modification near a semisimple pole applies to the $gE_{12}$ term, while $g'B\dd w$ is regular in the adapted model. Thus the cutoff of the full transition glues the core representative to the exterior character-horizontal representative with the energy bounds above.

\subsubsection{Vertical classes near branch points}
\label{subsec:vertical-classes}

Write $Y=\eta$ in $A_2$, where $\eta^2$ is the numerator in~\eqref{eq:familyA2}, and $Y=y$ in $D_3$, where $y^2=P_u(z)$ is the quartic obtained by clearing the denominator in~\eqref{eq:familyD3}. On each tame disk retain the original scaled coordinate $z=z_i+\varepsilon_i w/u$. Here $(z_i,\varepsilon_i)=(0,-1)$ in $A_2$ and $(0,1),(1,-1)$ in $D_3$. Put $F_i(w)=Y^2$. Then $F_i(0)=c_i^2$, where $c_i=c$ in $A_2$, and on $|w|<4R$
\[
 F_i'(w)=-1+
 \begin{cases}
 O_K(R|u|^{-2}) & (A_2),\\
 O_K((1+R)|u|^{-1}) & (D_3).
 \end{cases}
\]
On a fixed disk containing all $c_i^2$ in its interior, $F_i(w)\to c_i^2-w$ uniformly on $K$. Its nearby simple root $r_i$ is zero when $c_i=0$. Integration of $F_i'$ from $r_i$ to $w$ gives $F_i(w)=(r_i-w)B_i(w)$, with $B_i-1$ bounded by the corresponding derivative error above. In the spectral coordinate $s^2=r_i-w$ we have exactly
\[
 \frac{\dd z}{Y}=-\frac{2\varepsilon_i}{u\sqrt{B_i(r_i-s^2)}}\dd s .
\]

Let $\omega$ be the holomorphic differential with period one on the chosen primitive $A$-cycle, and put $\eta_v=\bar\omega$. In $A_2$ the normalizing multiplier of $\dd z/\eta$ is $O_K(|u|^{1/3})$ by the period scaling in Section~\ref{sec:periodend}. In $D_3$ take one lift of a fixed circle enclosing the two bounded branch values and excluding the two far ones. The lifts give opposite primitive classes, as seen by intersecting them with the lift of an arc from a bounded to a far branch value. On this circle $P_u(z)/u\to z(z-1)$ uniformly. Evaluating the limiting integral on a large circle gives $\sqrt u\oint_A\dd z/y\to\pm2\pi i$.

Thus $\eta_v=O_K(d_u)\dd\bar s$, where $d_u=|u|^{-2/3}$ in $A_2$ and $d_u=|u|^{-1/2}$ in $D_3$. Its primitive $f$ vanishing at $s=0$ has size $O_K(d_u\sqrt R)$ over $|w|<2R$ by radial integration in the slightly larger spectral disk. A cutoff $\chi$ equal to one on $|w|\le R$ and supported in $|w|<2R$ has derivative $O(R^{-1/2})$ in $s$. The form $\eta_v$ and its correction $\bar\partial(\chi f)$ each have squared energy $O_K(d_u^2R)$ there. Their difference vanishes near the branch and all lifts of the pole, and represents the same spectral Dolbeault class. Its direct image on the overlap is taken on the unramified annulus. The exact normal-form coordinate has uniformly comparable collar sizes and derivatives. For a vertical vector $b_v\eta_v$, the full semiflat squared norm is a fixed positive multiple of $|b_v|^2\Ima\tau$ by the Riemann bilinear relations. For these cycles, $\Ima\tau\asymp_K1$ in $A_2$ and $\Ima\tau\asymp_K\log|u|$ in $D_3$ by Section~\ref{sec:periodend}.

\subsubsection{Simple zeros away from the poles}

For $r=|z|$, we use the normalized model
\begin{equation*}
 \Phi_{\mathrm{sp},T_{\mathrm f}}=T_{\mathrm f}\begin{pmatrix}0&z\\1&0\end{pmatrix}\dd z,
 \qquad h_{T_{\mathrm f}}=\diag(r^{-1/2}e^{-\ell_{T_{\mathrm f}}(r)},r^{1/2}e^{\ell_{T_{\mathrm f}}(r)}).
\end{equation*}
Here $\ell_{T_{\mathrm f}}(r)=\psi_1(8T_{\mathrm f}r^{3/2}/3)\ge0$, where $\psi_1$ is the solution of the ordinary differential equation specified in~\cite[Definition~3.5(i), equation~(24)]{CL24}; our ordered frame is the reverse of the one used there. The argument of this solution has order one when the radius is comparable to $T_{\mathrm f}^{-2/3}$. We choose the radius $r_0=T_{\mathrm f}^{-1/3}$, with $T_{\mathrm f}=|u|^{2/3}$ for $A_2$ and $T_{\mathrm f}=|u|^{1/2}$ for $D_3$. A local horizontal class has a holomorphic representative
\begin{equation*}
 v_{\mathrm{hol}}=(0,\varphi_{\mathrm{hol}}),\qquad
 \varphi_{\mathrm{hol}}
 =T_{\mathrm f}\begin{pmatrix}0&f(z)\\0&0\end{pmatrix}\dd z,\qquad
 |f|+|f'|\le C|a_h|.
\end{equation*}
To obtain the bound on $f$, fix a base value $u_0$ and freeze the change of coordinates while differentiating. In $x=z/u_0^{1/3}$ in $A_2$ and $x=pz/\sqrt{u_0}$ in $D_3$ the variations of the quadratic differentials at $u_0$ are respectively
\[
 a_hu_0^{4/3}\frac{\dd x^2}{x},\qquad
 a_hu_0\frac{\dd x^2}{x(x-p/\sqrt{u_0})}.
\]
The ordinary roots under consideration converge to the nonzero simple roots $x^3=-1$ and $x^2=-1$, respectively. In particular the two far $D_3$ roots satisfy $z=(\pm i+o(1))\sqrt{u_0}/p$. On fixed larger disks there the displayed coefficients divided by $T_{\mathrm f}^2$ are bounded by $C|a_h|$. The normal-form chart at $u_0$ and its inverse have uniformly bounded derivatives on smaller fixed disks by simplicity of these limiting roots. Expressing the quadratic variation in this frozen chart gives exactly $T_{\mathrm f}^2 f(z)\dd z^2$. Cauchy's estimate also bounds $f'$.

On the punctured disk the Higgs field component of the semiflat variation and an infinitesimal complex gauge transformation satisfy
\begin{align}
 \varphi_{\sfm}&=\frac{f}{2z}\Phi_{\mathrm{sp},T_{\mathrm f}},
 \qquad S=\diag(-f/(4z),f/(4z)),\notag\\
 [\Phi_{\mathrm{sp},T_{\mathrm f}},S]&=\varphi_{\mathrm{hol}}-\varphi_{\sfm}.
 \label{eq:zero-commutator}
\end{align}
With the transition functions of the remaining flat unitary line bundle held constant as at the tame pole, the local semiflat horizontal class is represented by
\begin{equation*}
 v_{\sfm}=v_{\mathrm{hol}}-D''S
 =(-\bar\partial_E S,\varphi_{\sfm})=(0,\varphi_{\sfm}).
\end{equation*}
The first component records the variation of $\bar\partial_E$ in this holomorphic frame, where $S$ is holomorphic on the punctured disk. The unitary Chern connection also depends on the varying Hermitian metric.

We use a cutoff $\chi$ equal to one on $r\le r_0$ and zero on $r\ge2r_0$. On the annulus $r_0<r<2r_0$,
\begin{equation*}
 D''(\chi S)=((\bar\partial\chi)S,\chi[\Phi_{\mathrm{sp},T_{\mathrm f}},S]).
\end{equation*}
The local representatives $v_{\mathrm{hol}}$ on the inner disk, $v_{\sfm}+D''(\chi S)$ on the annulus, and $v_{\sfm}$ outside agree on their overlaps and define a closed representative. Near zero this representative equals $v_{\mathrm{hol}}$ and is smooth. Since $h_{11}/h_{22}\le r^{-1}$, the contribution of the Higgs field component to the squared norm on the replacement region is bounded by $CT_{\mathrm f}^2|a_h|^2r_0$. The diagonal cutoff contributes at most $C|a_h|^2r_0^{-2}$ to the integral of the squared norm of the connection component. The integral $E_{h,\mathrm{zero}}$ of the squared norm of this representative on $r<2r_0$ therefore satisfies
\begin{equation}
 E_{h,\mathrm{zero}}\le C|a_h|^2(T_{\mathrm f}^2r_0+r_0^{-2}).
 \label{eq:zeronorm}
\end{equation}
A vertical deformation of the unitary flat line bundle is represented by $b_v\eta_v$, with $\eta_v$ normalized as in Section~\ref{subsec:vertical-classes}. Its local primitive has size at most $C|b_v|\sqrt{r_0}$. Applying the cutoff before pushing forward from the spectral curve contributes at most $C|b_v|^2r_0$ to the squared norm of the representative. For an explicit bound on the ordinary-branch disks, set $Y=\eta/u_0^{2/3}$ in $A_2$ and $Y=py/u_0$ in $D_3$, using the frozen $x$ coordinates above. The limiting polynomials for $Y^2$ are $x^4+x$ and $x^2(x^2+1)$ respectively. On a fixed disk about any of the indicated nonzero simple roots the actual polynomial is $(x-x_*)U(x)$ with $U,U^{-1}$ uniformly bounded and $x_*$ its actual root. Thus $\dd x/Y$ is bounded in $s^2=x-x_*$. This is $u_0^{1/3}\dd z/\eta$ and $\sqrt{u_0}\dd z/y$, respectively, so the period normalizations of Section~\ref{subsec:vertical-classes} bound both normalized forms in $s$. This proves the primitive and cutoff estimates uniformly through the Picard fibers. For either complex coefficient $b_v$, the full semiflat vertical squared norm is a fixed multiple of $|b_v|^2\Ima\tau$, with $\Ima\tau$ bounded below in these cycle conventions. Thus the bounds on these disks control all real angular directions, including their linear combinations.

The horizontal semiflat squared norm is comparable to $T_{\mathrm f}^2|a_h|^2\Ima\tau$. By~\eqref{eq:zeronorm}, the fraction of this norm used on an ordinary-branch disk is at most
\[
 \frac{C_K}{\Ima\tau}(r_0+T_{\mathrm f}^{-2}r_0^{-2})
 =\frac{C_K}{\Ima\tau}(T_{\mathrm f}^{-1/3}+T_{\mathrm f}^{-4/3}).
\]
For a vertical vector the corresponding fraction is $O_K(r_0/\Ima\tau)$.
The error in replacing the fiducial solution by the decoupled metric on this annulus is $O(e^{-cT_{\mathrm f}r_0^{3/2}})=O(e^{-cT_{\mathrm f}^{1/2}})$.

\subsubsection{Gluing and comparison of tangent norms}

Write $h_{\mathrm{sp}}$ for the global decoupled metric and $h_{\mathrm{mod}},h_d$ for the harmonic and decoupled metrics of the rational local model. In a common affine eigenline frame over the tame spectral disk, the logarithm of $h_{\mathrm{sp},++}/h_{d,++}$ is harmonic across $\zeta=0,\pm c$, since the prescribed factors cancel, including their combined factor when $c=0$. With the common determinant normalization it is odd under $\zeta\mapsto-\zeta$. Write the ratio as $e^{-2\Rea f(\zeta)}$, choosing the unique holomorphic odd $f$ with $f(0)=0$. The matrix
\[
 G=\cosh f(\zeta)\,\Id-\frac{\sinh f(\zeta)}{\zeta}B
\]
descends holomorphically over the whole base disk because $\cosh f$ and $(\sinh f)/\zeta$ are even and regular at $\zeta=0$. Its eigenvalues $e^{\mp f}$ have product one. It commutes with the Higgs field, preserves the parabolic line (also at $c=0$), and satisfies $G^*h_d=h_{\mathrm{sp}}$. This is the fiberwise matching also obtained from the flat isometry used for the horizontal deformation above. Its scalar phase has no effect on the metrics or the trace-free deformations.

In each tame collar interpolate, relative to $h_d$, by multiplying $\log(h_d^{-1}h_{\mathrm{mod}})$ by a cutoff in $|w|/R$ and exponentiating. Pull this entire interpolation back by $G$. Its Chern curvature and its Higgs commutator are the conjugates of those before pullback, with the same pointwise norms in the respective metrics. Proposition~\ref{prop:globalmodel} gives decay faster than every inverse power of $R$ for the relative logarithm and its first two derivatives on this collar. In this model frame differentiation of the cutoff, the connection coefficients, and the Higgs coefficients introduces at most powers of $R$. At each ordinary branch point away from the tame disks insert the fiducial metric using the same logarithmic interpolation to its decoupled metric on $r_0<r<2r_0$. Its collar error and its first two derivatives are bounded by a polynomial in $T_{\mathrm f},r_0^{-1}$ times $e^{-cT_{\mathrm f}r_0^{3/2}}$. These constructions are isometric in the matched frames, independently of the point in the angle torus. Denote the glued metric by $h_*$. It has the prescribed determinant and agrees with the exact decoupled metric near the wild pole at infinity.

On a sector choose $x=z/u^{1/3}$ for $A_2$ and $x=pz/\sqrt u$ for $D_3$, and use a fixed smooth round metric $g_X$ on the $x$-sphere. The inverse coordinate lengths of the tame collars in $x$ are of order $|u|^{4/3}/R$ and $|u|^{3/2}/R$, respectively. The ordinary-branch normal-form coordinates and their inverses have uniformly bounded derivatives on fixed small disks in $x$, since the limiting roots used there are simple. Conversion to the $L^2(g_X)$ norm of the curvature contraction thus costs at most a fixed power of $|u|$. Since $R=|u|^{1/10}$ and $T_{\mathrm f}r_0^{3/2}=T_{\mathrm f}^{1/2}$, for every $M$
\[
 \|i\Lambda_{g_X}(F_{h_*}+[\Phi_{\mathrm{sp}},\Phi_{\mathrm{sp}}^{*h_*}])\|_2
 \le C_{K,M}|u|^{-M}.
\]

For a uniform spectral gap, the normalized quadratic differentials on compact subsets away from $x=0,\infty$ have limiting coefficients $x^2+x^{-1}$ and $1+x^{-2}$ respectively, up to an overall unit phase, after division by $T_{\mathrm f}^2$. Choose a fixed small annulus surrounding exactly one simple zero of the limiting coefficient, separated from every other zero and all poles. On this annulus $h_*=h_{\mathrm{sp}}$, the sheets are exchanged around the core, and the eigenvalue gap is at least $cT_{\mathrm f}$ in $g_X$. Lemma~\ref{lem:observation-form-gap} gives a global lower bound independent of the admissible parameters in $K$ and the point of the angle torus. At each marked point both $h_*$ and the actual harmonic metric $h$ are harmonic locally and induce the same good filtration. Thus their ratio and its inverse are bounded for each fixed filtered Higgs bundle by \cite[Theorem~2.21]{Moch21}, also at a strongly parabolic point. Lemma~\ref{lem:ad-small-defect}, with the defect estimate at twice any desired order, gives
\[
 \|\log(h_*^{-1}h)\|_\infty\le C_{K,M}|u|^{-M}
 \qquad\text{for every }M.
\]
This estimate transfers the squared tangent energies from $h_*$ to $h$ with a factor $1+O_{K,M}(|u|^{-M})$.

The horizontal representative is obtained from the parallel character family on the exterior and the cyclic core representatives. Their transitions are the holomorphic gauges calculated above. At an ordinary branch point the only logarithmic factor of the eigenline in the branch coordinate $s$ is $s^{1/2}$. Its logarithmic derivatives on the two sheets agree, so its trace-free commuting correction vanishes. Replacing these transitions by their cutoffs changes only the representative of the same deformation cocycle. For vertical classes the cutoff is the Dolbeault coboundary of the local primitive on the spectral curve, regular also at $c=0$.

The harmonic representative minimizes the $L^2$ norm in its class. Divide the entire local squared energies of the closed representatives by their respective horizontal or vertical semiflat squared norms. At tame poles,~\eqref{eq:weakenergy} and the vertical estimates give $O_K(|u|^{-4/3+1/10})$ in $A_2$ and $O_K(|u|^{-9/10}/\log|u|)$ in $D_3$. At the other simple zeros the bounds are $O_K(|u|^{-2/9})$ and $O_K(|u|^{-1/6}/\log|u|)$, respectively. Outside the replacement regions the semiflat horizontal and vertical parts are orthogonal. By the Cauchy--Schwarz inequality, the absolute value of the mixed term inside these regions is at most the product of the norms of the two representatives there. Thus
\begin{equation*}
 g_{L^2}\le(1+C|u|^{-1/10})g_{\sfm}.
\end{equation*}

Applying Lemma~\ref{lem:comparison-bootstrap}, we obtain
\begin{equation}
 (1+C|u|^{-1/10})^{-1}g_{\sfm}\leq g_{L^2}\leq(1+C|u|^{-1/10})g_{\sfm}
 \label{eq:weak-comparison}
\end{equation}
for $A_2,D_3$, uniformly on compact subsets of admissible weights and prescribed principal parts.

These estimates prove Theorem~\ref{thm:comparison}. For the ALG$^*$ families, after the spectral periods have identified the ordered model triple in Section~\ref{sec:periodend}, the sharper decay estimate in~\cite[Theorem~1.10]{CVZH23} gives coordinates at infinity in which the metric has order-$2$ decay to its model, with the same $\nu,\kappa_0,L$. The three ALG metrics at $a=0$ require the separate arguments in Section~\ref{subsec:centers}.

\section{Periods on integral two-cycles in the seven families}
\label{sec:integral}

\subsection{The intersection form and transport in homology}
\label{subsec:topology-transport}

For a rational elliptic surface $Y$ with section $O$ and fiber $D$ at infinity, we write $M=Y\setminus D$. We determine the fiber types from the vanishing orders of the coefficients and discriminant in a minimal Weierstrass equation using~\cite[Table~2]{CV24}. For the fibers at infinity considered below, $D$ is a tree of rational curves, so $H^1(D;\Z)=0$. The compactly supported cohomology sequence therefore injects $H_2(M;\Z)$ into the free abelian group $H_2(Y;\Z)$, with image orthogonal to every boundary component. Since $O\cdot F=1$, the fiber class is primitive. Lee--Lin computed integral bases for these complements in~\cite[Sections~3.1.5, 3.1.6, 3.1.8, and~3.1.9]{LL23}. The diffeomorphisms in~\cite[Theorem~4.1]{CV24} preserve the fiber class; after transporting those bases, we find that the intersection form has radical $\Z F$ and that its quotient lattices are
\begin{equation}
 \begin{array}{c|c|c|c}
 \text{Family}&D&H_2(M;\Z)/\Z F&|\det|\\\hline
 A_1&\mathrm{III}^*&\text{negative }A_1&2\\
 A_2&\mathrm{IV}^*&\text{negative }A_2&3\\
 D_2&I_2^*&\text{negative }A_1\oplus A_1&4\\
 D_3&I_1^*&\text{negative }A_3&4
 \end{array}\,.
 \label{eq:quotient-lattices}
\end{equation}

By adjunction, $C^2=-2$ for each smooth compact rational curve embedded in the Higgs moduli space, since the moduli space is a holomorphic symplectic surface. Given the compactification and its section, we use the intersections of the curves to determine whether their classes generate the quotient lattice. If their images span a full-rank sublattice $\Lambda'\subset\Lambda=H_2(M;\Z)/\Z F$, then
\begin{equation}
 |\det\Lambda'|=[\Lambda:\Lambda']^2|\det\Lambda|\,.
 \label{eq:lattice-index}
\end{equation}
Equality of the Gram determinants therefore gives an integral basis after adjoining $F$.

For the triangular families below, the characteristic classes come from the following calculation. For integers $n_1,n_2$ and $\Gamma\simeq\PP^1$, a bundle of the form
\begin{equation*}
 \mathcal E=\cO_{\PP^1_z}(-n_1)\oplus
 \bigl(\cO_{\PP^1_z}(-n_2)\boxtimes\cO_\Gamma(1)\bigr)
\end{equation*}
has summand classes $-n_1X$ and $-n_2X+H$, where $X,H$ are the positive point classes of the two factors. Since $X^2=H^2=0$ and $\int_{\PP^1\times\Gamma}XH=1$, we obtain
\begin{align}
 c_1(\mathcal E)&=-(n_1+n_2)X+H,\qquad
 c_2(\mathcal E)=-n_1XH,\notag\\
 \Delta(\mathcal E)&=-n_1XH-\frac14\bigl(-2(n_1+n_2)XH\bigr)
 =\frac{n_2-n_1}{2}XH\,.
 \label{eq:discriminant}
\end{align}
At a puncture where an eigenline bundle is a graph over the first summand, its degree on $\Gamma$ is zero, while the second summand has degree one. After subtracting half the determinant degree, these degrees are $(-1/2,1/2)$.

For a quartic $P(z)=Az^4+Bz^3+Cz^2+Dz+E$, we use the binary quartic invariants
\begin{align*}
 I(P)&=12AE-3BD+C^2,\\
 J(P)&=72ACE+9BCD-27AD^2-27B^2E-2C^3,
\end{align*}
so that $\operatorname{disc}_z(P)=(4I(P)^3-J(P)^2)/27$.
\begin{lemma}
\label{lem:marking-transport}
We consider a smooth family of Hitchin moduli spaces of complex dimension two equipped with identifications of the spectral fibrations at infinity that preserve the base coordinate $u$. We assume that along every compact parameter path there is a disk $|u|\le R$ containing all finite discriminant values of the spectral family, that the inverse image of this disk is proper over the path, and that the projections of this inverse image and its boundary to the path are smooth submersions. Under these hypotheses the given identifications of the fibrations at infinity extend to a smooth family of diffeomorphisms transporting the surfaces along the path. If the parameter domain is simply connected, the induced map on integral homology depends only on the initial and final parameter values.
\end{lemma}

\begin{proof}
Along a given compact parameter path, we choose $R$ large enough that an annulus around $|u|=R$ contains only regular values of the Hitchin map throughout the path. Identifying the fibrations above this annulus gives a lift of the tangent vector to the parameter path that keeps $u$ constant and is tangent to the boundary. A partition of unity extends the horizontal lift over the inverse image of the disk without changing its projection to the parameter path or its values on a smaller neighborhood of the boundary. Properness gives a flow for the full parameter interval on the compact part. The flow agrees with the prescribed transport on the overlap, so the interior and exterior maps glue. The integral homology groups form a local system, which has trivial monodromy on a simply connected parameter domain.
\end{proof}

We use the following local properness argument in the four families whose integral homology contains classes of self-intersection minus two. The compact spectral curves are reduced and locally planar. For a reduced spectral curve, coherent sheaves correspond to modules over its finite spectral algebra. The prescribed characteristic polynomial gives rank one on each component; a torsion-free spectral sheaf pushes forward to a vector bundle on the smooth base curve.
We check properness locally on the product of the parameter space and the Hitchin base. Over an integral spectral curve every rank-one torsion-free sheaf is stable. At a reducible curve the two components $C_+,C_-$ are graphs over $\PP^1$. All compact spectral curves in these four families are connected. Indeed, in the split case the two graphs meet along the zero divisor of their difference, a nonzero section of $\cO_{\PP^1}(2)$. Label them by the two distinct points at one unramified wild pole. For a rank-one torsion-free spectral sheaf $T$, let $T_\pm$ be its restriction to $C_\pm$ modulo torsion. The kernel $K_\pm$ of $T\to T_\mp$ pushes forward to the saturated invariant line for $C_\pm$, since the quotient is a line bundle on the other graph. Write $k_\pm=\deg(\pi_*K_\pm)$, where $\pi$ is the spectral covering. Then $\chi(T_\mp)=\chi(T)-k_\pm-1$.

At a pole with distinct labeled eigenvalues, assign to each sheet the weight of its eigenline. At a tame point with zero residue on the graded lines, temporarily assign the two weights to the two sheets in either complementary order. Write $w_\pm$ for the respective sums. The checks in the family calculations below show that every weight sum induced on an invariant line by a compatible flag gives the same strict integer bound on $k_\pm$ as this assignment. This includes the case when the two saturated lines have the same fiber at a tame point. We have $w_++w_-=-\deg E$ and $\chi(T)=2+\deg E$ for the chosen extensions. Thus
\[
 \chi(T_\mp)>1-w_\mp\quad\Longleftrightarrow\quad k_\pm+w_\pm<0 .
\]
On the admissible splitting loci the thresholds are nonintegral. Choose rational $q_\pm$ in the same open unit intervals as $1-w_\pm$, with $q_++q_-=\chi(T)$. These thresholds apply throughout a neighborhood of the chosen spectral curve, including any nearby splitting allowed by a vanishing tame residue. To realize this polarization on the relative curve, take the disjoint smooth sections $P_\pm$ above the chosen wild pole and an integer $N>0$ clearing the denominators. The bundle
\[
 \cO(-Nq_+P_+-Nq_-P_-)\oplus\cO^{N-1}
\]
has $-\deg_{C_\pm}/N=q_\pm$. The relative compactified Jacobian for this general polarization is proper by the relative theorem of Esteves used in~\cite[Remark~2.16 and proof of Theorem~5.4]{MRV17}. At a tame pole $p$, choose the line of the parabolic flag in the universal projective bundle $\PP(\mathcal E|_p)$. The equations requiring the residue to preserve this line and act on it by its prescribed scalar define a closed incidence subspace. Projectivization makes this construction independent of the local lift of the universal bundle. The integer bounds just established identify this closed subspace with the stable parabolic moduli space, also when the residue vanishes and the flag varies. Properness is local on the base, so restriction to a compact parameter path and the closed Hitchin disk gives the properness required above. The smooth local gauge slices of Section~\ref{sec:families} make the projection to the parameter path a submersion. Choosing $R$ beyond all finite discriminant values makes the boundary projection a submersion. On the exterior use fixed $u$ and the continued parallel angles of Section~\ref{sec:periodend} for the prescribed identifications in Lemma~\ref{lem:marking-transport}.

\subsection{The \texorpdfstring{$A_0$}{A0} family}
\label{sec:A0}

For~\eqref{eq:familyA0}, the ramified pole at infinity has order four and weights $(1/4,3/4)$. For the degree-$-1$ extension with the prescribed parabolic determinant, a spectral sheaf has degree one. The degree-one Abel map $Q\mapsto\mathcal I_Q^\vee$ identifies this relative compactified $\Pic^1$ with the Weierstrass family $y^2=z^3+az+u$, by~\cite[Fact~7.2]{MRV17} and the relative construction in~\cite[Sections~2.1--2.3]{EGK00}. For degree zero we tensor by $\cO(-O)$, where $O$ is the smooth section at infinity; these constructions commute with base change. The total surface is smooth because the derivative of its defining equation with respect to $u$ is nonzero. The spectral quadratic differential has odd pole order at infinity, so there is no invariant line subbundle. The hypotheses of Proposition~\ref{prop:ramified-quotient} therefore hold for every $a\in\C$.

The minimal elliptic compactification has a fiber of type $\mathrm{II}^*$ at infinity. Removing this fiber gives a surface with $H_2=\Z F$ by~\cite[Section~3.1.4]{LL23}, or equivalently by the compactly supported cohomology sequence for the complement of the tree of rational components. The periods on $F$ therefore determine the periods on all compact two-cycles:
\begin{equation}
 H_2(\cM_a;\Z)=\Z F,\qquad
 \int_F(\omega_I,\omega_J,\omega_K)=(4\pi^2,0,0).
 \label{eq:A0-periods}
\end{equation}
\subsection{The \texorpdfstring{$A_1$}{A1} family}
\label{sec:A1}

For $A_1$, with $w=z^{-1}$,
\begin{equation}
 \lambda_+=-\left(w^{-4}+\frac a2w^{-2}+\frac b2w^{-1}+O(1)\right)\dd w.
 \label{eq:polar}
\end{equation}
The weights $(\alpha,1-\alpha)$ correspond to the formal eigenline bundles in this order. The underlying holomorphic vector bundle in this extension has degree $-1$ for $0<\alpha<1$, and its parabolic degree is zero. At $\alpha=0$ or $1$, elementary modification gives weights $(0,0)$ and a holomorphic bundle of degree zero. We denote the moduli space by $\cM_{\alpha,b,a}$ and use the domain
\begin{equation}
 \mathcal P_{A_1}=\{(\alpha,b,a)\in[0,1]\times\C^2:(\alpha,b)\ne(0,0),(1,0)\}.
 \label{eq:domain}
\end{equation}
The omitted pairs are precisely the walls where the period vector of a root vanishes in this marking.

At $a=b=0$ two compact curves generate the second integral homology group. We construct them in the central fiber and transport their classes over~\eqref{eq:domain}.

\subsubsection{Two spheres in the type $\mathrm{III}$ fiber}

At $a=b=u=0$ and $0<\alpha<1$, we use the family of Higgs bundles
\begin{equation}
 E=\cO\oplus\cO(-1),\qquad
 \Phi_{\mathrm{sp},[v_0:v_1]}=
 \begin{pmatrix}z^2&0\\v_0z+v_1&-z^2\end{pmatrix}\dd z,
 \qquad [v_0:v_1]\in\PP^1.
 \label{eq:sphere}
\end{equation}
On nonzero $(v_0,v_1)$, the diagonal automorphisms act through simultaneous scalar multiplication of both entries. The invariant second summand has parabolic degree $-\alpha<0$. The other saturated invariant eigenline bundle has degree $-2$ for $v_1\ne0$ and degree $-1$ for $v_1=0$, so its parabolic degree is also negative. Thus~\eqref{eq:sphere} gives a compact sphere $\Gamma_+$ of stable Higgs bundles. The family with the two diagonal entries interchanged gives $\Gamma_-$.

The two spheres exhaust the central fiber: any splitting has the form $\cO(k)\oplus\cO(-1-k)$ with $k\ge0$, and $k\ge1$ would force the lower-left entry of the Higgs field to belong to $H^0(\cO(1-2k))=0$, making $\cO(k)$ an unstable invariant line subbundle. Thus $k=0$, and stability requires the lower-left entry to be a nonzero polynomial of degree at most one. In the chart where this polynomial is $z-p$, an upper-triangular gauge transformation makes the diagonal entry a constant $a_0$, and the determinant equation is
\begin{equation}
 a_0^2=p^4+ap^2+bp+u.
 \label{eq:surfacechart}
\end{equation}
The charts with constant lower-left entry add the two points at infinity. At $a=b=u=0$, the branches $a_0=\pm p^2$ are rational curves meeting with intersection multiplicity two. This is the type $\mathrm{III}$ fiber, and
\begin{equation}
 F=\Gamma_++\Gamma_-,\qquad
 \Gamma_+^2=\Gamma_-^2=-2,\qquad \Gamma_+\cdot\Gamma_-=2.
 \label{eq:typeIII}
\end{equation}

The Hitchin fibration at these parameter values has the section
\begin{equation*}
 \Phi_{\mathrm{sp},u}=
 \begin{pmatrix}z^2&u\\1&-z^2\end{pmatrix}\dd z
 \quad\text{on }\cO\oplus\cO(-1).
\end{equation*}
At $u=0$ the section passes through the stable point $[v_0:v_1]=[0:1]$ of~\eqref{eq:sphere}. For $u\ne0$ the spectral curve is irreducible. The Jacobian has Weierstrass equation $Y^2=X^3-4uX$ with minimal vanishing orders $(3,\ge5,9)$ at infinity, so its minimal rational elliptic compactification adds a fiber of type $\mathrm{III}^*$. The class $F$ is primitive, and~\eqref{eq:quotient-lattices} gives a quotient lattice whose intersection matrix has determinant of absolute value two. Since $\Gamma_+^2=-2$, equation~\eqref{eq:lattice-index} gives
\begin{equation}
 H_2=\Z\Gamma_+\oplus\Z\Gamma_-,\qquad
 H_2=\Z\Gamma\oplus\Z F\quad\text{with }\Gamma=\Gamma_+.
 \label{eq:integralbasis}
\end{equation}

A universal bundle for the family parametrized by~\eqref{eq:sphere} is
\begin{equation}
 \mathcal E=\cO\oplus
 \bigl(\cO_{\PP^1_z}(-1)\boxtimes\cO_{\Gamma_+}(1)\bigr).
 \label{eq:universalsphere}
\end{equation}
For~\eqref{eq:universalsphere}, $(n_1,n_2)=(0,1)$, so $\int_{\PP^1\times\Gamma_+}\Delta(\mathcal E)=1/2$. At infinity the eigenline $L_+$ is the graph of $v_0/2$ over the first summand, and the eigenline $L_-$ is the second summand. Thus
\begin{equation}
 \sum_j\alpha_j\deg L_j^0
 =-\frac\alpha2+\frac{1-\alpha}{2}=\frac12-\alpha.
 \label{eq:weightedspheredegree}
\end{equation}
Proposition~\ref{prop:CW}, together with $F=\Gamma_++\Gamma_-$ and~\eqref{eq:fiberperiod}, gives the $\omega_I$-areas $4\pi^2\alpha$ and $4\pi^2(1-\alpha)$, respectively.

Over a smooth spectral curve, the whole $\Pic^1(\Sigma_u)$ is the Hitchin fiber because $\Pic^1(\PP^1)$ is a point, giving the angle lattice of~\eqref{eq:fiberperiod}. The determinant in~\eqref{eq:universalsphere} has degree one on $\Gamma_+$ and hence admits no square root there. We therefore use the projective classes~\eqref{eq:projectiveclasses}.

At infinity the two spectral branches are smooth and distinct. For $b\ne0$ the nearby curves are irreducible. At $b=0$ the strict inequality $0<\alpha<1$ makes the stability thresholds locally constant. Section~\ref{subsec:topology-transport} thus applies, including the integer weights checked below. The domain~\eqref{eq:domain} contracts by moving $\alpha$ linearly to $1/2$ before contracting $(a,b)$ to zero. Lemma~\ref{lem:marking-transport} gives a path-independent marking.

The real period on $\Gamma=\Gamma_+$ is
\begin{equation}
 \int_\Gamma\omega_I=4\pi^2\alpha.
 \label{eq:realperiod}
\end{equation}
The residues of the two eigenvalue one-forms at infinity are $r_+=-b/2$ and $r_-=b/2$. Pairing the normalized first Chern classes of the two eigenline bundles with the curve gives $(-1/2,1/2)$, and hence
\begin{equation}
 \int_\Gamma\Omega_I
 =8\pi^2i\left((-b/2)(-1/2)+(b/2)(1/2)\right)=4\pi^2ib.
 \label{eq:complexperiod}
\end{equation}

\subsubsection{Equal and integer weights}

For $0<\alpha<1$, an invariant line subbundle has parabolic degree equal to an integer plus either $\alpha$ or $1-\alpha$, so its parabolic degree cannot be zero. This includes $\alpha=1/2$: the regular semisimple leading coefficient still labels the formal eigenline bundles even though the two weights coincide. At $\alpha=0$ or $1$ and $b\ne0$, the elementary modification gives weights $(0,0)$ and degree zero. A reducible polystable Higgs bundle would have a global meromorphic eigenvalue one-form on $\PP^1$ whose only pole is at infinity and whose residue is $\pm b/2$, contrary to the residue theorem. The moduli spaces therefore have no reducible polystable points at any parameter in~\eqref{eq:domain}. By the Biquard--Boalch correspondence and Theorem~5.4 of~\cite{BB04}, these moduli spaces are smooth and complete there.

Lemma~\ref{lem:weight-domains} gives a common gauge domain and smooth dependence on the weights, including at equal and integer weights. At an integer endpoint the elementary modification of Section~\ref{subsec:calibration} preserves the period formula on the normalized extension.

At the two excluded endpoints with $b=0$ a root has zero period and reducible polystable Higgs bundles occur.

Since $\Gamma^2=-2$ and $F^2=\Gamma\cdot F=0$, the roots of the lattice~\eqref{eq:integralbasis} are
\begin{equation}
 \{\delta\in H_2:\delta^2=-2\}
 =\{\pm\Gamma+kF:k\in\Z\}.
 \label{eq:roots}
\end{equation}
The two embedded spheres of the type $\mathrm{III}$ fiber give compactly supported Dehn twists~\cite[Sections~9.1.1--9.1.2]{FMSW26}. On homology the induced reflections act by
\begin{equation}
 \Gamma\longmapsto-\Gamma,\qquad
 \Gamma\longmapsto2F-\Gamma.
 \label{eq:Weyl}
\end{equation}
On homology the reflections generate translations by $2F$, and on the corresponding real period they act by reflections. The orbit of every normalized real period $r=(4\pi^2)^{-1}\int_\Gamma\omega_I$ meets $[0,1]$. A reflection also changes the sign of the corresponding complex period, which poses no restriction since $b$ ranges over $\C$. The diffeomorphisms are compactly supported, so they do not change the asymptotic coefficient.

\subsection{The \texorpdfstring{$A_2$}{A2} family}
\label{sec:A2}

For the family~\eqref{eq:familyA2}, we write the quartic numerator as $\eta^2$, so that $\lambda=\eta\dd z/z$. The labeled eigenvalues at the two poles are
\begin{equation*}
 \lambda_{0,+}=c\frac{\dd z}{z}+\cdots,
 \qquad
 \lambda_{\infty,+}=-\left(w^{-3}+\frac a2w^{-2}
                              +\mu w^{-1}\right)\dd w+\cdots.
\end{equation*}
At $z=0$ we assign weights $(\alpha_0,1-\alpha_0)$, with $0\le\alpha_0\le1/2$ and residue $-c$ on the parabolic line of weight $1-\alpha_0$; at infinity we allow the real representatives of the weights $(\alpha_\infty,1-\alpha_\infty)$. The extension has $\deg E=-2$ before the weights are normalized, and the generic Hitchin fibers are $\Pic^0$. We use the period coordinates
\begin{equation*}
 x=\alpha_0+\alpha_\infty-1,\qquad y=1-2\alpha_0,
 \qquad Q_1=2i(\mu-c),\qquad Q_2=4ic.
\end{equation*}
We use the following closed alcove and require the indicated periods on root classes to be nonzero:
\begin{equation}
 x,y\ge0,\quad x+y\le1,
 \qquad
 (x,Q_1),(y,Q_2),(x+y,Q_1+Q_2)\notin\Z\times\{0\}.
 \label{eq:A2domain}
\end{equation}

\subsubsection{Two rational curves and their periods}

For $a=\mu=c=u=0$ and weights in the interior of the alcove, we use the family on $E=\cO(-1)^2$ with Higgs fields
\begin{equation*}
 \Phi_{\mathrm{sp},[v_0:v_1]}=
 \begin{pmatrix}z^2&0\\v_0z+v_1&-z^2\end{pmatrix}\frac{\dd z}{z},
 \qquad[v_0:v_1]\in\PP^1,
\end{equation*}
with the fiber of the second summand as the flag at $z=0$. The invariant second summand has parabolic degree $-x$. The other saturated eigenline bundle, generated by $(2z^2,v_0z+v_1)$, has degree $-3$ for $v_1\ne0$ and $-2$ for $v_1=0$. Its fiber at $z=0$ agrees with that of the second summand in both cases, and its parabolic degrees are $x+y-2$ and $x+y-1$, respectively. All these parabolic degrees are negative. Constant lower-triangular automorphisms preserving the flag translate the lower-left coefficient by a multiple of $z^2$ and rescale it. The remaining parameters $[v_0:v_1]$ give the embedded rational curve $\Gamma_1$.

A universal bundle for this family is
\begin{equation*}
 \mathcal E_1=\cO(-1)\oplus
       (\cO(-1)\boxtimes\cO_{\Gamma_1}(1)).
\end{equation*}
For $(n_1,n_2)=(1,1)$, the discriminant in~\eqref{eq:discriminant} vanishes. At $z=0$ the second summand is the parabolic line $F_0$ of weight $1-\alpha_0$, and the quotient $E_0/F_0$ has weight $\alpha_0$. Their degrees in increasing weight order are $0,1$, including where the nilpotent residue vanishes. At infinity the eigenline $L_+$ is a graph over the first summand and the eigenline $L_-$ is the second summand. Both ordered pairs thus have normalized degrees $(-1/2,+1/2)$. Substituting the discriminant integral and these degrees in the real Chern--Weil formula and in \eqref{eq:complexclass} gives
\begin{equation}
 (r_{\Gamma_1},Q_{\Gamma_1})
   =(\alpha_0+\alpha_\infty-1,\,2i(\mu-c))=(x,Q_1).
 \label{eq:A2firstperiod}
\end{equation}
We compute the characteristic numbers on the rational curve just constructed and use the transported classes to evaluate the residue formula at the other parameters in the family.

At $[v_0:v_1]=[1:0]$ the residue vanishes, so keeping the pair fixed and varying the flag gives $\Gamma_2=\PP(E_0)$. The invariant second summand has parabolic degree $-x$ when its fiber at $z=0$ is the chosen flag and $-x-y$ otherwise. The other line subbundle has degree $-2$, with parabolic degrees $x+y-1$ and $x-1$. Every choice of flag gives a stable parabolic Higgs bundle, and the pair has only scalar automorphisms. The underlying universal bundle is pulled back from $\PP^1_z$, so its discriminant vanishes. At $z=0$ the tautological parabolic line bundle is $\cO(-1)$ and the quotient line bundle is $\cO(1)$, giving normalized degrees $(+1,-1)$ in increasing weight order; the degrees at infinity vanish. Thus
\begin{equation}
 (r_{\Gamma_2},Q_{\Gamma_2})=(1-2\alpha_0,4ic)=(y,Q_2).
 \label{eq:A2flagperiod}
\end{equation}
The two graded lines at $z=0$ contribute $4ic$ to the normalized complex period.

\subsubsection{The type $\mathrm{IV}$ fiber and the classes of its components}

The quartic $P(z)=z^4+uz$ at these parameter values has binary quartic invariants $I(P)=0$, $J(P)=-27u^2$, from which the minimal discriminant has vanishing order four. The fibration has a section on $\cO\oplus\cO(-2)$ given by
\begin{equation*}
 \begin{pmatrix}0&z^4+uz\\1&0\end{pmatrix}\frac{\dd z}{z}
\end{equation*}
with the fiber of the second summand as the flag. At $u=0$ the Higgs bundle on the section has two invariant line subbundles of degree $-2$ and parabolic degrees $x+y-1,-x-1$. For $u\ne0$ the spectral curve is irreducible. Every point of the section thus represents a stable Higgs bundle. The smooth symplectic surface is relatively minimal, so its central fiber has type $\mathrm{IV}$. At infinity the vanishing orders in a minimal Weierstrass equation are $(\ge3,4,8)$, giving a fiber of type $\mathrm{IV}^*$ in the rational elliptic compactification.

In a local chart near the intersection of the components, choose the first basis vector in the leading eigenline at infinity for which the upper-left coefficient is monic, and normalize the lower-left entry to $z-z_*$. A constant upper-triangular unipotent gauge transformation then removes the linear term of the upper-left entry, leaving $z^2+a_0$. The upper-right entry $B$, of degree at most one, satisfies
\begin{equation*}
 B=-2a_0z+u-2a_0z_*,\qquad a_0^2+2z_*^2a_0=uz_*.
\end{equation*}
The parabolic line at $z=0$ has the form $\operatorname{span}(r,1)$. Since the residue preserves the flag and acts trivially on both graded pieces, we have
\begin{equation*}
 a_0=-z_*r,\qquad a_0r+u-2a_0z_*=0,
 \qquad u=z_*r(r-2z_*).
\end{equation*}
Thus the central fiber $z_*r(r-2z_*)=0$ has three smooth reduced branches with distinct tangent lines. Denoting the third component by $\Gamma_3$, we obtain
\begin{equation*}
 F=\Gamma_1+\Gamma_2+\Gamma_3,
 \quad\Gamma_i^2=-2,
 \quad\Gamma_i\cdot\Gamma_j=1\quad(i\ne j).
\end{equation*}
The section has intersection number one with $F$, so the fiber class is primitive. The quotient of the intersection lattice by the fiber class for a fiber of type $\mathrm{IV}^*$ at infinity has determinant three by~\eqref{eq:quotient-lattices}. The Gram matrix of $\Gamma_1,\Gamma_2$ has the same determinant, so~\eqref{eq:lattice-index} gives
\begin{equation*}
 H_2=\Z F\oplus\Z\Gamma_1\oplus\Z\Gamma_2,
 \qquad
 (\Gamma_i\cdot\Gamma_j)_{1\le i,j\le2}
 =\begin{pmatrix}-2&1\\1&-2\end{pmatrix}.
\end{equation*}
The classes of square $-2$ in the lattice are precisely $\pm\Gamma_1+nF$, $\pm\Gamma_2+nF$ and $\pm(\Gamma_1+\Gamma_2)+nF$. Requiring the period vector to be nonzero on each such class gives \eqref{eq:A2domain}.

\subsubsection{Excluded parameters and the global marking}

A possible global rational eigenvalue one-form is
\begin{equation*}
 \rho=\left(z+\frac a2+\frac{\varepsilon c}{z}\right)\dd z,
 \qquad\varepsilon=\pm1.
\end{equation*}
Its square equals the quadratic differential in~\eqref{eq:familyA2} precisely when $\mu=\varepsilon c$, $u=a\varepsilon c$. If $k\in\Z$ is the degree of the underlying line bundle, the two possible parabolic degrees are $k+x+1$ and $k+x+y+1$. A reducible polystable Higgs bundle exists exactly when $(x,Q_1)$ or $(x+y,Q_1+Q_2)$ lies in $\Z\times\{0\}$; conversely the direct sum of a line bundle of the required integer degree and the complementary line bundle determined by the determinant gives a reducible polystable parabolic Higgs bundle. The condition on $(y,Q_2)$ ensures that the weights and residue at the tame pole are admissible. If $c=0$, then $0<y<1$ and hence $0<\alpha_0<1/2$. At $y=0$ or $1$ we have $c\ne0$, so the residue contributes a positive term to the off-diagonal indicial operator. The moduli space with these tame and irregular pole conditions is smooth and its metric is complete throughout~\eqref{eq:A2domain}.

Contract the domain by moving the real periods into the alcove interior before contracting the complex parameters. At any positive time the real periods on positive roots lie in $(0,1)$. Near a reducible spectral curve the invariant-line degrees retain the same strict bounds. At a splitting with $c=0$ one has $\mu=0$, so $Q_1=Q_2=0$. The exclusions in~\eqref{eq:A2domain} then place both $x$ and $x+y$ in $(0,1)$, giving the same integer threshold for either flag assignment. These are the local hypotheses of Section~\ref{subsec:topology-transport}.

Across $\alpha_\infty=1$ the leading eigenvalues remain distinct. Lemma~\ref{lem:weight-domains} and the modification in Section~\ref{subsec:calibration} apply. Lemma~\ref{lem:marking-transport} now extends~\eqref{eq:A2firstperiod}--\eqref{eq:A2flagperiod} to the contractible domain.

\subsubsection{Affine Weyl action}

The two simple-root spheres $\Gamma_1,\Gamma_2$ and the third component $\Gamma_3=F-\Gamma_1-\Gamma_2$ of the type $\mathrm{IV}$ fiber realize the three affine simple reflections. Their compactly supported Dehn twists act by
\begin{equation*}
 r_1(x,y)=(-x,x+y),\quad r_2(x,y)=(x+y,-y),
 \quad r_3(x,y)=(1-y,1-x).
\end{equation*}
The three reflections induced by the Dehn twists generate the affine Weyl group of type $A_2$ with the closed fundamental alcove described above. They act simultaneously on the complex periods by
\begin{equation*}
 (-Q_1,Q_1+Q_2),\quad(Q_1+Q_2,-Q_2),\quad(-Q_2,-Q_1),
\end{equation*}
respectively, and are the identity outside a compact set. The inverse formulas in Section~\ref{subsec:inverse-periods} use this marking.
\subsection{The \texorpdfstring{$D_0$}{D0} family}
\label{sec:D0}

In~\eqref{eq:familyD0} the Higgs field has ramified poles of order two at $z=0,\infty$, both with weights $(1/4,3/4)$; the chosen extension has degree $-2$ and fixed parabolic determinant. A dilation $z=\varkappa z'$ sends $(p_0,p_\infty)$ to $(p_0/\varkappa,p_\infty\varkappa)$, so $h=p_0p_\infty$ is invariant. For the realization argument we use $h>0$ and $p_0=p_\infty=\sqrt h$.

Since the spectral quadratic differential has poles of odd order at $z=0,\infty$, no rational eigenvalue one-form exists for any $u$. By Proposition~\ref{prop:ramified-quotient}, the metric is smooth and complete for all $p_0,p_\infty\ne0$. The spectral sheaf has degree zero, and the Hitchin fiber is the compactified Jacobian $\overline{\Pic}^{\,0}$ of the corresponding integral spectral curve of arithmetic genus one. The minimal elliptic compactification adds an $I_4^*$ fiber, whose complement has $H_2=\Z F$ by~\cite[Section~3.1.11]{LL23}. Together with the periods of the hyperk\"ahler forms on the fiber, this gives
\begin{equation}
 H_2(\cM_h;\Z)=\Z F,\qquad
 \int_F(\omega_I,\omega_J,\omega_K)=(4\pi^2,0,0).
 \label{eq:D0-data}
\end{equation}
The fiber periods are independent of the positive parameter $h$. In Section~\ref{sec:periodend} we compute how this parameter enters the constant term of the expansion at infinity.
\subsection{The \texorpdfstring{$D_1$}{D1} family}
\label{sec:D1}

For the spectral family~\eqref{eq:familyD1} we write $\lambda_0=pb$. The unramified pole at $z=0$ has labeled eigenvalue one-forms
\begin{equation}
 \lambda_+=(pz^{-2}+mz^{-1}+O(1))\dd z\,,\qquad
 \lambda_-=-\lambda_+\,.
 \label{eq:D1-labels}
\end{equation}
At $z=0$ the two eigenline bundles have weights~$(\alpha,1-\alpha)$, with~$\alpha\in\R$. At the ramified pole of order two at infinity the weights are necessarily~$(1/4,3/4)$. We work with degree~$-2$ extensions and the prescribed parabolic determinant. The real representatives of the weights specify the local holomorphic lattices of the eigenline bundles. The elementary modifications below give extensions with normalized weights. For realization we take~$p=1$ and~$b=\lambda_0>0$.

\begin{proposition}
\label{prop:D1-dictionary}
For every~$\alpha\in\R$, $m\in\C$ and~$\lambda_0>0$, the moduli space equipped with its $L^2$ metric is a smooth complete ALG$^*$ hyperk\"ahler four-manifold with~$\nu=3$. With respect to the integral homology basis
\begin{equation}
 H_2(M;\Z)=\Z\Gamma\oplus\Z F\,,\qquad
 \Gamma^2=-4\,,\qquad \Gamma\cdot F=F^2=0\,,
 \label{eq:D1-lattice}
\end{equation}
where $\Gamma\cdot O=0$ for the zero section $O$ of the elliptic compactification and $\int_\Gamma c_1(L_+^0)=1$, the periods are
\begin{align}
 \int_F\omega_I&=4\pi^2\,,&\int_F\Omega_I&=0\,,\notag\\
 \int_\Gamma\omega_I&=4\pi^2(1-2\alpha)\,,&
 \int_\Gamma\Omega_I&=16\pi^2i m\,.
 \label{eq:D1-periods}
\end{align}
\end{proposition}

\subsubsection{Stability and the relative Jacobian}

The pole of $\lambda^2$ at infinity has odd order, so this quadratic differential cannot be the square of a rational one-form. Every filtered Higgs bundle is therefore stable for all~$\alpha,m,u$. By Proposition~\ref{prop:ramified-quotient}, the moduli space is smooth, its metric is complete, and the solutions depend smoothly on the finite-dimensional parameters. The local model at the ramified pole stays fixed, and the weighted Sobolev domains at the unramified pole of order two remain equivalent as the weights vary. In particular, the domain calculation in Lemma~\ref{lem:weight-domains} covers equal and integer weights.

For~$\chi=z^2\lambda/\dd z$, the spectral equation is
\begin{equation}
 \chi^2=p^2+2pmz+uz^2+bz^3\,,\qquad
 \lambda=\chi\frac{\dd z}{z^2}\,,\qquad
 \partial_u\lambda=\frac{\dd z}{2\chi}\,.
 \label{eq:D1-spectral-cubic}
\end{equation}
For every~$u$ the spectral curve is integral of arithmetic genus one. The spectral sheaves have degree zero, and the norm map takes values in~$\Pic^0(\PP^1)$, a point. Thus the Hitchin fiber is the compactified Jacobian $\overline{\Pic}^{\,0}$ of the spectral curve. The section at infinity induces an isomorphism from the relative compactified Jacobian to the Weierstrass surface
\begin{equation}
 Y^2=X^3+uX^2+2\lambda_0 mX+\lambda_0^2\,,\qquad X=bz\,,\quad Y=b\chi\,.
 \label{eq:D1-Weierstrass}
\end{equation}
In this Weierstrass family, $O$ is the section at infinity. The coordinates
\begin{equation}
 P=(0,\lambda_0)\,,\qquad -P=(0,-\lambda_0)
 \label{eq:D1-sections}
\end{equation}
describe the two smooth sections through the points of the spectral curve above the unramified pole. The total space is nonsingular over every $u\in\C$. Vanishing of the $u$-derivative forces $X=0$ and hence $Y=\pm\lambda_0\ne0$, where the $Y$-derivative is nonzero.

\subsubsection{The primitive class of square \texorpdfstring{$-4$}{-4}}

The minimal rational elliptic compactification has a fiber of type~$I_3^*$ at infinity. The nonidentity components of that fiber have intersection graph~$D_7$. We determine which component meets $P$ and compute the orthogonal projection of the section class.

In the coordinates $t=u^{-1}$, $x=t^2X$ and $y=t^3Y$, the minimal equation is
\begin{equation}
 y^2=x^3+t x^2+2\lambda_0 m t^4x+\lambda_0^2t^6\,.
 \label{eq:D1-minimal-infinity}
\end{equation}
The chord-and-tangent formula gives
\begin{equation}
 2P=(m^2-u,\,mu-m^3-\lambda_0)\,.
 \label{eq:D1-doubleP}
\end{equation}
Both~$P$ and~$2P$ specialize to the singular point of the minimal Weierstrass fiber. For~$m\ne0$, a second doubling gives
\begin{equation}
 X(4P)=\frac{u^2}{4m^2}+O(u)\,,\qquad
 Y(4P)=-\frac{u^3}{8m^3}+O(u^2)\,.
 \label{eq:D1-fourP}
\end{equation}
Thus~$4P$ specializes to a nonsingular point of~$y^2=x^3$. The component group is~$\Z/4$, and~$P$ meets one of the two components at the ends of the short arms of the Dynkin diagram of the nonidentity components. Each represents an element of order four in the component group. The component met by the section is locally constant in~$(m,\lambda_0)$ as long as~$\lambda_0\ne0$, so the same description holds at~$m=0$.

We number the~$D_7$ chain~$1$--$2$--$3$--$4$--$5$, with branches~$5$--$6$ and~$5$--$7$. The section~$P$ meets the component with label~$7$. Let $C$ be the Cartan matrix of $D_7$ in this numbering. Direct inversion gives
\begin{equation}
 (c_1,\ldots,c_7)^T=C^{-1}e_7=(1/2,1,3/2,2,5/2,5/4,7/4)^T\,.
 \label{eq:D1-Cartan}
\end{equation}
Here~$c_i$ denotes the entry corresponding to the component~$\Theta_i$. Since~$P\cdot O=0$, $O^2=-1$ and~$P\cdot F=O\cdot F=1$, the class
\begin{equation}
 v=P-O-F+\sum_{i=1}^7c_i\Theta_i
 \label{eq:D1-Shioda}
\end{equation}
is orthogonal to~$O,F,\Theta_i$. Its square is
\begin{equation}
 v^2=-2+c_7=-\frac14\,.
\end{equation}
The projection is the one described in \cite[Lemmas~8.1--8.2 and equation~(8.2)]{Shioda90}. The same self-intersection number follows from Shioda's height formula~\cite[Theorem~8.6]{Shioda90}. In particular
\begin{align}
 \Gamma_0=4v={}&4P-4O-4F+2\Theta_1+4\Theta_2+6\Theta_3\notag\\
 &+8\Theta_4+10\Theta_5+5\Theta_6+7\Theta_7
 \label{eq:D1-Gamma0}
\end{align}
is integral and orthogonal to every boundary component, with
\begin{equation}
 \Gamma_0^2=-4\,,\qquad \Gamma_0\cdot O=0\,,\qquad
 \Gamma_0\cdot P=-1\,.
 \label{eq:D1-intersections}
\end{equation}
The projection defines a homomorphism from the group of sections, so~$v(-P)=-v(P)$. Hence~$\Gamma_0\cdot(-P)=1$, and elliptic inversion sends~$\Gamma_0$ to~$-\Gamma_0$.

The pairing of $F,\Gamma_0$ with the integral divisors~$O,P-O$ gives the unimodular diagonal matrix~$\diag(1,-1)$, so the span of $F,\Gamma_0$ is a primitive sublattice of $H_2(Y;\Z)$, where $Y$ is the elliptic compactification. Since the boundary is a tree of rational curves, the compactly supported cohomology sequence for the complement identifies~$H_2(M;\Z)$ with the sublattice of $H_2(Y;\Z)$ orthogonal to every boundary component. By \cite[Section~3.1.10]{LL23} this sublattice has rank two, so~$F,\Gamma_0$ form an integral basis.

\subsubsection{Eigenline bundle degrees and periods}

The integral spectral curves form a family $\mathcal S\to\C_u$ with relative compactified Jacobian $\cM\to\C_u$. We use the relative universal rank-one torsion-free sheaf on $\mathcal S\times_{\C_u}\cM$, rigidified along $O$. Since the sections $P$ and $-P$ meet the smooth locus of every spectral curve, restriction to them gives line bundles over the whole moduli space, including the fibers over singular curves. These are the eigenline bundles of the universal projective Higgs bundle at $z=0$. The normalized rational first Chern class
\begin{equation}
 c_1(L_+^0)=\frac12\bigl(c_1(L_+)-c_1(L_-)\bigr)
 \label{eq:D1-projective-eigenline}
\end{equation}
is unchanged when the universal family is tensored by a line bundle pulled back from the parameter space. Under the map $Q\mapsto\cO(Q-O)$ on the smooth generic fiber, the restrictions to the two sections have divisor classes~$P-O$ and~$-P-O$, respectively, up to classes pulled back from the base or supported on the compactification boundary. The finite fibers are irreducible, so any vertical correction there is a multiple of a full fiber. Each of these corrections pairs trivially with $\Gamma_0$. Therefore
\begin{equation}
 \int_{\Gamma_0}c_1(L_+^0)
 =\frac12\Gamma_0\cdot(P-(-P))=-1\,.
 \label{eq:D1-eigenline-degree}
\end{equation}
With the orientation~$\Gamma=-\Gamma_0$, we have~$\int_\Gamma c_1(L_+^0)=1$. Then~$\int_\Gamma c_1(L_-^0)=-1$. Elliptic inversion fixes $O$, so both $\Gamma$ and its image under inversion are orthogonal to $O$. Since $F\cdot O=1$, the image is exactly $-\Gamma$.

For Corollary~\ref{cor:ramified-periods} take the degree-two map
$\PP^1_\zeta\to\PP^1_z$, $z=\zeta^2$, branched at $0$ and $\infty$.
At $\zeta=0$, pulling back~\eqref{eq:D1-labels} gives residues $2m,-2m$ and real weight representatives $2\alpha,2(1-\alpha)$. Restriction of the pulled-back family to $\{\zeta=0\}\times\Gamma$ gives the original eigenlines over $\{z=0\}\times\Gamma$. Their normalized degrees on $\Gamma$ are therefore $1,-1$. At $\zeta=\infty$, the two adapted eigenlines both have weight $1/2$ by the ramified model of order two following~\eqref{eq:ramified-local-model}. Their residues vanish. The complex period on the base is half the expression on the cover,
\begin{equation}
 \int_\Gamma\Omega_I
 =\frac{8\pi^2 i}{2}\bigl((2m)(1)+(-2m)(-1)\bigr)
 =16\pi^2im\,.
 \label{eq:D1-complex-period}
\end{equation}
For the real period, use the same adapted lattices with the indicated real representatives as $\alpha$ varies. On the transported cycle the discriminant integral of this fixed topological extension is locally constant. The contribution at $\zeta=\infty$ is fixed as well. The weight term at $\zeta=0$ in Proposition~\ref{prop:CW} therefore gives
\begin{equation}
 \frac{\dd}{\dd\alpha}\left(\frac1{4\pi^2}\int_\Gamma\omega_I\right)
 =-\frac12\frac{\dd}{\dd\alpha}
     \bigl((2\alpha)(1)+2(1-\alpha)(-1)\bigr)=-2\,.
 \label{eq:D1-real-slope}
\end{equation}
If a weight is normalized by an elementary modification, Section~\ref{subsec:calibration} gives the compensating change in the discriminant integral. Both calculations thus extend across integer weights.
The quantity $r_\Gamma+2\alpha$ is independent of~$\alpha,m,\lambda_0$. At~$\alpha=1/2,m=0$ the map~$\Phi_{\mathrm{sp}}\mapsto-\Phi_{\mathrm{sp}}$ preserves the unordered formal types and the weights and induces an $L^2$ isometry preserving~$\omega_I$. On the degree-zero compactified Jacobian it acts by elliptic inversion, which fixes~$O$ and sends~$\Gamma$ to~$-\Gamma$. The real period on the class sent to its negative therefore vanishes at these parameter values. Equation~\eqref{eq:D1-real-slope} proves
\begin{equation}
 \int_\Gamma\omega_I=4\pi^2(1-2\alpha)\,.
 \label{eq:D1-real-period}
\end{equation}
The fiber periods are~\eqref{eq:fiberperiod}.

\subsubsection{Integer weights and the range of the real period}

The operation of tensoring a spectral sheaf by
\begin{equation}
 \cO\bigl(kP-k(-P)\bigr)\,,\qquad k\in\Z\,,
 \label{eq:D1-modification}
\end{equation}
is translation by~$2kP$ on the compactified Jacobian. At~$P$ the modification multiplies a local frame of the spectral sheaf by~$z^{-k}$. Since $h\sim|z|^{2\alpha}$, the weight changes from~$\alpha$ to~$\alpha-k$, with the opposite shift on the other eigenline bundle. The determinant and eigenvalue one-forms are unchanged. The modification changes the extension of the same filtered harmonic bundle and induces an isometry of the $L^2$ metrics. It extends over every finite singular fiber because both marked points $P,-P$ lie in the smooth locus of that spectral curve. In a unitary frame the corresponding diagonal loop has determinant one and extends across an annulus to the identity.

For each section $Q$, the map $T_Q$ is fiberwise translation by $Q$. Orthogonal projection to $\langle O,F,\Theta_1,\ldots,\Theta_7\rangle^\perp$ sends $T_{Q*}\Gamma_0$ to $4(v(P+Q)-v(Q))=4v(P)$. For translation by twice the marked section, the difference $T_{2P*}\Gamma_0-\Gamma_0$ is therefore a rational multiple of $F$. Moreover
\begin{equation}
 (T_{2P*}\Gamma_0)\cdot O
 =\Gamma_0\cdot(-2P)=2\,.
\end{equation}
Together with~$F\cdot O=1$ this gives
\begin{equation}
 T_{2P*}\Gamma_0=\Gamma_0+2F\,,\qquad
 T_{2P*}\Gamma=\Gamma-2F\,.
 \label{eq:D1-translation}
\end{equation}
The isometry between the moduli spaces with weights $\alpha$ and $\alpha-1$ gives $\int_{\Gamma-2F}\omega_{I,\alpha-1}=\int_\Gamma\omega_{I,\alpha}$. Since $r_F=1$, this is precisely $r_\Gamma(\alpha-1)=r_\Gamma(\alpha)+2$. An interval of weights of length one gives an interval of real periods of length two, and the translates of this interval under the changes of marking cover $\R$. Since $(n\Gamma+kF)^2=-4n^2$, this lattice has no class of square $-2$. The normalized period vector on $\Gamma$ is consequently unrestricted in $\R\times\C$.

\subsection{The \texorpdfstring{$D_2$}{D2} family}
\label{sec:D2}

\subsubsection{Parameters and periods}

We write $h=p_0p_\infty$ in the spectral family~\eqref{eq:familyD2}. The variable~$u\in\C$ is the coordinate on the Hitchin base when the coefficients~$m_0,m_\infty\in\C$ of the simple-pole terms and the leading coefficients are fixed. We label the formal eigenvalue one-forms by the expansions
\begin{align}
 \lambda_{0,+}&=(p_0z^{-2}+m_0z^{-1}+O(1))\dd z\,,\notag\\
 \lambda_{\infty,+}&=-(p_\infty w^{-2}+m_\infty w^{-1}+O(1))\dd w\,,
 \qquad w=z^{-1}\,.
 \label{eq:D2-labels}
\end{align}
The weights on the eigenline bundles $L_+,L_-$ are~$(\alpha_0,1-\alpha_0)$ at $z=0$ and~$(\alpha_\infty,1-\alpha_\infty)$ at infinity. The chosen degree~$-2$ extension has fixed parabolic determinant and parabolic degree zero. The real representatives of the weights specify the local holomorphic lattices of the eigenline bundles. Elementary modifications shift the weights by integers to give the usual normalized representatives.

The normalized real and complex periods are
\begin{align}
 x&=\alpha_0+\alpha_\infty-1\,,&
 Q_1&=2i(m_\infty-m_0)\,,\notag\\
 y&=\alpha_0-\alpha_\infty\,,&
 Q_2&=-2i(m_0+m_\infty)\,.
 \label{eq:D2-period-coordinates}
\end{align}
Their admissible domain is
\begin{equation}
 \mathcal P_{D_2}
 =\left\{(\alpha_0,\alpha_\infty,m_0,m_\infty)\in\R^2\times\C^2:
 (x,Q_1),(y,Q_2)\notin\Z\times\{0\}\right\}\,.
 \label{eq:D2-domain}
\end{equation}
For the global marking and the realization argument, we take~$h>0$ and set~$p_0=p_\infty=\sqrt h$.

\begin{proposition}
\label{prop:D2-dictionary}
For every parameter in~$\mathcal P_{D_2}\times\R_{>0}$, the moduli space with its $L^2$ metric is a smooth complete ALG$^*$ hyperk\"ahler four-manifold with~$\nu=2$. There is a global integral marking
\begin{equation}
 \begin{gathered}
 H_2(M;\Z)=\Z F\oplus\Z\Gamma_1\oplus\Z\Gamma_2\,,\qquad
 F^2=F\cdot\Gamma_i=0\,,\\
 (\Gamma_i\cdot\Gamma_j)_{i,j=1}^2=
 \begin{pmatrix}-2&0\\0&-2\end{pmatrix}\,,
 \end{gathered}
 \label{eq:D2-lattice}
\end{equation}
in which the periods are
\begin{equation}
 \begin{gathered}
 \frac1{4\pi^2}\int_F(\omega_I,\Omega_I)=(1,0)\,,\\
 \frac1{4\pi^2}\int_{\Gamma_1}(\omega_I,\Omega_I)=(x,Q_1)\,,\qquad
 \frac1{4\pi^2}\int_{\Gamma_2}(\omega_I,\Omega_I)=(y,Q_2)\,.
 \end{gathered}
 \label{eq:D2-periods}
\end{equation}
\end{proposition}

\subsubsection{Stability and the integral lattice}

An invariant line in the generic fiber of the bundle determines a rational eigenvalue one-form. After exchanging the two eigenvalues if necessary, we may take its leading coefficient at $z=0$ to be $p_0$. It then has one of the forms
\begin{equation}
 \rho_\varepsilon=\left(\frac{p_0}{z^2}+\frac{m_0}{z}
 +\varepsilon p_\infty\right)\dd z\,,\qquad \varepsilon\in\{1,-1\}\,.
 \label{eq:D2-rational-eigenform}
\end{equation}
Comparison of its square with~\eqref{eq:familyD2} gives
\begin{equation}
 m_\infty=\varepsilon m_0\,,\qquad u=m_0^2+2\varepsilon h\,.
 \label{eq:D2-square-loci}
\end{equation}
For~$\varepsilon=1$, the parabolic degree of a line bundle of degree~$k\in\Z$ is~$k+x+1$; for~$\varepsilon=-1$, it is~$k+y+1$. A reducible polystable Higgs bundle can occur only when~$(x,Q_1)\in\Z\times\{0\}$ or~$(y,Q_2)\in\Z\times\{0\}$. Conversely, on either locus we choose a line bundle of the required degree. On this line bundle and its complement determined by the determinant, the Higgs fields~$\rho_\varepsilon$ and~$-\rho_\varepsilon$ give a polystable direct sum. Thus~\eqref{eq:D2-domain} is precisely the domain on which every polystable Higgs bundle is irreducible. Since both poles have regular semisimple leading terms and there are no reducible polystable points on this domain, the Biquard--Boalch construction~\cite[Theorem~5.4]{BB04} gives smooth moduli spaces with complete $L^2$ metrics. Proposition~\ref{prop:unramified-slices} and Lemma~\ref{lem:weight-domains} give local slices and smooth dependence on the real representatives of the weights, including equal and integer weights.

For generic~$u$, we write~$\eta=z^2\lambda/\dd z$. The compact spectral curve has affine equation
\begin{equation}
 \eta^2=p_\infty^2z^4+2p_\infty m_\infty z^3+uz^2
 +2p_0m_0z+p_0^2\,.
 \label{eq:D2-quartic}
\end{equation}
It is a smooth genus-one double cover. Degree-zero spectral line bundles push forward to bundles of degree~$-2$. Since the norm map takes values in~$\Pic^0(\PP^1)$, a point, the Hitchin fiber is $\Pic^0$ of the spectral curve. In particular, the moduli spaces are nonempty and the angles have the $2\pi$ periods of the Picard fiber.

We construct the compact cycles at~$m_0=m_\infty=0$ and weights satisfying
\begin{equation}
 x>0\,,\qquad y>0\,,\qquad x+y<1\,.
 \label{eq:D2-reference-chamber}
\end{equation}
At~$u=2h$ and~$u=-2h$, respectively, the quartic is~$P_1^2$ and~$P_2^2$, where
\begin{equation}
 P_1=p_\infty z^2+p_0\,,\qquad
 P_2=p_0-p_\infty z^2\,.
\end{equation}
We choose representatives~$\ell_{i,0},\ell_{i,1}$ of a basis of~$H^0(\cO(2))/\C P_i$. On~$E=\cO(-1)\oplus\cO(-1)$ the family of Higgs fields is given by
\begin{equation}
 \Phi_{\mathrm{sp},i,[v_0:v_1]}=
 \begin{pmatrix}P_i&0\\v_0\ell_{i,0}+v_1\ell_{i,1}&-P_i\end{pmatrix}
 \frac{\dd z}{z^2}\,,\qquad [v_0:v_1]\in\PP^1\,.
 \label{eq:D2-spheres}
\end{equation}
A constant lower-triangular unipotent bundle automorphism adds a multiple of $P_i$ to the lower-left entry, and a diagonal automorphism rescales that entry. Thus the lower-left entry gives a point of the projective line obtained from the quotient space of quadratics above.

The eigenline bundle given by the second summand has degree~$-1$ and parabolic degree~$-x$ for~$i=1$ and~$-y$ for~$i=2$. The other saturated eigenline bundle is generated meromorphically by~$(2P_i,v_0\ell_{i,0}+v_1\ell_{i,1})$. It has degree~$-3$ if the two entries have no common zero and~$-2$ if they share one zero. Since the lower-left entry is nonzero modulo~$P_i$, the entries cannot share both zeros. The corresponding parabolic degrees are~$x-2,x-1$ for~$i=1$ and~$y-2,y-1$ for~$i=2$. All these degrees are negative in~\eqref{eq:D2-reference-chamber}, so both $\PP^1$-families consist of stable Higgs bundles, including the members for which the two entries share a zero.

The invariant line subbundle of degree~$-1$ is unique, so every isomorphism between two members of~\eqref{eq:D2-spheres} preserves it, and the remaining constant lower-triangular automorphisms act by translations by multiples of $P_i$ and scalar rescalings, as used to form the quotient. For a constant infinitesimal gauge $g$, the upper-right entry of $[g,\Phi_{\mathrm{sp},i}]$ is $-2P_i g_{12}\,\dd z/z^2$. A variation within the displayed family has zero upper-right entry, so $g_{12}=0$. The remaining gauge directions in the lower-left entry are precisely addition of a multiple of $P_i$ and rescaling. The induced map of the projective quotient is therefore an injective immersion; compactness gives embedded curves~$\Gamma_1,\Gamma_2\cong\PP^1$. By adjunction, $\Gamma_i^2=-2$, and~$\Gamma_1\cdot\Gamma_2=0$ because the curves lie in different fibers.

For these parameter values a section of the Hitchin fibration on~$\cO\oplus\cO(-2)$ is given by
\begin{equation}
 \Phi_{\mathrm{sp},u}=
 \begin{pmatrix}0&p_\infty^2z^4+uz^2+p_0^2\\1&0\end{pmatrix}
 \frac{\dd z}{z^2}\,.
 \label{eq:D2-section}
\end{equation}
Away from~$u=\pm2h$ the spectral curve is irreducible. At either exceptional value the Higgs bundle on the section has two invariant line subbundles of degree~$-2$, with parabolic degrees~$x-1,-x-1$ at~$2h$ and~$y-1,-y-1$ at~$-2h$. Thus every point of the section represents a stable Higgs bundle. For $P(z)=p_\infty^2z^4+uz^2+p_0^2$, these give
\begin{equation}
 I(P)=u^2+12h^2\,,\qquad
 \operatorname{disc}_z(P)=16h^2(u^2-4h^2)^2\,.
 \label{eq:D2-reference-invariants}
\end{equation}
Both singular fibers at finite values have type~$I_2$. At infinity the minimal Weierstrass equation of the rational elliptic compactification has vanishing orders~$(2,3,8)$, giving a fiber of type~$I_2^*$. Since the moduli space is a smooth holomorphic symplectic surface, it is relatively minimal, so the two fibers just described are its singular fibers over $\C_u$.

The section~\eqref{eq:D2-section} meets $F$ once, so $F$ is primitive. For the complement of a fiber of type $I_2^*$, the quotient lattice in~\eqref{eq:quotient-lattices} has determinant four. The embedded spheres~\eqref{eq:D2-spheres} have Gram matrix $\diag(-2,-2)$ with the same determinant, so the index formula~\eqref{eq:lattice-index} gives~\eqref{eq:D2-lattice}. The classes of square~$-2$ in integral homology are exactly
\begin{equation}
 \pm\Gamma_1+nF\,,\qquad \pm\Gamma_2+nF\,,\qquad n\in\Z\,.
 \label{eq:D2-roots}
\end{equation}

\subsubsection{Compact periods}

Both poles have regular semisimple leading terms, as required in Proposition~\ref{prop:CW} and~\eqref{eq:complexclass}.

A universal bundle for either sphere is
\begin{equation}
 \mathcal E_i=\cO(-1)\oplus
 \bigl(\cO(-1)\boxtimes\cO_{\Gamma_i}(1)\bigr)\,.
\end{equation}
Here $(n_1,n_2)=(1,1)$ in~\eqref{eq:discriminant}, so $\Delta(\mathcal E_i)=0$. At either puncture, the eigenline transverse to the second summand is a graph over the first. This eigenline and the second summand have normalized degrees $(-1/2,1/2)$ by the calculation following~\eqref{eq:discriminant}. In the eigenvalue order~\eqref{eq:D2-labels}, these degrees are
\begin{equation}
 \begin{array}{c|cc}
  &z=0\quad(+,-)&z=\infty\quad(+,-)\\\hline
  \Gamma_1&(-1/2,1/2)&(-1/2,1/2)\\
  \Gamma_2&(-1/2,1/2)&(1/2,-1/2)
 \end{array}\,.
 \label{eq:D2-degree-table}
\end{equation}

The residues at $0$ and $\infty$ are $(m_0,-m_0)$ and $(-m_\infty,m_\infty)$. Pairing these and the weights with~\eqref{eq:D2-degree-table} in~\eqref{eq:CW} and~\eqref{eq:complexclass} gives~\eqref{eq:D2-periods} with the four values in~\eqref{eq:D2-period-coordinates}; the fiber periods are~\eqref{eq:fiberperiod}.

\subsubsection{Global marking and integer weights}

The affine change of coordinates~\eqref{eq:D2-period-coordinates} is a diffeomorphism from the parameter domain onto
\begin{equation}
 \mathcal P_{D_2}\cong
 \left((\R\times\C)\setminus(\Z\times\{0\})\right)^2\,.
\end{equation}
A loop in either factor bounds a disk in the product of the real line and the complex plane. The disk meets at most finitely many deleted points, and general position moves its interior off them. Thus each factor is simply connected. The additional parameter $h$ ranges over $\R_{>0}$.

If $Q_i\ne0$, the corresponding splitting is absent locally. At $Q_i=0$ its real coordinate stays between consecutive integers. Thus near each reducible curve in the product of the pole-parameter space and Hitchin base, fixed strict integer thresholds determine a general rational polarization. The two branches over each pole stay distinct, so Section~\ref{subsec:topology-transport} gives local properness.

Lemma~\ref{lem:weight-domains} supplies the common weighted domains at both poles for arbitrary real representatives. Lemma~\ref{lem:marking-transport} and simple connectivity extend the marking over $\mathcal P_{D_2}\times\R_{>0}$. The discriminant number and normalized eigenline degrees are locally constant in the transported classes, so the same pairings give~\eqref{eq:D2-periods} throughout this domain.

For any integers~$n_0,n_\infty$, the weights at the two poles admit the simultaneous shifts
\begin{equation}
 \alpha_0\longmapsto\alpha_0+n_0\,,\qquad
 \alpha_\infty\longmapsto\alpha_\infty+n_\infty\,.
 \label{eq:D2-weight-shifts}
\end{equation}
In a unitary frame of formal eigenvectors, the shift at a puncture is represented by~$\diag(e^{-in\theta},e^{in\theta})$, with the sign of the connection form chosen to increase the weight of $L_+$ by~$n$. The loop is nullhomotopic in~$\SU(2)$, so it extends to a unitary transformation equal to the identity outside a small disk about the puncture. This transformation identifies the two weighted configuration spaces and their gauge groups while preserving the $L^2$ pairing and all three complex structures. On the holomorphic extensions it acts by opposite elementary modifications. On spectral sheaves, increasing the weight of $L_+$ by $n$ corresponds to tensoring by $\cO\bigl(nP_- - nP_+\bigr)$, where $P_+$ and $P_-$ are the two labeled points over the puncture. Both points lie in the smooth locus of every compact spectral curve.

The isometry from the shifted weights back to the original weights acts on the marked homology classes by
\begin{equation}
 \Gamma_1\longmapsto\Gamma_1+(n_0+n_\infty)F\,,\qquad
 \Gamma_2\longmapsto\Gamma_2+(n_0-n_\infty)F\,.
 \label{eq:D2-marking-shifts}
\end{equation}
To determine the action, we vary~$Q_1,Q_2$ independently. Since the complex periods are preserved, each~$\Gamma_i$ maps to itself plus an integral multiple of~$F$. The real periods in~\eqref{eq:D2-periods} determine these multiples.

The translation lattice in the~$(x,y)$-plane is
\begin{equation}
 \Lambda_{D_2}=\{(n_0+n_\infty,n_0-n_\infty):n_0,n_\infty\in\Z\}\,.
\end{equation}
The image of the closed unit square in the weight plane is~$|x|+|y|\le1$. Its translates by~$\Lambda_{D_2}$ cover~$\R^2$, with disjoint interiors, because they are the images of the integer translates of the unit square under the same invertible affine map. Once one takes the isometries acting on the marking by~\eqref{eq:D2-marking-shifts} into account, normalized weights give the entire marked period domain.

At integer weights Section~\ref{subsec:calibration} identifies the $(0,1)$ extension with its $(0,0)$ normalization, preserving the trace-free metric.

\subsection{The \texorpdfstring{$D_3$}{D3} family}
\label{sec:D3}

In~\eqref{eq:familyD3}, the residues on the graded lines at $z=0,1$ are $(c_i,-c_i)$ with weights $(\alpha_i,1-\alpha_i)$, and the residue acts by $-c_i$ on the parabolic line of weight $1-\alpha_i$. At infinity the eigenvalue one-form $\lambda_{\infty,+}$ is $-(pw^{-2}+mw^{-1}+O(1))\dd w$, with weights $(\alpha_\infty,1-\alpha_\infty)$ on the eigenline bundles $L_+,L_-$ in that order; the extension has $\deg E=-3$ before the weights are normalized. A generic Hitchin fiber is $\Pic^{-1}$ of the spectral curve, with angular periods $2\pi$, and we use the period coordinates
\begin{equation*}
 x=\alpha_0+\alpha_1+\alpha_\infty-1,
 \quad y_i=1-2\alpha_i,
 \quad Q=2i(m-c_0-c_1),\quad Q_i=4ic_i.
\end{equation*}
We take the real periods in the closed alcove
\begin{equation}
 x,y_0,y_1\ge0,\qquad x+y_0+y_1\le1.
 \label{eq:D3alcove}
\end{equation}
The six pairs of periods
\begin{gather}
 (y_0,Q_0),\ (x,Q),\ (y_1,Q_1),\quad
 (x+y_0,Q+Q_0),\ (x+y_1,Q+Q_1),\notag\\
 (x+y_0+y_1,Q+Q_0+Q_1)
 \label{eq:D3-exclusions}
\end{gather}
must all avoid $\Z\times\{0\}$.

\subsubsection{Four components of the central fiber}

At $m=c_0=c_1=u=0$ we work with weights in the interior of the alcove and use the notation $P=pz(z-1)$. On $E=\cO(-1)\oplus\cO(-2)$ the two families
\begin{equation*}
 \Phi_{\mathrm{sp},\pm,[v_0:v_1]}=
 \begin{pmatrix}\pm P&0\\v_0z+v_1&\mp P\end{pmatrix}\frac{\dd z}{z(z-1)}
\end{equation*}
have flags given by the fibers of the second summand at both tame points. For $\Phi_{\mathrm{sp},+}$ this summand has parabolic degree $-x$. The other saturated line subbundle has degree $-3$ generically and $-2$ when $v_0z+v_1$ vanishes at $z=0$ or $1$; its parabolic degrees are $x+y_0+y_1-2$, $x+y_1-1$, $x+y_0-1$, respectively. For $\Phi_{\mathrm{sp},-}$ the second summand has parabolic degree $x+y_0+y_1-1$, and the other saturated line subbundle has parabolic degrees $-x-1,-x-y_0,-x-y_1$. All the parabolic degrees in the two families are negative. Constant diagonal automorphisms act by scalar multiplication on the parameters, so the two projective parameter lines give compact spheres $\Gamma_c,\Gamma_4$.

At $v_0z+v_1=z$ the residue at $z=0$ vanishes, and varying the flag with the Higgs pair fixed defines $\Gamma_0=\PP(E_0)$. The second summand has parabolic degree $-x$ or $-x-y_0$, while the other invariant line subbundle has degree $-2$ and parabolic degree $x+y_0+y_1-1$ when its fiber is the chosen flag and $x+y_1-1$ otherwise, so the parabolic Higgs bundle is stable for every flag. The construction with $v_0z+v_1=z-1$ gives $\Gamma_1$, and for either special polynomial $\ell=z,z-1$, the matrix
\begin{equation*}
 g=\begin{pmatrix}1&-2P/\ell\\0&1\end{pmatrix}
\end{equation*}
is holomorphic and sends the pair represented by the lower-triangular matrix with $P$ in the upper-left entry to the pair with $-P$ in that entry, taking the flag at the zero of $\ell$ to the fiber of the other saturated invariant line subbundle. Each curve obtained by varying a flag therefore meets both curves given by the triangular Higgs matrices.

These four curves are all the components of the central fiber. For a splitting $\cO(k)\oplus\cO(-3-k)$ with $k\ge0$, the lower-left entry of the Higgs field vanishes because its line bundle has degree $-1-2k$, and the Higgs bundle is unstable. In the remaining splitting the matrix has the form
$\left(\begin{smallmatrix}A&B\\\ell&-A\end{smallmatrix}\right)$,
where $\deg A\le2$, $\deg B\le3$, $\deg\ell\le1$. The case $\ell=0$ is again unstable: the minimum parabolic degrees of the first line subbundle in the chosen splitting are $x$ and $1-x-y_0-y_1$. If $\ell$ is linear with root $z_*$, the equation $A^2+B\ell=P^2$ gives $A(z_*)=\pm P(z_*)$; an upper-triangular bundle automorphism reduces $A$ to $\pm P$ and $B$ to zero. For constant $\ell$, the leading coefficient of $A$ is $\pm p$, so the same upper-triangular change of frame applies. The remaining choices of flag are precisely those described above.

To compute the intersection multiplicities, we choose a frame with $\ell=z-z_*$ and $A=P+a_0$. The determinant equation gives
\begin{equation*}
 B=(u-2pa_0)(z+z_*-1),\qquad
 a_0^2+2pz_*(z_*-1)a_0=uz_*(z_*-1).
\end{equation*}
At the tame pole $z=0$ the parabolic line has the form $\operatorname{span}(r,1)$. The conditions on this line give $a_0=-z_*r$ and express the base coordinate in the following smooth local chart near $z_*=0$
\begin{equation*}
 u=\frac{z_*r[r-2p(z_*-1)]}{z_*-1}.
\end{equation*}
The branch $z_*=0$ obtained by varying the flag meets the two branches represented by triangular Higgs matrices at $r=0,-2p$. At each intersection the remaining factor is a unit. At the tame pole $z=1$, we write the flag as $\operatorname{span}(r,1)$. Then $a_0=(1-z_*)r$, and the smooth local chart near $z_*=1$ has equation
\begin{equation*}
 u=\frac{(z_*-1)r(r-2pz_*)}{z_*}\,.
\end{equation*}
The flag component $z_*=1$ meets the other two branches at $(z_*,r)=(1,0),(1,2p)$. Since $p\ne0$, the remaining factor is a unit at each intersection. Thus all four intersections are transverse with multiplicity one. By adjunction each sphere has self-intersection $-2$. With the components of the type $I_4$ fiber listed in cyclic order, its class is given by
\begin{equation*}
 F=\Gamma_0+\Gamma_c+\Gamma_1+\Gamma_4.
\end{equation*}

\subsubsection{Compact periods}

On a curve parametrizing triangular Higgs matrices, we take $\mathcal E=\cO(-1)\oplus(\cO(-2)\boxtimes\cO(1))$. Substituting $(n_1,n_2)=(1,2)$ in~\eqref{eq:discriminant} gives $\int\Delta=1/2$. At each tame pole, the quotient by the flag and the flag itself have normalized degrees $(-1/2,+1/2)$ on both triangular families. At infinity the labeled eigenlines have these degrees on $\Gamma_c$ and the reverse degrees on $\Gamma_4$. For $\Gamma_0$ and $\Gamma_1$, the calculation for a varying flag preceding~\eqref{eq:A2flagperiod} gives vanishing discriminant and normalized degrees $(+1,-1)$ at the pole where the flag varies; all other eigenline bundle degrees vanish. Consequently
\begin{align}
 (r_{\Gamma_c},Q_{\Gamma_c})&=(x,Q),\notag\\
 (r_{\Gamma_i},Q_{\Gamma_i})&=(y_i,Q_i),\quad i=0,1,
 \label{eq:D3periods}\\
 (r_{\Gamma_4},Q_{\Gamma_4})
   &=(1-x-y_0-y_1,-Q-Q_0-Q_1).\notag
\end{align}
Their sum gives $(r_F,Q_F)=(1,0)$.

\subsubsection{The intersection lattice and admissible parameters}

At $m=c_0=c_1=0$, the quartic is $P_u=p^2z^2(z-1)^2+uz(z-1)$, whose discriminant is $p^2(p^2-4u)u^4$. There is one additional $I_1$ fiber at $u=p^2/4$, and the minimal rational elliptic compactification has a fiber of type $I_1^*$ at infinity. A section is
\begin{equation*}
 \begin{pmatrix}P&uz(z-1)\\1&-P\end{pmatrix}\frac{\dd z}{z(z-1)}
\end{equation*}
on $\cO(-1)\oplus\cO(-2)$, with the fiber of the second summand as the parabolic line at $z=0,1$. Indeed, at both points the residue is nonzero nilpotent with that kernel. At $u=0$ the section lies in $\Gamma_c$; at every other $u$ the spectral curve is irreducible, including the nodal case. Since the section meets $F$ once, $F$ is primitive. The embedded spheres $\Gamma_0,\Gamma_c,\Gamma_1$ form a chain with Gram matrix
\begin{equation*}
 \begin{pmatrix}-2&1&0\\1&-2&1\\0&1&-2\end{pmatrix},
 \qquad |\det|=4.
\end{equation*}
This equals the determinant of the quotient lattice for a fiber of type $I_1^*$ at infinity in~\eqref{eq:quotient-lattices}. Equation~\eqref{eq:lattice-index} gives
\begin{equation*}
 H_2=\Z F\oplus\Z\Gamma_0\oplus\Z\Gamma_c\oplus\Z\Gamma_1.
\end{equation*}
The quotient is the negative $A_3$ lattice, whose roots are the signed sums of consecutive simple roots. Adding $nF$ to the corresponding sums in integral homology gives the classes whose periods must be nonzero, as expressed by the six conditions in~\eqref{eq:D3-exclusions}.

A global rational eigenvalue one-form must be
\begin{equation*}
 \rho=\left(p+\frac{\varepsilon_0c_0}{z}
                  +\frac{\varepsilon_1c_1}{z-1}\right)\dd z,
 \qquad\varepsilon_i\in\{-1,1\}.
\end{equation*}
Comparing the square of this one-form with the quadratic differential in~\eqref{eq:familyD3} gives
\begin{equation*}
 m=\varepsilon_0c_0+\varepsilon_1c_1,
 \qquad u=2\varepsilon_1c_1(p+\varepsilon_0c_0).
\end{equation*}
For $n_i=(1-\varepsilon_i)/2$, an invariant line subbundle of parabolic degree zero can occur for the splitting with these signs precisely when $x+n_0y_0+n_1y_1$ is integral, giving the four nonvanishing conditions involving $Q$. The other two conditions concern the tame poles. If $c_i=0$, then $0<y_i<1$ and hence $0<\alpha_i<1/2$. At $y_i=0$ or $1$ the residue $c_i$ is nonzero, so the off-diagonal indicial operator has a positive term. The quotient is smooth and its metric is complete.

At fixed $p\ne0$, contract the alcove coordinates to an interior point, then contract $(m,c_0,c_1)$. At a splitting, $Q+n_0Q_0+n_1Q_1=0$. If $c_i=0$, then $Q_i=0$ and either choice of $n_i$ satisfies this equation. The alcove inequalities and~\eqref{eq:D3-exclusions} thus put every corresponding $x+n_0y_0+n_1y_1$ in $(0,1)$. Hence all compatible flags give the same strict integer threshold locally. Section~\ref{subsec:topology-transport} gives local properness with these thresholds. At $\alpha_\infty=0,1$ the distinct leading eigenvalues give the common gauge domains of Lemma~\ref{lem:weight-domains}, and Section~\ref{subsec:calibration} preserves the periods under normalization. Lemma~\ref{lem:marking-transport} now extends the integral basis and period formulas over the contractible domain.

\subsubsection{Affine Weyl action}

The chain $\Gamma_0,\Gamma_c,\Gamma_1$ and the fourth component $\Gamma_4=F-\Gamma_0-\Gamma_c-\Gamma_1$ of the type $I_4$ fiber realize the four affine simple reflections. Their compactly supported Dehn twists induce the following reflections on $(y_0,x,y_1)$:
\begin{gather*}
 (-y_0,x+y_0,y_1),\qquad(y_0+x,-x,y_1+x),\\
 (y_0,x+y_1,-y_1),\qquad(1-x-y_1,x,1-x-y_0).
\end{gather*}
The reflections generate the affine Weyl group of type $A_3$. Each orbit of the real period vector meets the alcove~\eqref{eq:D3alcove}. The action on the complex periods comes from the same homology reflections, with $Q_F=0$ in place of $r_F=1$. The weights and residues in this marking are computed in Section~\ref{subsec:inverse-periods}.

\subsection{Weights and residues from compact periods}
\label{subsec:inverse-periods}

Use a common affine Weyl marking for real and complex periods in $A_1,A_2,D_3$. For $D_1,D_2$ the table uses unnormalized real weights and their induced integral markings.

{\small
\setlength{\tabcolsep}{4pt}
\renewcommand{\arraystretch}{1.25}
\begin{longtable}{@{}p{.08\textwidth}p{.43\textwidth}p{.43\textwidth}@{}}
\caption{Compact periods and inverse formulas.}\label{tab:compact-periods}\\
\toprule
Type & Integral basis and normalized periods & Weights and residues\\
\midrule
\endfirsthead
\multicolumn{3}{l}{\small Table~\thetable\ continued.}\\
\toprule
Type & Integral basis and normalized periods & Weights and residues\\
\midrule
\endhead
\bottomrule
\endfoot
$A_0$ & $F$; no additional compact period. & No weight or residue parameter to recover.\\
\midrule
$A_1$ & $(F,\Gamma)$, $\Gamma^2=-2$;\newline $(r,Q)=(\alpha,ib)$. & \refstepcounter{equation}\label{eq:A1-inverse}$\alpha=r,\quad b=-iQ$.\hfill(\theequation)\\*
\multicolumn{3}{@{}p{.96\textwidth}@{}}{\footnotesize Domain: $0\le r\le1$, with $(r,Q)\ne(0,0),(1,0)$.}\\
\midrule
$A_2$ & $(F,\Gamma_1,\Gamma_2)$;\newline $(x,y)=(\alpha_0+\alpha_\infty-1,1-2\alpha_0)$,\newline $(Q_1,Q_2)=(2i(\mu-c),4ic)$. & \refstepcounter{equation}\label{eq:A2inverse}$\begin{aligned}[t]
\alpha_0&=(1-y)/2,\\
\alpha_\infty&=x+(1+y)/2,\\
c&=-iQ_2/4,\\
\mu&=-i(2Q_1+Q_2)/4.
\end{aligned}$\par\hfill(\theequation)\\*
\multicolumn{3}{@{}p{.96\textwidth}@{}}{\footnotesize Domain: $x,y\ge0$, $x+y\le1$; the pairs $(x,Q_1)$, $(y,Q_2)$, $(x+y,Q_1+Q_2)$ avoid $\Z\times\{0\}$. The Gram matrix of $(\Gamma_1,\Gamma_2)$ is the negative $A_2$ matrix.}\\
\midrule
$D_0$ & $F$; no additional compact period. & No weight or residue parameter to recover.\\
\midrule
$D_1$ & $(F,\Gamma)$, $\Gamma^2=-4$, $\Gamma\cdot O=0$, $\int_\Gamma c_1(L_+^0)=1$;\newline $(r,Q)=(1-2\alpha,4im)$. & \refstepcounter{equation}\label{eq:D1-inverse}$\alpha=(1-r)/2,\quad m=-iQ/4$.\par\hfill(\theequation)\\*
\multicolumn{3}{@{}p{.96\textwidth}@{}}{\footnotesize Domain: $(r,Q)\in\R\times\C$, including $(0,0)$. This lattice has no class of square $-2$.}\\
\midrule
$D_2$ & $(F,\Gamma_1,\Gamma_2)$;\newline $x=\alpha_0+\alpha_\infty-1$, $y=\alpha_0-\alpha_\infty$,\newline $Q_1=2i(m_\infty-m_0)$,\newline $Q_2=-2i(m_0+m_\infty)$. & \refstepcounter{equation}\label{eq:D2-inverse}$\begin{aligned}[t]
\alpha_0&=(x+y+1)/2,\\
\alpha_\infty&=(x-y+1)/2,\\
m_0&=i(Q_1+Q_2)/4,\\
m_\infty&=i(Q_2-Q_1)/4.
\end{aligned}$\par\hfill(\theequation)\\*
\multicolumn{3}{@{}p{.96\textwidth}@{}}{\footnotesize Domain: $x,y\in\R$ and $(x,Q_1),(y,Q_2)\notin\Z\times\{0\}$. The Gram matrix of $(\Gamma_1,\Gamma_2)$ is $\diag(-2,-2)$.}\\
\midrule
$D_3$ & $(F,\Gamma_0,\Gamma_c,\Gamma_1)$ in chain order;\newline $(r_{\Gamma_0},r_{\Gamma_c},r_{\Gamma_1})=(y_0,x,y_1)$,\newline $(Q_{\Gamma_0},Q_{\Gamma_c},Q_{\Gamma_1})=(Q_0,Q,Q_1)$. & \refstepcounter{equation}\label{eq:D3inverse}$\begin{aligned}[t]
\alpha_i&=(1-y_i)/2,\quad i=0,1,\\
\alpha_\infty&=x+(y_0+y_1)/2,\\
c_i&=-iQ_i/4,\\
m&=-i(2Q+Q_0+Q_1)/4.
\end{aligned}$\par\hfill(\theequation)\\*
\multicolumn{3}{@{}p{.96\textwidth}@{}}{\footnotesize Here $x=\alpha_0+\alpha_1+\alpha_\infty-1$, $y_i=1-2\alpha_i$, $Q=2i(m-c_0-c_1)$ and $Q_i=4ic_i$. The real domain is $x,y_0,y_1\ge0$, $x+y_0+y_1\le1$; the six pairs in~\eqref{eq:D3-exclusions} avoid $\Z\times\{0\}$. The displayed chain has negative $A_3$ Gram matrix.}\\
\end{longtable}
}
{\footnotesize We use the normalization~\eqref{eq:compact-normalization}, with $(r_F,Q_F)=(1,0)$. Each excluded pair rules out zero period on an integral translate of a root by $F$. The weights specify the chosen lattices, so $\alpha_\infty$ may lie outside $[0,1]$.\par}

For $A_2,D_3$, the exclusions imply $0<\alpha_i<1/2$ if $c_i=0$, as required in~\eqref{eq:tame-domain}.

\section{Spectral periods and the parameters at infinity}
\label{sec:periodend}

The period ratio determines the logarithmic potential in the four $D$ families. In the three $A$ families, two linear combinations of the spectral periods determine the cone coordinate and the coefficient of the deformation at infinity.

\subsection{Logarithmic asymptotics of the period ratio}

For a quartic $P$, the invariants introduced in Section~\ref{subsec:topology-transport} give $j=256I(P)^3/\operatorname{disc}_z(P)$. For a cubic in Weierstrass form we use $j=c_4^3/\Delta$, where $c_4$ and $\Delta$ are the standard Weierstrass invariant and discriminant. Both formulas use the normalization $j(\tau)=e^{-2\pi i\tau}+744+O(e^{2\pi i\tau})$ at the cusp.

\subsubsection{Type $D_0$}

Setting $\eta=z^2\lambda/\dd z$, the spectral curve takes the form
\begin{equation}
 \eta^2=p_\infty z^3+uz^2+p_0z,\qquad
 \lambda=\eta\frac{\dd z}{z^2},\qquad
 \partial_u\lambda=\frac{\dd z}{2\eta}.
 \label{eq:D0-cubic}
\end{equation}
With $X=p_\infty z$, $Y=p_\infty\eta$, and $h=p_0p_\infty$, the equation is $Y^2=X^3+uX^2+hX$. Its invariants are
\begin{equation}
 c_4=16(u^2-3h),\qquad
 \Delta=16h^2(u^2-4h),\qquad
 j=\frac{256(u^2-3h)^3}{h^2(u^2-4h)}.
 \label{eq:D0-invariants}
\end{equation}
The expansion of the elliptic modular function at the cusp gives
\begin{equation}
 \tau=\frac{2i}{\pi}\log\frac{4u}{\sqrt h}+O(u^{-2}).
 \label{eq:D0-tau}
\end{equation}

\subsubsection{Type $D_1$}

For the Weierstrass equation~\eqref{eq:D1-Weierstrass},
\begin{align}
 c_4&=16(u^2-6\lambda_0 m),\notag\\
 \Delta&=-16\lambda_0^2(4u^3-4m^2u^2-36\lambda_0 mu+32\lambda_0 m^3+27\lambda_0^2).
 \label{eq:D1-invariants}
\end{align}
It follows that
\begin{align}
 j&=-\frac{64u^3}{\lambda_0^2}\left(1+\frac{m^2}{u}+O(u^{-2})\right),\notag\\
 \tau&=\frac{i}{2\pi}\log\frac{-64u^3}{\lambda_0^2}
       +\frac{im^2}{2\pi u}+O(u^{-2}).
 \label{eq:D1-j-tau}
\end{align}

\subsubsection{Type $D_2$}

The quartic polynomial $P$ in~\eqref{eq:D2-quartic} has
\begin{align*}
 I(P)&=u^2+12h^2-12hm_0m_\infty,\\
 \operatorname{disc}(P)&=16h^2u^4-16h^2(m_0^2+m_\infty^2)u^3+O(u^2).
\end{align*}
Indeed, the term of order $u^4$ in the discriminant comes from $16AC^4E$, and the terms of order $u^3$ come from $-4AC^3D^2-4B^2C^3E$. Thus
\begin{align}
 j&=\frac{16u^2}{h^2}\left(1+\frac{m_0^2+m_\infty^2}{u}+O(u^{-2})\right),\notag\\
 \tau&=\frac i\pi\log\frac{4u}{h}
       +\frac{i(m_0^2+m_\infty^2)}{2\pi u}+O(u^{-2}).
 \label{eq:D2-tau}
\end{align}

\subsubsection{Type $D_3$}

Multiplying the coefficient of the spectral quadratic differential in the given coordinate by $z^2(z-1)^2$ gives the quartic
\begin{equation*}
 P_u=p^2z^2(z-1)^2+2pmz(z-1)^2+uz(z-1)+c_0^2(z-1)^2+c_1^2z^2.
\end{equation*}
Its coefficients are
\begin{equation*}
 (A,B,C,D,E)=\bigl(p^2,-2p^2+2pm,u+p^2-4pm+c_0^2+c_1^2,-u+2pm-2c_0^2,c_0^2\bigr).
\end{equation*}
These coefficients give $I(P_u)=u^2+O(u)$. Only the discriminant term $-4AC^3D^2$ contributes at order $u^5$, and its leading term is $-4p^2u^5$. Thus, for arbitrary fixed $m,c_0,c_1$,

\begin{equation}
 j=-\frac{64u}{p^2}(1+O(u^{-1})),\qquad
 \tau=\frac{i}{2\pi}\log\frac{-64u}{p^2}+O(u^{-1}).
 \label{eq:D3-tau}
\end{equation}

Choose a fixed circle enclosing the two inner branch points and excluding the other two for large $|u|$. In $D_0,D_1,D_2$ we may take a circle about zero. In $D_3$ the two inner branch values tend to $0,1$, while the outer branch values have modulus asymptotic to $\sqrt{|u|}/|p|$. A fixed circle of sufficiently large radius $R>2$ in the separating annulus $1\ll |z|\ll\sqrt{|u|}$ encloses the inner pair, uniformly for parameters in a compact admissible set. The two lifts of the circle are parallel simple essential curves: the double cover of either complementary disk branched at two points is an annulus. Choose a lift $A$ and its orientation so that the leading differential on $A$ is $\sqrt u\,\dd z/z$ in the first three cases and $\sqrt u\,\dd z/\sqrt{z(z-1)}$ in $D_3$, with $\sqrt{z(z-1)}\sim z$ at infinity. The integral of each leading differential is $2\pi i\sqrt u$. Expanding the square root on this fixed contour gives a convergent power series in $u^{-1}$ after $\sqrt u$ is factored out. In particular,
\[
 Z_A=2\pi i\sqrt u\,(1+O(u^{-1})).
\]
This expansion can be differentiated term by term and is uniform on compact parameter sets. In the normalization of~\eqref{eq:factor2}, put $\zeta=-2i Z_A$ on the double cover and $r=|\zeta|$. Then
\begin{equation}
 r\sim4\pi\sqrt{|u|},\qquad V=\Ima\tau=\kappa_0+\frac\nu\pi\log r+o(1).
 \label{eq:logarithmic-radius}
\end{equation}
Substituting this asymptotic expression for the radius into the four period-ratio expansions gives the constants $\kappa_0$ and the positive choices of leading parameters listed in Table~\ref{tab:end-parameters}.
For instance, in $D_2$, $\Ima\tau=\pi^{-1}\log|4u/h|+o(1)$ and $\log|u|=2\log r-2\log(4\pi)+o(1)$, giving $\kappa_0=-\pi^{-1}\log(4\pi^2|h|)$. The other three cusp expansions give the remaining entries. In the parallel angles of~\eqref{eq:factor2} the semiflat metric is
\[
 g_{\sfm}=(\Ima\tau)\bigl(|\dd(2Z_A)|^2+\dd\theta_e^2\bigr)
     +(\Ima\tau)^{-1}\bigl(\dd\theta_m-(\Rea\tau)\dd\theta_e\bigr)^2 .
\]
Theorem~\ref{thm:comparison} therefore gives the same logarithmic coefficient and radial constant for the $L^2$ metric. In the radial coordinate above, varying the arguments of the leading parameters changes the constant term in $\Rea\tau$. It leaves the leading term $2\pi i\sqrt u$ of $Z_A$ fixed. On a sector, write
\[
 \Rea\tau=-\frac{\nu}{\pi}\arg\zeta+t+O(r^{-2}),\qquad
 \Ima\tau=\kappa_0+\frac{\nu}{\pi}\log r+O(r^{-2}).
\]
The choice $B\mapsto B+nA$ adds $n$ to $t$. For the positive leading parameters in the table, $t=0\pmod{\Z}$ in $D_0,D_2$ and $t=1/2\pmod{\Z}$ in $D_1,D_3$. Take parallel angle coordinates with the origin specified below, choose a lift of $t$, and set
\[
 \theta_1=\arg\zeta-\pi t/\nu,\qquad
 \theta_2=\theta_e,\qquad
 \theta_3=(\pi/\nu)\theta_m+\theta_1\theta_e .
\]
Adding $2\pi$ to $\theta_e$ adds $2\pi\theta_1$ to $\theta_3$, and adding $2\pi$ to $\theta_m$ adds $2\pi^2/\nu$ to $\theta_3$. On one turn in the original $u$-plane the marked angles return as $(-\theta_e,-\theta_m+\nu\theta_e)$. Hence the deck involution in these coordinates is $(\theta_1,\theta_2,\theta_3)\mapsto(\theta_1+\pi,-\theta_2,-\theta_3)$. With $\Theta=(\nu/\pi)(\dd\theta_3-\theta_2\dd\theta_1)$, one has $\dd\theta_m-(\Rea\tau)\dd\theta_e=\Theta+O(r^{-2})\dd\theta_2$. Equation~\eqref{eq:factor2} reads
\[
 \Omega_I=\dd\zeta\wedge\bigl((\Ima\tau)\dd\theta_e+
                         i(\dd\theta_m-(\Rea\tau)\dd\theta_e)\bigr).
\]
Thus the ordered model in~\cite[Section~2.1]{CV24} uses the base variable $r e^{i\theta_1}=e^{-i\pi t/\nu}\zeta$ and the same constant phase $e^{-i\pi t/\nu}$ on $\Omega_I$. This rotation about $\omega_I$ depends only on the chosen leading parameters and integral frame. In applying Torelli, periods initially given in the standard model frame are entered in Table~\ref{tab:compact-periods} as $Q_\Gamma^{\mathrm{sp}}=e^{i\pi t/\nu}Q_\Gamma^{\mathrm{CV}}$. Both ordered triples are then compared in the spectral frame; applying the common rotation recovers the standard frame. The differentiated expansions above give positive-order decay in the model distance $r\sqrt{\log r}$; factors of $\log r$ can be absorbed by decreasing the exponent.

\subsection{Periods of the spectral differential in the three ALG families}

The cycles below are primitive integral cycles on a smooth spectral curve, oriented so that $A\cdot B=1$. After subtracting the $u$-independent additive constants from their periods and choosing branches of the fractional powers of $u$ on a sector, we define $U,W$ as linear combinations of the spectral periods. These functions depend on the base coordinate and deformation parameters, while $v$ is a complex-valued linear combination of the angular coordinates on the universal cover of each torus fiber. The combinations specified below for each type satisfy
\begin{equation}
 \Omega_I=\dd U\wedge\dd v+\dd W\wedge\dd\bar v,
 \qquad U=\ell_j u^{\beta_j}+\cdots,\qquad W=-d_j a u^{1-\beta_j}+\cdots.
 \label{eq:normalized-periods}
\end{equation}
We normalize these combinations so that the coefficients $\ell_j,d_j$ are positive. A subscript $a$ or $u$ on a spectral period denotes differentiation with respect to that parameter with the cycle transported by the Gauss--Manin connection.

\subsubsection{Type $A_0$}

With $\omega=e^{2\pi i/3}$ and~$\tau_2=\Ima\omega=\sqrt3/2$, we use the integral
\begin{equation}
 I_{A_0}=\int_0^1(1-s^3)^{1/2}\dd s=\frac13B(1/3,3/2)\,.
 \label{eq:A0-beta-integrals}
\end{equation}
At~$a=0$ we use the scaled coordinate $z/u^{1/3}$. The cycle~$A$ follows the path from~$\omega e^{\pi i/3}$ through the origin to~$e^{\pi i/3}$ on one sheet and returns on the other; rotating it by~$\omega$ gives~$B$. The corresponding spectral periods are
\begin{equation}
 Z_A=2e^{\pi i/3}(1-\omega)I_{A_0}u^{5/6}\,,\qquad
 Z_B=\omega Z_A\,.
 \label{eq:A0-period-integrals}
\end{equation}
The rotation multiplies $\partial_a\lambda=z\dd z/(2y)$ by~$\omega^2$, so $Z_{B,a}=\bar\omega Z_{A,a}$ at~$a=0$.

First set
\begin{equation}
 v=\frac{\theta_m-\omega\theta_e}{\sqrt{\tau_2}}\,,\quad
 U=\frac{Z_B-\bar\omega Z_A}{i\sqrt{\tau_2}}\,,\quad
 W=\frac{\omega Z_A-Z_B}{i\sqrt{\tau_2}}\,.
 \label{eq:A0-period-frame}
\end{equation}
Equation~\eqref{eq:factor2} reads $\Omega_I=\dd U\wedge\dd v+\dd W\wedge\dd\bar v$. We then replace the coordinates $U,v,W$ by $e^{-\pi i/6}U,e^{\pi i/6}v,e^{\pi i/6}W$. There is no term linear in $a$ in the expansion of~$U$ at the reference curve. Homogeneity and Taylor expansion in~$au^{-2/3}$ give
\begin{equation}
 U=\ell_0u^{5/6}+O(a^2u^{-1/2})\,,\qquad
 W=-d_0a u^{1/6}+O(a^2u^{-1/2})\,,
 \label{eq:A0-leading-periods}
\end{equation}
where
\begin{equation}
 \ell_0=4\sqrt{3\tau_2}\,I_{A_0}>0\,.
 \label{eq:A0-LM}
\end{equation}

\subsubsection{Type $A_1$}

At $a=b=0$, we use the coordinate
\begin{equation}
 z=e^{\pi i/4}u^{1/4}x\,.
\end{equation}
We choose continuous branches of the fractional powers on a sector in the $u$-plane. The rescaled spectral curve is $y^2=u(1-x^4)$, and the leading period integrals reduce to
\begin{equation}
 I_{A_1}=\int_0^1\sqrt{1-x^4}\,\dd x
     =\frac14 B(1/4,3/2)\,.
 \label{eq:I0}
\end{equation}
On one sheet the cycle $A$ follows the path from $x=i$ to $x=1$ through the origin, and the cycle $B$ follows the path from $x=-1$ to $x=i$ through the origin. Each returns along the same path on the other sheet. With these orientations we obtain
\begin{align}
 Z_A&=2e^{\pi i/4}(1-i)I_{A_1}\,u^{3/4}
     =2\sqrt2 I_{A_1}\,u^{3/4}\,,\label{eq:ZA0}\\
 Z_B&=2e^{\pi i/4}(1+i)I_{A_1}\,u^{3/4}
     =2i\sqrt2 I_{A_1}\,u^{3/4}\,.\label{eq:ZB0}
\end{align}
The lifted automorphism $(x,y)\mapsto(ix,y)$ sends $A$ to $B$ and multiplies $\partial_a\lambda=z^2\dd z/(2y)$ by~$-i$. Hence $Z_{B,a}=-iZ_{A,a}$ at $a=b=0$.

The beta-integral expression for the complete elliptic integral gives
\begin{equation}
 \ell_1=4\sqrt2I_{A_1}=\frac{2\Gamma(1/4)^2}{3\sqrt\pi}=\frac83\mathbf K(1/\sqrt2)\,,
 \label{eq:A1-integral-constants}
\end{equation}
where $\mathbf K$ denotes the complete elliptic integral of the first kind, with modulus as argument. With
\begin{equation}
 U=Z_A-iZ_B,\qquad W=Z_A+iZ_B,\qquad v=\theta_m-i\theta_e\,,
 \label{eq:UW}
\end{equation}
the coordinates have the normalization in~\eqref{eq:normalized-periods} and the cone parameter is $\beta_1=3/4$.

Factoring out $u^{3/4}$ makes the periods analytic in the scaled parameters $(au^{-1/2},bu^{-3/4})$ near the smooth reference quartic. The term linear in $a$ cancels in $U$, and the linear $b$-terms in the spectral periods are additive constants. After subtracting these constants, the remainder is $O(|u|^{-1/4})$, with the $a^2$ term contributing at this order. With $\rho\asymp|u|^{3/4}$, this gives, uniformly for $(a,b)$ in compact subsets of $\C^2$,
\begin{equation}
 U=\ell_1u^{3/4}+\Op(\rho^{-1/3}),\qquad
 W=-d_1a u^{1/4}+\Op(\rho^{-1/3}).
 \label{eq:A1-leading-periods}
\end{equation}

\subsubsection{Type $A_2$}

At $a=\mu=c=0$, the scaled coordinates and constants are
\begin{equation*}
 z=-u^{1/3}x,\qquad\eta=iu^{2/3}Y,
 \quad Y^2=x(1-x^3),\quad
 \tau_0=e^{\pi i/3},\quad\xi=e^{2\pi i/3},\quad \tau_2=\frac{\sqrt3}{2}\,.
\end{equation*}
The cycles $A$ and $B$ lie over $[0,1]$ and $[0,\xi]$, respectively. Each runs from $0$ to the other endpoint on one sheet and returns on the other. In the local coordinate $\zeta=\sqrt x$ at the ramification point over zero, the tangent directions of the lifted cycles are $1$ and $e^{\pi i/3}$. Hence $A\cdot B=1$, and both cycles are primitive integral cycles. Direct integration gives
\begin{equation*}
 Z_A=2iI_{A_2}u^{2/3},\qquad Z_B=\tau_0Z_A,\qquad
 I_{A_2}=\int_0^1x^{-1/2}(1-x^3)^{1/2}\dd x
   =\frac13B(1/6,3/2)\,.
\end{equation*}
In this family $\partial_a\lambda=(z^2+az/2)\dd z/(2\eta)$. On the reference curve the lifted automorphism $(x,Y)\mapsto(\xi x,\tau_0Y)$ sends $A$ to $B$ and multiplies $z^2\dd z/(2\eta)$ by~$\bar\tau_0$. Thus $Z_{B,a}=\bar\tau_0Z_{A,a}$ at the reference curve. The normalized coordinates are
\begin{equation*}
 U=\frac{\bar\tau_0 Z_A-Z_B}{\sqrt \tau_2},\quad
 W=\frac{\tau_0 Z_A-Z_B}{\sqrt \tau_2},\quad
 v=\frac{i(\theta_m-\tau_0\theta_e)}{\sqrt \tau_2}\,.
\end{equation*}
Then
\begin{equation*}
 \Omega_I=\dd U\wedge\dd v+\dd W\wedge\dd\bar v,
 \qquad\frac i2\dd v\wedge\dd\bar v
       =\dd\theta_e\wedge\dd\theta_m,
\end{equation*}
and
\begin{equation*}
 U=\ell_2u^{2/3}+\cdots,
 \qquad W=-d_2a u^{1/3}+\cdots,
 \qquad \ell_2=4\sqrt \tau_2\,I_{A_2}>0\,.
\end{equation*}
To estimate the remainders, fix a branch $t=u^{-1/3}$ and set $P=at$, $M=\mu t^2$, $C=ct^2$. In the same scaled coordinates the full family is
\[
 Y^2=x(1-x^3)+Px^3-(P^2/4+2M)x^2-C^2,
 \qquad t^2\lambda=iY\,\dd x/x .
\]
Represent $A$ and $B$ by contours on the reference curve avoiding the points over $x=0,\infty$ and all ramification points. Since the four limiting roots are simple, the projected contours still avoid them for small $(P,M,C)$, and their closed lifts continue with the same choice of sheet. Their integrals of $t^2\lambda$ are holomorphic in $(P,M,C^2)$. Other parallel lifts to the punctured curve can add peripheral periods, linear in $M,C$ after this scaling.

The rotation above gives the linear terms in $P$. Terms of weight two in $t$ give constants in the unscaled periods. The possible $a^2$ coefficient is zero, since on $\eta^2=z^4+uz$ one has
\[
 \left.\partial_a^2\lambda\right|_{a=\mu=c=0}
 =\frac{uz^2}{4\eta^3}\dd z=\frac16\dd(z^2/\eta).
\]
The term proportional to $\mu$ can survive. On the same curve,
\[
 \left.\partial_\mu\lambda\right|_{a=\mu=c=0}
 =\frac{z\,\dd z}{\eta}
 =-\frac13\,\dd\log\!\left(\frac{\eta-z^2}{\eta+z^2}\right).
\]
We subtract these constant coefficients, including any peripheral contributions, as affine translations of the period coordinates. They are independent of $a$ and $u$, so the period derivatives and~\eqref{eq:bilinear} are unchanged. All remaining terms have weight at least three. For $\rho\asymp |u|^{2/3}$ this gives
\begin{equation}
 U=\ell_2u^{2/3}+\Op(\rho^{-1/2}),\qquad
 W=-d_2a u^{1/3}+\Op(\rho^{-1/2}).
 \label{eq:A2-period-remainders}
\end{equation}
On closed subsectors, the convergent Taylor series may be differentiated termwise. Each derivative in the flat base coordinate $u^{2/3}$ costs one power of $\rho$, uniformly for parameters in compact sets. Dividing $\partial_uW$ by $\partial_uU$ gives
\[
 W_U=-\frac{d_2a}{2\ell_2}u^{-1/3}+\Op(\rho^{-3/2}).
\]
When $a=0$, only even powers of $t$ occur before multiplication by $t^{-2}$. After the constant terms have been removed, both remainders in~\eqref{eq:A2-period-remainders} are $\Op(\rho^{-1})$, and $W_U=\Op(\rho^{-2})$.

In particular $S=U-\bar W=\ell_2u^{2/3}+\Op(\rho^{1/2})$. The correction divided by $\ell_2u^{2/3}$ is $\Op(\rho^{-1/2})$, and the derivative of the change from $\ell_2u^{2/3}$ to $S$ differs from the identity by $\Op(\rho^{-1/2})$. Inverting on a closed subsector and using the corresponding branch of the power gives
\begin{equation}
 u=(S/\ell_2)^{3/2}+\Op(\rho).
 \label{eq:A2-base-inverse}
\end{equation}
Taking $B'=B-A$ as the second basis cycle gives the same hexagonal lattice with period ratio $e^{2\pi i/3}$. If we also replace $\theta_m$ by $\theta_m-\theta_e$, the coordinates $U,W,v$ remain unchanged.

\subsection{The Riemann bilinear relation}

\begin{proposition}
\label{prop:bilinear}
For $A_j$, $j=0,1,2$, extend the symplectic basis $(A,B)$ above by parallel transport over a simply connected subset of the complement of the discriminant. Use the punctured-cycle lifts of Section~\ref{sec:families} for the periods. Then the spectral periods and the associated normalized coordinates satisfy
\begin{equation}
 Z_{A,a}Z_{B,u}-Z_{B,a}Z_{A,u}=-\frac{iK_j}{4},\qquad
 \left.\partial_aW\right|_U=-\frac{K_j}{2}\partial_Uu.
 \label{eq:bilinear}
\end{equation}
\end{proposition}

\begin{proof}
We write $\eta_u=\partial_u\lambda$ and $\eta_a=\partial_a\lambda$, and let $f_a$ denote a local primitive of $\eta_a$ at each pole. In $A_1$ and $A_2$, the sign $\varepsilon=\pm1$ of the eigenvalue one-form distinguishes the two points over infinity. At these points the two differentials have expansions
\begin{equation*}
 \eta_a=-\frac\varepsilon2w^{-2}\dd w+O(1)\dd w,
 \qquad\eta_u=-\frac\varepsilon2\dd w+O(w)\dd w.
\end{equation*}
Taking $f_a=\varepsilon/(2w)+O(w)$ gives $\Res(f_a\eta_u)=-1/4$ at each point. In $A_0$, the coordinate $z=\zeta^{-2}$ at the single spectral point at infinity gives
\begin{equation*}
 \eta_u=-\dd \zeta+O(\zeta^2)\dd \zeta,\qquad
 \eta_a=-\zeta^{-2}\dd \zeta+O(1)\dd \zeta.
\end{equation*}
Here $f_a=\zeta^{-1}+O(\zeta)$, so $\Res(f_a\eta_u)=-1$ at this point. There are no other contributions, since both differentials are regular at the ordinary ramification points. In $A_2$, $\eta_u=\dd z/(2\eta)$ and $\eta_a=(z^2+az/2)\dd z/(2\eta)$ are regular also above the tame pole, using $z=s^2$ if that point is ramified. By the Riemann bilinear relation, the period pairing is $-\pi i$ for $A_1,A_2$ and $-2\pi i$ for $A_0$, which is the first identity in~\eqref{eq:bilinear}.

On each smooth elliptic spectral curve, $\eta_u$ is a nonzero holomorphic differential, so $Z_{A,u}\ne0$ and $\tau(u)=Z_{B,u}/Z_{A,u}$ lies in the upper half-plane. In each of the three coordinate systems, $U_u$ is a nonzero constant times $Z_{A,u}(\tau(u)-\bar\tau_0)$, where $\tau_0$ is the limiting period ratio in the chosen cycle basis. Since $\Ima\tau_0>0$, we have $U_u\ne0$, and $U$ is a local coordinate on the complement of the discriminant in the Hitchin base.

Each normalized linear change from $(Z_A,Z_B)$ to $(U,W)$ has determinant $2i$. Therefore $U_aW_u-U_uW_a=K_j/2$. With $U$ held constant,
\begin{equation*}
 \left.\partial_aW\right|_U
 =W_a-W_u\frac{U_a}{U_u}
 =-\frac{K_j}{2U_u}
 =-\frac{K_j}{2}\partial_Uu.
\end{equation*}
\end{proof}

For each of the three reference curves with vanishing deformation parameters, $U=\ell_j u^{\beta_j}$ and $W=0$. Homogeneity gives $W_a=-d_j u^{1-\beta_j}$ for a coefficient~$d_j$. Since $U_u=\beta_j\ell_j u^{\beta_j-1}$ and $W_u=0$, the identity $U_aW_u-U_uW_a=K_j/2$ gives $\beta_j\ell_jd_j=K_j/2$. Hence
\begin{equation}
 d_j=\frac{K_j}{2\beta_j\ell_j},\qquad
 \ell_jd_j=\frac{K_j}{2\beta_j}.
 \label{eq:period-product}
\end{equation}

\begin{table}[ht]
\centering
\small
\setlength{\tabcolsep}{5pt}
\begin{tabular}{@{}ccclc@{}}
\toprule
Type & $\beta_j$ & $\tau_0$ & $\ell_j$ & $K_j$\\
\midrule
$A_0$ & $5/6$ & $e^{2\pi i/3}$ & $\frac43\sqrt{3\sqrt3/2}\,B(1/3,3/2)$ & $8\pi$\\[3pt]
$A_1$ & $3/4$ & $i$ & $\frac83\mathbf K(1/\sqrt2)$ & $4\pi$\\[3pt]
$A_2$ & $2/3$ & $e^{\pi i/3}$ & $\frac43\sqrt{\sqrt3/2}\,B(1/6,3/2)$ & $4\pi$\\
\bottomrule
\end{tabular}
\par\medskip
\begin{tabular}{@{}ccll@{}}
\toprule
Type & $\nu$ & $\kappa_0$ & Positive leading parameter\\
\midrule
$D_0$ & $4$ & $-\pi^{-1}\log(16\pi^4|h|)$ & $h=e^{-\pi\kappa_0}/(16\pi^4)$\\
$D_1$ & $3$ & $-\pi^{-1}\log(8\pi^3|\lambda_0|)$ & $\lambda_0=e^{-\pi\kappa_0}/(8\pi^3)$\\
$D_2$ & $2$ & $-\pi^{-1}\log(4\pi^2|h|)$ & $h=e^{-\pi\kappa_0}/(4\pi^2)$\\
$D_3$ & $1$ & $\pi^{-1}\log(2/(\pi|p|))$ & $p=(2/\pi)e^{-\pi\kappa_0}$\\
\bottomrule
\end{tabular}
\caption{The parameters at infinity and their normalizations.}
\label{tab:end-parameters}
\par\smallskip
\begin{minipage}{.97\textwidth}
\footnotesize
The ALG rows use the cycles of Section~\ref{sec:periodend} and the ALG$^*$ rows use the radius of~\eqref{eq:logarithmic-radius} and $L=1$. We take $p_0=p_\infty=\sqrt h$ in $D_0,D_2$ and $p=1$, $b=\lambda_0$ in $D_1$.

\end{minipage}
\end{table}

\subsection{The semiflat hyperk\"ahler forms}

We use on each end a single choice of the parallel reference introduced in Section~\ref{sec:families}. Indeed, on the original end annulus the return on its angle torus is affine, $t\mapsto Mt+b$, since the peripheral characters are fixed. For $A_0,A_1,A_2$ the determinants of $\Id-M$ are respectively $1,2,3$. For $D_0,D_1,D_2,D_3$ the $I_\nu^*$ monodromy is conjugate, up to the cycle conventions, to minus a unipotent matrix with both eigenvalues one, so $\det(\Id-M)=4$. The endomorphism $\Id-M$ of the real angle torus is therefore surjective. Choose $t_0$ with $(\Id-M)t_0=b$; transport of this point gives a parallel origin on the entire end. Such choices can be continued along any parameter interval by path lifting for the finite isogeny. In the $D$ cases we make the choice before taking the twofold cover used in the asymptotic model. We use these parallel angles in the forms below and in the end identifications of Sections~\ref{sec:finite} and~\ref{sec:realization}. In the marking transport of Lemma~\ref{lem:marking-transport}, we identify the ends at fixed $u$ and fixed parallel angles along the parameter path. An algebraic section used to present a compactified Picard variety need not serve as this metric origin.

We write $S=U-\bar W$ and $W_U=\partial_UW$, taking the derivative with the parabolic weights and the principal parts of the eigenvalue one-forms held constant. The semiflat hyperk\"ahler forms are
\begin{align}
 \omega_{I,\sfm}&=\frac i2(\dd S\wedge\dd\bar S+\dd v\wedge\dd\bar v),\notag\\
 \omega_{K,\sfm}&=\Ima(\dd S\wedge\dd v),\notag\\
 \omega_{J,\sfm}&=\frac{1+|W_U|^2}{1-|W_U|^2}\Rea(\dd S\wedge\dd v)
 +\frac2{1-|W_U|^2}\Rea(W_U\dd S\wedge\dd\bar v).
 \label{eq:semiflat-forms}
\end{align}
Indeed, substituting $\dd W=W_U\dd U$ and $\dd S=\dd U-\overline{W_U}\dd\bar U$ into~\eqref{eq:normalized-periods} gives the last two forms. For the first form, the identity $\dd W=W_U\dd U$ gives
\begin{equation*}
 \frac i2\dd S\wedge\dd\bar S
 =\frac i2(\dd U\wedge\dd\bar U-\dd W\wedge\dd\bar W).
\end{equation*}
We use $S,v$ as normalized coordinates on the cover of the end in Section~\ref{sec:finite}.

\subsection{The ALG metrics at $a=0$}
\label{subsec:centers}

\subsubsection{Type $A_0$}

At $a=0$, equation~\eqref{eq:A0-period-integrals} gives $Z_B=\omega Z_A$ for every $u\ne0$, so the normalized coordinate $W$ vanishes identically. Thus~\eqref{eq:semiflat-forms} is the flat model, and~\eqref{eq:A0-actual-comparison} with its derivative bounds gives order-$2$ decay.

\subsubsection{Type $A_1$}

We differentiate the holomorphic differential with respect to $b$. At $a=b=0$,
\begin{equation}
 \partial_b\left(\frac{\dd z}{2y}\right)
 =-\frac{z\,\dd z}{4y^3}
 =-\frac1{8u}\dd\left(\frac{z^2}{y}\right).
 \label{eq:bexact}
\end{equation}

Exactness in~\eqref{eq:bexact} makes the first $b$-derivative of every period vanish there. The period ratio is analytic in $bu^{-3/4}$, so its departure from the limiting value is $O(b^2u^{-3/2})$. Thus the semiflat metric differs from the flat model in the normalized period coordinates by $O(\rho^{-2})$ in the model norm. The terms linear in $b$ in the periods of $\lambda$ are independent of the base parameter and were subtracted in~\eqref{eq:UW}. The same scaling gives the corresponding bounds for all covariant derivatives of the difference between the semiflat and model metrics on the end. By Proposition~\ref{prop:comparison}, the $L^2$ metric also approaches the flat model with order-$2$ decay.

\subsubsection{Type $A_2$}

Pullback by $f(z)=-z$ sends $(a,u,\mu,c)$ to $(-a,-u,\mu,c)$, preserves the punctures, labeled eigenline bundles, parabolic weights, and residues, and commutes with the elementary modifications used to normalize the weights. At $a=0$ it therefore induces a triholomorphic isometric involution $\iota$ satisfying $u\circ\iota=-u$.

On the reference curve the leading periods acquire the factor $\xi=e^{2\pi i/3}$ along a positive half-turn in the $u$-plane. Since $Z_B=e^{\pi i/3}Z_A$, comparison of the endpoint with the starting fiber using $f$ gives the integral action
\begin{equation*}
 A\longmapsto B-A,\qquad B\longmapsto-A,
 \qquad P=\begin{pmatrix}-1&1\\-1&0\end{pmatrix}.
\end{equation*}
For any fixed $\mu,c$, the same matrix gives the action on integral homology. If $Z=(Z_A,Z_B)^t$ and $C_Z$ is its constant Taylor coefficient, the affine return $Z\mapsto PZ+b_Z$ satisfies $b_Z=(\Id-P)C_Z$, by comparison of the constant terms. Thus the translation used in~\eqref{eq:A2-period-remainders} makes this return linear. For the angle coordinates the half-turn preserves the peripheral characters and acts affinely on their torsor. Since the induced integral action on the character torus has $\det(\Id-P)=3$, choose a parallel origin fixed by this affine half-turn. Its square, the full return, fixes the same origin. In the coordinates at infinity defined by the spectral periods, the action is
\begin{equation*}
 U\longmapsto\xi U,\quad W\longmapsto\xi^{-1}W,
 \quad v\longmapsto\xi^{-1}v.
\end{equation*}
It follows that $(S,v)\mapsto(\xi S,\xi^{-1}v)$ for $S=U-\bar W$. In~\eqref{eq:semiflat-forms} the flat model is obtained by setting $W_U=0$. The coefficients of the ordered triple are smooth functions of $W_U,\overline{W_U}$ near zero, and the metric is a smooth algebraic function of a normalized hyperk\"ahler triple. The differentiated estimates following~\eqref{eq:A2-period-remainders} thus give $g_{\sfm}-g_{\mathcal C}=\Op(\rho^{-2})$ when $a=0$, and $\Op(\rho^{-1/2})$ for general $a$. Theorem~\ref{thm:comparison} and Lemma~\ref{lem:comparison-bootstrap} then give positive-order decay, with all spatial derivatives, for $g_{L^2}-g_{\mathcal C}$ in these spectral coordinates.

On the logarithmic cover of the base, the chosen lift shifts $\arg S$ by $2\pi/3$, and its square is the ALG deck transformation of angle $4\pi/3$. On the threefold cover that trivializes the torus monodromy, the angular period is $4\pi$ and the lift has order six. Since $u\circ\iota=-u$, the involution $\iota$ is free for sufficiently large $|u|$, and the quotient end has cone parameter $\beta'=\tfrac12\cdot\tfrac23=\tfrac13$. The ordered hyperk\"ahler forms descend to this quotient, with the same decay estimates for their differences from the model forms and for covariant derivatives of every order.

The exterior argument for the ALG rate in Chen--Chen~\cite[Theorem~A; arXiv version, Section~3, Theorem~4.4 and proof of Theorem~4.6]{CC19} applies to this quotient end. The smallest positive decay exponent for closed anti-self-dual forms on this quotient end is
\begin{equation*}
 \min_{n\in\Z,\,n<2\beta'}\frac{2\beta'-n}{\beta'}
 =\min_{n\le0}(2-3n)=2.
\end{equation*}
The proof first uses an exterior right inverse for the gauge operator to improve the self-dual error. In the nonconstant torus Fourier modes, subtracting the exterior inverse applied to the inhomogeneous term leaves an exponentially decaying harmonic remainder. The invariant closed anti-self-dual modes have the exponents just computed. The iteration and its vector fields are defined on exterior domains, so fixed points of the involution in the compact part do not enter. The gauge vector fields have size $\Op(\rho^{1-\epsilon})$ for some $\epsilon>0$ at the first step and improve thereafter. Their flows lift to the original twofold cover, are isotopic to the identity, and preserve the leading coefficient of $S^{1/\beta_2}$. In the lifted coordinates the $L^2$ metric of $A_2$ at $a=0$ approaches the flat model with order-$2$ decay for every admissible choice of compact periods.

\section{Isomonodromic transport on ALG ends}
\label{sec:finite}

Integrating the variation of the K\"ahler form along the filtered isomonodromic transport gives a global potential. Its leading term, $-2K_j\Rea(a u_0)$, determines the two real Chen--Viaclovsky parameters.

\subsection{Isomonodromic transport of filtered connections}

\begin{proposition}
\label{prop:finite-transport}
For every admissible choice of weights and residues in $A_j$, $j=0,1,2$, there are global diffeomorphisms $T_a:\cM_0\to\cM_a$ for all $a\in\C$, with $T_0=\Id$. They are holomorphic for the de Rham complex structure $J$ and preserve $\omega_I$ and $\omega_K$. Along every smooth parameter path $a(s)$, they satisfy
\begin{equation}
 \frac{\dd}{\dd s}T_{a(s)}^*\omega_{J,a(s)}
 =-K_j\ddc\Rea\bigl(\dot a(s)u_{a(s)}\circ T_{a(s)}\bigr).
 \label{eq:actualvariation}
\end{equation}
They are defined on the entire stable locus, including points over singular spectral curves.
\end{proposition}

\begin{proof}
For $j=0,1,2$, the germs $f_{j,c}$ at $w=0$ satisfy $f_{j,c}\circ f_{j,d}=f_{j,c+d}$. They are
\begin{equation}
 f_{0,c}(w)=\frac{w}{(1+2cw^2)^{1/2}},\qquad
 f_{1,c}(w)=\frac{w}{(1+cw^2)^{1/2}},\qquad
 f_{2,c}(w)=\frac{w}{1+cw/2}.
 \label{eq:ALG-germs}
\end{equation}
Here the roots in the denominators take the value $1$ at $w=0$. The pullbacks change the coefficient $a$ to $a+c$ and preserve the leading coefficient and every residue. For example,
\begin{align*}
 f_{1,c}^*(w^{-4}\dd w)&=\left(w^{-4}+\frac c2w^{-2}+O(1)\right)\dd w,\\
 f_{1,c}^*(w^{-2}\dd w)&=w^{-2}\dd w+O(1)\dd w,\\
 f_{1,c}^*(w^{-1}\dd w)&=w^{-1}\dd w+O(w)\dd w.
\end{align*}
The odd lift of the first germ to the double cover on which the eigenvalue one-forms become single-valued is $\zeta(1+2c\zeta^4)^{-1/4}$. For the third germ,
\begin{equation*}
 f_{2,c}^*(w^{-3}\dd w)=\left(w^{-3}+\frac c2w^{-2}\right)\dd w,
 \qquad f_{2,c}^*(w^{-2}\dd w)=w^{-2}\dd w.
\end{equation*}
These formulas also show that the leading Stokes directions are unchanged.

For $a$ in a compact subset of $\C$, the germs are defined on a common disk at infinity. Extend them to orientation-preserving diffeomorphisms $F_{j,a}$ equal to the identity outside a larger disk, choosing isotopies from the identity continuously in $a$. In type $A_2$ the larger disk avoids the tame pole. Pullback gives a filtered flat connection with the prescribed principal parts. Its $(0,1)$ component defines the holomorphic bundle and flatness makes its $(1,0)$ component a meromorphic connection. On the flat deformation complex the de Rham complex structure sends $\delta\cD$ to $i\delta\cD$. Since $F_{j,a}^*$ is complex linear on endomorphism-valued one-forms and commutes with the flat differential, its induced tangent map is complex linear. Thus the map on smooth de Rham moduli is $J$-holomorphic.

Any two extensions in the chosen class are isotopic relative to open neighborhoods of the poles on which the extensions agree. Parallel transport for the flat connection along such an isotopy gives a gauge isomorphism preserving the filtration and the labeled Stokes data. Within this class, the induced map on moduli therefore depends only on the coordinate germs. Composition and inversion of these maps on moduli are induced by $f_{j,c}\circ f_{j,d}=f_{j,c+d}$ and $f_{j,-c}$, respectively. For $A_1,A_2$, Boalch constructed isomonodromic transport for admissible deformations of the wild character varieties for the prescribed unramified irregular types~\cite[Definition~10.1 and Theorem~10.2]{Boalch14}. Here the differences between the entries of the irregular type have constant pole order, so the $a$-plane parametrizes an admissible deformation. For $A_0$, the odd lift of $f_{0,c}$ commutes with the deck involution, so pullback descends to the original curve and preserves the filtration and the labeled Stokes data.

Pullback preserves the flat subbundles with their induced filtrations and parabolic degrees. By the wild nonabelian Hodge correspondence and the smooth local gauge slices of Section~\ref{sec:families}, pullback induces global diffeomorphisms
\begin{equation}
 T_a:\cM_0\longrightarrow\cM_a,
 \qquad T_0=\Id.
 \label{eq:global-transport}
\end{equation}
In the normalization of Section~\ref{sec:variation}, their de Rham holomorphic symplectic form is the convergent pairing
\begin{equation*}
 \omega_K+i\omega_I=-i\int_X\tr(\delta\cD_1^0\wedge\delta\cD_2^0).
\end{equation*}
Pullback by the orientation-preserving diffeomorphism leaves this integral invariant, so $T_a$ preserves both real forms.

Transporting the flat connection back to the fixed marked surface replaces $j_0$ by $(F_{j,a}^{-1})^*j_0$. It agrees with $j_0$ near the poles and outside the annulus where the germ was extended. Theorem~\ref{thm:variation} and~\eqref{eq:ALG-Hamiltonians} give~\eqref{eq:actualvariation} for the pulled-back form along every smooth path $a(s)$.
\end{proof}

\subsection{Bounds for the Hitchin map and its derivatives}

For $j=0,1,2$, $\beta_j$ denotes the cone parameter of $A_j$. With a choice of basepoint $x_0\in\cM_0$, we define $\rho_a(x)=1+\operatorname{dist}_{g_a}(x,T_a(x_0))$. When the fixed weights and residues vary, choose $x_0$ locally along a smooth section of the family at $a=0$. Finitely many such choices suffice over a compact parameter set. Covariant derivatives and tensor norms in the following estimates are taken with respect to $g_a$.
\begin{lemma}
\label{lem:hitchin-map-bounds}
On compact sets of admissible weights and prescribed principal parts, the Hitchin map satisfies
\begin{equation}
 |\nabla^ku_a|\leq C_k\rho_a^{1/\beta_j-k},\qquad
 |\nabla^k\dd\dd_{J_a}^cu_a|\leq C_k\rho_a^{1/\beta_j-2-k}.
 \label{eq:Hessianbounds}
\end{equation}
The constants depend on the compact parameter set and on $k$.
\end{lemma}

\begin{proof}
On the end, the differentiated spectral expansions of Section~\ref{sec:periodend} give $|\nabla^ku_a|\le C_k\rho_a^{1/\beta_j-k}$ in the model coordinates. The metric comparison and Lemma~\ref{lem:comparison-bootstrap} transfer these estimates to $g_a$. Since $\nabla^{g_a}J_a=0$ and $\dd^c_{J_a}u_a=-\dd u_a\circ J_a$, the second bound follows from that for $\nabla^{k+2}u_a$. On a bounded closed disk of the Hitchin base the relative properness established in Section~\ref{sec:integral} gives a compact inverse image over each compact set of admissible parameters. The smooth slices of Section~\ref{sec:families} then bound the derivatives on that inverse image. In particular the complex Hessian decays with exponent $2-1/\beta_j$, equal to $4/5$, $2/3$, and $1/2$ for $j=0,1,2$.
\end{proof}

\subsection{Variation of the pulled-back metric}

Along any smooth path $a(s)$, $0\leq s\leq1$, with $a(0)=0$, we write
\begin{equation}
 \widehat g_s=T_{a(s)}^*g_{a(s)},\quad \widehat g_0=g_0,\qquad
 p_j=2-1/\beta_j.
 \label{eq:GsHs}
\end{equation}
\begin{proposition}
\label{prop:transport-metric}
Along every smooth parameter path $a(s)$ on a compact interval with $a(0)=0$, the metrics $\widehat g_s$ are uniformly equivalent to $g_0$, and their distance functions satisfy
\begin{equation}
 \rho_{a(s)}(T_{a(s)}x)\asymp\rho_0(x).
 \label{eq:radialcomparison}
\end{equation}
For every $k\ge0$, the $k$th $g_0$-covariant derivatives of $\widehat g_s-g_0$ and $\nabla^{\widehat g_s}-\nabla^{g_0}$ have $g_0$-norms bounded uniformly by $C_k\rho_0^{-p_j-k}$ and $C_k\rho_0^{-p_j-1-k}$, respectively.

\end{proposition}

\begin{proof}
The complex structure $J$ is independent of $s$. We contract~\eqref{eq:actualvariation} with $J$ and use~\eqref{eq:Hessianbounds} to obtain
\begin{equation}
 |\dot{\widehat g}_s|_{\widehat g_s}(x)\leq C\rho_{a(s)}(T_{a(s)}x)^{-p_j}\leq C.
 \label{eq:metricderivative}
\end{equation}
Gronwall's inequality gives $e^{-C}g_0\leq\widehat g_s\leq e^Cg_0$. Since $T_{a(s)}$ is an isometry for the pullback metric and carries the chosen basepoint to the target basepoint, this gives~\eqref{eq:radialcomparison}. Integration of~\eqref{eq:metricderivative} now gives $|\widehat g_s-g_0|_{g_0}\le C\rho_0^{-p_j}$.

To estimate covariant derivatives on the moduli space, we differentiate the Levi-Civita connection of the pulled-back metric. For $\mathcal A_s=\nabla^{\widehat g_s}-\nabla^{g_0}$,
\begin{align}
 2\widehat g_s(\dot{\mathcal A}_s(X,Y),Z)
 &=(\nabla_X^{\widehat g_s}\dot{\widehat g}_s)(Y,Z)+(\nabla_Y^{\widehat g_s}\dot{\widehat g}_s)(X,Z)\notag\\
 &\quad-(\nabla_Z^{\widehat g_s}\dot{\widehat g}_s)(X,Y).
 \label{eq:connectionvariation}
\end{align}
The bounds~\eqref{eq:Hessianbounds} for $g_{a(s)}$ pull back by the isometry $T_{a(s)}$ and give
\begin{equation*}
 |(\nabla^{\widehat g_s})^k\dot{\widehat g}_s|_{\widehat g_s}\leq C_k\rho_0^{-p_j-k}.
\end{equation*}
Write $D_0=\nabla^{g_0}$ and $D_s=\nabla^{\widehat g_s}$. We already have the order-zero estimate for $\widehat g_s-g_0$ and, by uniform equivalence, the same estimate for $\widehat g_s^{-1}-g_0^{-1}$. Suppose these estimates hold through $k$ derivatives and that $D_0^m\mathcal A_s=O(\rho_0^{-p_j-1-m})$ for $m<k$. In $D_0^k(D_s\dot{\widehat g}_s)$, replacing $D_0$ by $D_s$ introduces at most $k-1$ derivatives of $\mathcal A_s$. The pulled-back Hessian bounds therefore give $O(\rho_0^{-p_j-k-1})$. Apply $D_0^k$ to~\eqref{eq:connectionvariation} after raising its last index with $\widehat g_s^{-1}$. Integration from $\mathcal A_0=0$ yields that same bound for $D_0^k\mathcal A_s$. Now $D_0^{k+1}\dot{\widehat g}_s$ involves at most $k$ derivatives of $\mathcal A_s$. Integration yields $D_0^{k+1}(\widehat g_s-g_0)=O(\rho_0^{-p_j-k-1})$, and differentiation of the inverse gives the corresponding estimate for $\widehat g_s^{-1}-g_0^{-1}$. Thus
\begin{equation}
 \widehat g_s-g_0=\Op(\rho_0^{-p_j}),\qquad
 \mathcal A_s=\Op(\rho_0^{-p_j-1}).
 \label{eq:strainall}
\end{equation}
\end{proof}

\subsection{Asymptotics of the transport map at infinity}

\begin{lemma}
\label{lem:endrigidity}
Let $g_{\mathcal C}$ be the normalized flat ALG model metric, with radial coordinate $r$ on its base cone, and let $T$ be a proper diffeomorphism between two marked ALG ends. Assume that the pullback of the target radius is comparable to the source radius, that the two metrics have the form $g_{\mathcal C}+\Op(r^{-p})$ for $p>0$ in their respective model coordinates, and that $T$ is an isometry pulling back two forms of the target hyperk\"ahler triple to their corresponding source forms, with matching model limits. Then, on the model cover, $\dd T$ tends to a constant orthogonal linear map $R_\infty$ preserving the model hyperk\"ahler triple. If $p>1$, a lift to the model cover satisfies
\begin{equation*}
 \dd T-R_\infty=\Op(r^{-p}),\qquad T(x)-R_\infty x=O(1)\,.
\end{equation*}
In a continuous family starting at the identity and preserving the marked lattice of the torus fibers at infinity, $R_\infty=\Id$. In such a family, if $0<p<1$, then
\begin{equation}
 \dd T-\Id=\Op(r^{-p}),\qquad
 T-\Id=\Op(r^{1-p}),
 \label{eq:rigidity}
\end{equation}
where the second estimate is unchanged by a bounded translation along the torus fiber.
\end{lemma}

\begin{proof}
Uniform equivalence of both metrics to the model metric bounds $\dd T$. In flat Cartesian coordinates on the model cover, the transformation law for the Levi-Civita connection is
\begin{equation}
 \dd^2T=\dd T\,\Gamma_{\mathrm{source}}
          -\Gamma_{\mathrm{target}}(\dd T,\dd T)
        =O(r^{-p-1}).
 \label{eq:D2T}
\end{equation}
Comparability of the two radii gives the same order of decay for the Christoffel symbols of the target metric evaluated at the image point. Along each ray, $\dd T$ has a limit $R_\infty$, with $\dd T-R_\infty=O(r^{-p})$. At radius $r$ an angular arc of length $O(r)$ joins the two rays. By~\eqref{eq:D2T} the values of the differential where the arc meets the two rays differ by $O(r^{-p})$, so the limits along the rays agree. A path across a fundamental parallelogram of the torus has bounded length, so this estimate shows that the limit is independent of the point of the torus. Thus $R_\infty$ is constant on the model cover and intertwines the linear parts of its deck transformations.

The limits of the metric and two-form identities show that $R_\infty$ is orthogonal and preserves two model symplectic forms, hence all three model symplectic forms. For each cone angle $2\pi\beta<2\pi$ occurring here, the fundamental group of a model end is $\Lambda\rtimes_M\Z$, where $\Lambda\cong\Z^2$ is the fiber lattice and $1$ is not an eigenvalue of the monodromy $M$. The commutator subgroup $(M-\Id)\Lambda$ has finite index in $\Lambda$. Its isolator, the set of elements having some positive power in that subgroup, is exactly $\Lambda$. Thus every end diffeomorphism preserves the fiber subgroup. More generally, the fiber subgroup is also detected by loops of uniformly bounded length based arbitrarily far out: a loop of bounded length at radius $r$ stays at radius $r-O(1)$, where nonzero angular winding costs length of order $r$. The metric equivalence and the radial bounds preserve this characterization. Integrating $\dd T$ along representatives of the two fiber generators and taking $r\to\infty$ shows that $R_\infty$ restricts to the induced lattice automorphism on the torus plane. The action on the marked lattice is locally constant along a continuous family, so the family starting at the identity fixes the torus plane pointwise in the limit. By orthogonality, $R_\infty$ is block diagonal. Since it preserves the holomorphic symplectic form pairing the two planes, it also acts as the identity on the base plane.

In the continuous family, the integral of $\dd T-\Id=O(r^{-p})$ along rays gives $T-\Id=O(r^{1-p})$ for $0<p<1$. Differentiating~\eqref{eq:D2T} and using the $O'$ bounds for the two metrics gives the higher derivative estimates by induction. For $p>1$ the radial integral of $\dd T-R_\infty$ converges. Its integral across an angular arc of bounded opening is $O(r^{1-p})$, tending to zero, and across a fiber parallelogram it is $O(r^{-p})$. Thus on the chosen lift $T(x)-R_\infty x$ tends to a constant translation and $\dd T-R_\infty=\Op(r^{-p})$.
\end{proof}

Use the coordinates $(S_a,v_a)$ of Section~\ref{sec:periodend} with parallel origins continued along the parameter path. Theorem~\ref{thm:comparison}, Lemma~\ref{lem:comparison-bootstrap}, and the differentiated period expansions give $g_a-g_{\mathcal C}=\Op(\rho^{-\delta_j})$, with $\delta_j=p_j$ for $j=0,1$. For $A_2$ take $\delta_2=\min(\varepsilon_*,1/2)$, where $\varepsilon_*>0$ is the common Hitchin-to-semiflat exponent on the path and $1/2$ is the semiflat-to-model exponent from~\eqref{eq:A2-period-remainders} and~\eqref{eq:semiflat-forms}.

In the source coordinates, $\widehat g_s$ differs from the flat model by $\Op(\rho^{-\delta_j})$ by~\eqref{eq:strainall} and the corresponding bound for $g_0$. Equation~\eqref{eq:radialcomparison} compares the source and target radii. For each of the preserved forms $\omega_I,\omega_K$, the source and target have the same normalization at infinity. The torus marking varies continuously along the path. The displacement in the base coordinate therefore satisfies, by Lemma~\ref{lem:endrigidity},
\begin{equation}
 S_{a(s)}\circ T_{a(s)}=S_0+\Op(\rho_0^{1-\delta_j}).
 \label{eq:Sdisplacement}
\end{equation}

On each lifted sector the period expansions give $u_a=(S_a/\ell_j)^{1/\beta_j}+\Op(\rho_a^{1/\beta_j-\delta_j})$. In type $A_2$, equation~\eqref{eq:A2-base-inverse} gives an $\Op(\rho_a)$ remainder, sufficient since $\delta_2\le 1/2$. The derivative of the leading term has size $O(\rho^{1/\beta_j-1})$, so~\eqref{eq:Sdisplacement} implies
\begin{equation}
 u_a\circ T_a=u_0+\Op(\rho_0^{1/\beta_j-\delta_j}).
 \label{eq:Hdifference}
\end{equation}
For $A_0$ the remainder is $\Op(\rho^{2/5})$, for $A_1$ it is $\Op(\rho^{2/3})$, and for $A_2$ it is $\Op(\rho^{3/2-\delta_2})$.

\subsection{A global potential for the difference of the K\"ahler forms}

\begin{proof}[Completion of the proof of Theorem~\ref{thm:transport}]
We integrate~\eqref{eq:actualvariation} along the straight path $sa$ from zero to the given parameter. The globally defined smooth real function
\begin{equation}
 \phi_a=-2K_j\Rea\int_0^1a\,(u_{sa}\circ T_{sa})\,\dd s
 \label{eq:potential}
\end{equation}
satisfies
\begin{align}
 T_a^*\omega_{J,a}&=\omega_{J,0}+i\partial_J\bar\partial_J\phi_a,
 \label{eq:finitevariation}\\
 \phi_a&=-2K_j\Rea(a u_0)+\Op(\rho_0^{1/\beta_j-\delta_j}).
 \label{eq:potentialasymp}
\end{align}
The first equality follows from $\ddc=2i\partial_J\bar\partial_J$, and the second from~\eqref{eq:Hdifference}.

The $L^2$ metric at $a=0$ approaches the ALG model with order-$2$ decay in each of the three cases by Section~\ref{subsec:centers}. To apply the uniqueness statement in~\cite[Proposition~5.3]{CV24}, we use the normal form constructed in~\cite[Section~5, equations~(5.13)--(5.14)]{CV24}. In that chart Chen--Viaclovsky use a global $I$-holomorphic function asymptotic to the $1/\beta_j$ power of the model base coordinate. With the present choice of phase it is
\begin{equation}
 z_{\mathrm{CV}}=\ell_j^{1/\beta_j}u_0+C,
 \label{eq:CV-holomorphic-coordinate}
\end{equation}
where $C\in\C$ is constant. This holomorphic function is constant on each compact connected Hitchin fiber, so it factors through $u_0$. As a function of $u_0$, it is entire with at most linear growth and hence affine; the normalization of the spectral periods fixes its linear coefficient. Changing the additive constant leaves the complex Hessian of the potential unchanged.

We rotate the ordered triple to $(\omega_I,\omega_K,-\omega_J)$. Its third complex structure is $-J$, and $\partial_{-J}\bar\partial_{-J}=-\partial_J\bar\partial_J$ on real functions. The same potential therefore changes the third form by $i\partial_{-J}\bar\partial_{-J}\phi_a$. With
\begin{equation}
 \gamma_j=K_j\ell_j^{-1/\beta_j},\qquad
 x_{\mathrm{CV}}=\Rea z_{\mathrm{CV}},\qquad y_{\mathrm{CV}}=\Ima z_{\mathrm{CV}},
 \label{eq:gamma-definition}
\end{equation}
the potential has the asymptotic expansion
\begin{equation}
 \phi_a=2\bigl((\gamma_j\Ima a)y_{\mathrm{CV}}-(\gamma_j\Rea a)x_{\mathrm{CV}}\bigr)
       +\Op(\rho_0^{1/\beta_j-\delta_j}).
 \label{eq:CVpotentialmatch}
\end{equation}
We denote the two real parameters in~\cite[Proposition~5.3]{CV24} by $a_{\mathrm{CV}},b_{\mathrm{CV}}$, so that the leading potential is $2(a_{\mathrm{CV}}y_{\mathrm{CV}}-b_{\mathrm{CV}}x_{\mathrm{CV}})$. We encode this pair by $c_{\mathrm{CV}}=a_{\mathrm{CV}}+ib_{\mathrm{CV}}$.
In the order-$2$ chart of the reference metric, take $-iv$ as the model fiber coordinate for the rotated triple. Its lattice is $L(\Z+\tau_0\Z)$ in each case: the extra factors $e^{-\pi i/3}$ in $A_0$ and $-i$ in $A_1$ preserve the corresponding lattices. In $A_2$ the lifted change of coordinates described in Section~\ref{subsec:centers} preserves the leading base variable. Write $p_j=2-1/\beta_j$. For any $0<\epsilon<\min(2-p_j,\delta_j)$, the $\Op(\rho^{-2})$ reference errors satisfy (5.13)--(5.14) of~\cite{CV24} with $a_0=b_0=0$. The first two rotated forms and $-J$ are fixed by transport. On taking two derivatives of the remainder in~\eqref{eq:CVpotentialmatch}, the third form has the stated leading Hessian and an $\Op(\rho^{-p_j-\epsilon})$ error in the same chart. Proposition~\ref{prop:transport-metric} gives uniform equivalence to the complete metric $g_0$ and the estimate $T_a^*g_a-g_0=\Op(\rho_0^{-p_j})$. The pulled-back metric is complete and positive, and pullback preserves the hyperk\"ahler volume identity. The potential thus satisfies (5.15)--(5.16) of~\cite[Proposition~5.3]{CV24} in its required fixed coordinate system. Uniqueness therefore identifies the pulled-back metric with the member of the two-parameter family having
\begin{equation}
 c_{\mathrm{CV}}=i\gamma_j\bar a.
 \label{eq:cvcomputed}
\end{equation}
The spectral integrals defining $\ell_j$ are positive, so $a\mapsto i\gamma_j\bar a$ is a real-linear bijection of $\C$.
\end{proof}

\section{Realization of the metrics of the seven types}
\label{sec:realization}

\subsection{The homology marking and coordinates at infinity}

For a smooth gravitational instanton $(M,g,\omega_I,\omega_J,\omega_K)$ of any of the seven types in Table~\ref{tab:types}, we rescale the metric and rotate the hyperk\"ahler triple so that its period vector on the fiber is $(4\pi^2,0,0)$. For an ALG$^*$ metric, the change of coordinates in~\cite[Theorem~1.10]{CVZH23} gives order-$2$ decay at infinity and preserves $\nu,\kappa_0,L$. In the ALG case, the two-parameter family constructed in~\cite[Theorem~F]{CV24} contains a metric with order-$2$ decay and the same periods on compact two-cycles.

By~\cite[Theorems~C and~E]{CV24}, the manifolds in the Torelli theorems of~\cite[Theorems~1.5 and~1.10]{CVZ24} admit diffeomorphisms compatible with the coordinates at infinity. For the reference Hitchin moduli spaces we use the elliptic compactifications and integral lattices calculated in Section~\ref{sec:integral}. In the ALG$^*$ case, we first choose a positive real value of the leading coefficient or product of leading coefficients from Table~\ref{tab:end-parameters}, so that the reference Hitchin metric and the given metric have the same $\kappa_0$; this choice is independent of their compact periods. Use such a diffeomorphism to transport the metric to the marked reference manifold, taking the order-$2$ member in the ALG case. Every integral class of square $-2$ has a nonzero period vector by~\cite[Theorems~1.7 and~1.12]{CVZ24}. For $A_1,A_2,D_3$, the simple roots are represented by the compact components of the reducible fibers in Sections~\ref{sec:A1}, \ref{sec:A2} and~\ref{sec:D3}. The remaining component in each fiber represents the affine simple root: $\Gamma_-=F-\Gamma_+$ in $A_1$, $\Gamma_3=F-\Gamma_1-\Gamma_2$ in $A_2$, and $\Gamma_4=F-\Gamma_0-\Gamma_c-\Gamma_1$ in $D_3$. Each component is a smooth embedded sphere of self-intersection $-2$. A Dehn twist in its disk normal bundle induces the Picard--Lefschetz reflection $R_\delta(\eta)=\eta+(\eta\cdot\delta)\delta$. These twists generate the affine reflections needed to place the real periods in the closed alcoves of Section~\ref{sec:integral}. Their compact supports preserve the end chart, and each reflection acts on the complex periods in the same marking.

By Lemma~\ref{lem:marking-transport} and its applications in Section~\ref{sec:integral}, the marking of integral homology extends through each admissible family by diffeomorphisms that agree with the chosen identifications of the spectral fibrations near infinity. In the ALG$^*$ families the chosen leading coefficient or product of leading coefficients remains fixed along this transport. At the endpoints we use the same leading model, with the common phase in the $D$ cases determined in Section~\ref{sec:periodend}. At fixed $u$ and fixed parallel angles the period expansions give base changes asymptotic to the identity with derivative tending to the identity. The angle marking is constant.

Sufficiently far out, the end changes are isotopic to the identity by interpolation. The gauges improving the decay to order $2$ are flows isotopic to the identity on a smaller exterior domain. Multiplying the generating vector fields by a cutoff extends these isotopies over the manifold, equal to the specified changes far out and equal to the identity on a compact core. Composing the marking with the endpoint adjustments therefore makes the two order-$2$ charts agree, up to the model automorphism in~\cite[Theorems~C and~E]{CV24}, without changing its action on integral homology. This gives the fixed coordinate system required in the two Torelli theorems.

\subsection{Explicit parameter choices}

Applying the inverse formulas of Section~\ref{subsec:inverse-periods} to the marked compact periods gives admissible weights and residues, since every root has a nonzero period vector.

For ALG$^*$, retain the positive leading parameter already chosen from Table~\ref{tab:end-parameters}. Sections~\ref{sec:families} and~\ref{sec:integral} give smoothness on the admissible domain. Properness over bounded Hitchin disks and the semiflat lower length bound $c\int^r\sqrt{\log t}\,\dd t$ for curves reaching radius $r$ give completeness. Theorem~\ref{thm:comparison} and the order-$2$ improvement of~\cite[Theorem~1.10]{CVZH23} give the same model in the common end chart. The periods of the triples agree on the marked integral bases, so by the ALG$^*$ Torelli theorem~\cite[Theorem~1.10]{CVZ24} there is an isometry preserving the ordered triples.

\subsection{ALG metrics with order-\texorpdfstring{$2$}{2} decay and their deformations}

For each $A_j$, $j=0,1,2$, we take the weights and residues computed in Section~\ref{subsec:inverse-periods} and specialize the family to $a=0$. By Section~\ref{subsec:centers}, the $L^2$ Hitchin metric approaches the model with order-$2$ decay. The member with this decay rate in the family of the given metric has the same cone parameter, limiting torus, and fiber area, and their compact periods agree in the chosen marking. By the ALG Torelli theorem~\cite[Theorem~1.5]{CVZ24}, there is an isometry between the two metrics that preserves their ordered hyperk\"ahler triples.

Regard this isometry as a map from the Hitchin metric at $a=0$ to the order-$2$ member of the family containing the given metric, keeping the normalized Hitchin coordinates fixed. Lemma~\ref{lem:endrigidity} with $p=2$ gives a limiting model isometry $R_\infty$ preserving the ordered hyperk\"ahler triple and the torus lattice. Composing the target coordinates with $R_\infty^{-1}$ adjusts the base and torus coordinates together by this model isometry. On the lifted end the isometry now has bounded displacement and differential $\Id+\Op(\rho^{-2})$. Pullback changes a leading term of degree $1/\beta_j$ by $\Op(\rho^{1/\beta_j-1})$, so the leading coefficient and the order of growth of the potential required in~\cite[Proposition~5.3]{CV24} are preserved. The isometry also preserves the two forms and the complex structure fixed in that construction. By uniqueness in~\cite[Proposition~5.3]{CV24}, pullback by this isometry identifies the two-parameter family associated with the given metric with that of the Hitchin metric at $a=0$.

In the transported coordinates at infinity, $c_{\mathrm{CV}}$ denotes the coefficient of the leading term in the potential for the given metric. By Theorem~\ref{thm:transport}, the Hitchin metric realizing this member has parameter
\begin{equation}
 a=\frac{i\overline{c_{\mathrm{CV}}}}{\gamma_j}.
 \label{eq:CV-inverse}
\end{equation}
Since $\gamma_j>0$, the coefficient of the leading term of the potential takes every complex value. By~\eqref{eq:finitevariation} the compact periods remain constant throughout this two-parameter family, which completes the proof of Theorem~\ref{thm:realization}.
\appendix
\section{Further applications of the variation formula}
\label{app:variation}

We use the residue formula~\eqref{eq:generalresiduevariation} to compute Hamiltonians for changes of the leading ALG$^*$ coefficients and for a rank-two family with a pole of order five.

\subsection{Leading-coefficient variations in the \texorpdfstring{ALG$^*$}{ALG-star} families}

Applying Corollary~\ref{cor:residuevariation} to paths of nonzero leading coefficients in the four ALG$^*$ families, we use the infinitesimal generators $v_0,v_\infty$ of the coordinate changes at $z=0$ and $z=\infty$, with $w=z^{-1}$ at infinity:
\begin{equation*}
 \begin{array}{c|c|c}
 \text{Type}&v_0&v_\infty\\\hline
 D_0&-(\dot p_0/p_0)z\partial_z&-(\dot p_\infty/p_\infty)w\partial_w\\
 D_1&-(\dot p/p)z\partial_z&-(\dot b/b)w\partial_w\\
 D_2&-(\dot p_0/p_0)z\partial_z&-(\dot p_\infty/p_\infty)w\partial_w\\
 D_3&0&-(\dot p/p)w\partial_w
 \end{array}
\end{equation*}
Dots denote path derivatives. In $D_0,D_1,D_2$, the coefficients of $z^{-2}\dd z^2$ at zero and $w^{-2}\dd w^2$ at infinity in $2\lambda^2$ are both $2u$, so the sum of the two contraction residues is $-2u\,\partial_s\log(p_0p_\infty)$ in $D_0,D_2$ and $-2u\,\partial_s\log(pb)$ in $D_1$.

In $D_3$, the $w^{-2}$ coefficient is $2(u+c_0^2+c_1^2)$, so the residue is $-2(u+c_0^2+c_1^2)(\dot p/p)$. The terms involving the residues are constant on the moduli space, so their $\ddc$ vanishes. Thus, with $h_{\mathrm{lead}}=p_0p_\infty,pb,p_0p_\infty,p$ in these four cases, respectively,
\begin{equation}
 \dot\omega_J=-8\pi\ddc\Rea\bigl(u\,(\dot h_{\mathrm{lead}}/h_{\mathrm{lead}})\bigr).
 \label{eq:logarithmic-Hamiltonian}
\end{equation}
In~\eqref{eq:logarithmic-Hamiltonian}, the forms and the Hitchin base function $u$ are pulled back to the initial moduli space by the filtration-preserving isomonodromic transport along the chosen path. The nonzero coefficients lie on the chosen cover of their parameter space, with residues and weights fixed along the path.

\subsection{A rank-two family with a pole of order five}

We apply the variation theorem to compute the Hamiltonians on the wild moduli space of rank-two filtered Higgs bundles with a pole of order five at infinity. For any $r\in\C$ and weights $(\alpha,1-\alpha)$ with $0<\alpha<1$, we prescribe the parabolic determinant line bundle and use an underlying holomorphic bundle of degree $-1$. The spectral curve is
\begin{equation}
 y^2=z^6+b_4 z^4+b_3 z^3+\left(\frac{b_4^2}{4}+2r\right)z^2+u_1z+u_0,
 \qquad \lambda=y\,\dd z.
 \label{eq:orderfivecurve}
\end{equation}
Here $(b_4,b_3)$ determine the terms of orders three and two in the principal part of the eigenvalue one-form at infinity and $(u_1,u_0)$ are affine coordinates on the Hitchin base. The coefficient of $z^2$ makes the residues of the eigenvalue one-forms $\pm y\,\dd z$ at infinity equal to $\mp r$, independently of $(b_4,b_3,u_1,u_0)$.

\begin{corollary}
\label{cor:orderfive}
On the smooth locus of the moduli space of stable filtered Higgs bundles with spectral equation~\eqref{eq:orderfivecurve}, the relation between the two Higgs fields is given by~\eqref{eq:phase}. Filtered isomonodromic transport gives
\begin{equation}
 \dot\omega_J
 =-2\pi\ddc\Rea(\dot b_4\,u_0)
  -4\pi\ddc\Rea(\dot b_3\,u_1).
 \label{eq:orderfivevariation}
\end{equation}
The moduli space has complex dimension four, and hence real dimension eight.
\end{corollary}

\begin{proof}
At $w=z^{-1}=0$, the binomial expansion gives
\begin{equation*}
 \lambda_+=-\left(w^{-5}+\frac{b_4}{2}w^{-3}
                +\frac{b_3}{2}w^{-2}+r w^{-1}+O(1)\right)\dd w.
\end{equation*}
Indeed, the coefficient of $w^4$ in the square root is $\frac12(b_4^2/4+2r)-b_4^2/8=r$. The vector field
\begin{equation*}
 v_\infty=-\frac{\dot b_4}{4}w^3\partial_w
      -\frac{\dot b_3}{2}w^4\partial_w
\end{equation*}
gives the principal part $-(\dot b_4/2)w^{-3}\dd w-(\dot b_3/2)w^{-2}\dd w$ for the Lie derivative of $\lambda_+$. Taking the Lie derivative along $w^3\partial_w$ gives no residue contribution from $w^{-3}\dd w$. Thus the displayed vector field induces both parameter variations and preserves the residue.

For $q_{\mathrm{sp}}=2y^2\dd z^2$, contraction with the displayed vector field gives
\begin{equation*}
 \Res_{w=0}(\iota_{v_\infty}q_{\mathrm{sp}})
 =-\frac{\dot b_4}{2}u_0-\dot b_3\,u_1.
\end{equation*}
The two Higgs fields differ by a factor of the imaginary unit, so $\tr(\Phi^2)=-q_{\mathrm{sp}}$. By Corollary~\ref{cor:residuevariation}, the variation is therefore given by~\eqref{eq:orderfivevariation}.

A generic spectral curve has six simple branch points and two unramified points over infinity, so by Riemann--Hurwitz we have $2g-2=-4+6=2$ and $g=2$. For a spectral line bundle $\mathcal L$, by Riemann--Roch and $\deg E=-1$ we have $\chi(\mathcal L)=\chi(E)=1$ and hence $\deg\mathcal L=2$. Since $\operatorname{Pic}^2(\PP^1)$ is a point, the norm imposes no further condition on the two-dimensional spectral Picard variety. The base also has dimension two, giving a moduli space of complex dimension four. The metric variation formula holds throughout the smooth stable locus.
\end{proof}

\raggedbottom

\end{document}